# ADAPTIVE POISSON DISORDER PROBLEM

By Erhan Bayraktar,[1] Savas Dayanik[2] and Ioannis Karatzas[3]

*University of Michigan, Princeton University and Columbia University*

We study the quickest detection problem of a sudden change in the arrival rate of a Poisson process from a known value to an *unknown and unobservable* value at an *unknown and unobservable* disorder time. Our objective is to design an alarm time which is adapted to the history of the arrival process and detects the disorder time as soon as possible.

In previous solvable versions of the Poisson disorder problem, the arrival rate after the disorder has been assumed a known constant. In reality, however, we may at most have some prior information about the likely values of the new arrival rate before the disorder actually happens, and insufficient estimates of the new rate after the disorder happens. Consequently, we assume in this paper that the new arrival rate after the disorder is a random variable.

The detection problem is shown to admit a finite-dimensional Markovian sufficient statistic, if the new rate has a discrete distribution with finitely many atoms. Furthermore, the detection problem is cast as a discounted optimal stopping problem with running cost for a finite-dimensional piecewise-deterministic Markov process.

This optimal stopping problem is studied in detail in the special case where the new arrival rate has Bernoulli distribution. This is a nontrivial optimal stopping problem for a two-dimensional piecewise-deterministic Markov process driven by the same point process. Using a suitable single-jump operator, we solve it fully, describe the analytic properties of the value function and the stopping region, and present methods for their numerical calculation. We provide a concrete example where the value function does not satisfy the smooth-fit principle on a proper subset of the connected, continuously differentiable optimal stopping boundary, whereas it does on the complement of this set.

Received February 2005; revised December 2005.
[1]Supported in part by U.S. Army Pantheon Project and NSF Grant DMS-06-04491.
[2]Supported by NSF Grant DMI-04-23327 and AFOSR Grant FA 9550-06-1-0496.
[3]Supported by NSF Grant DMS-06-01774.

**1. Introduction and synopsis.** Suppose that arrivals of certain events constitute a Poisson process $N = \{N_t; t \geq 0\}$ with a known rate $\mu > 0$. At some time $\theta$, the arrival rate suddenly changes from $\mu$ to $\Lambda$. Both the *disorder time* $\theta$ and the *post-disorder arrival rate* $\Lambda$ of the Poisson process are unknown and unobservable quantities. Our problem is to find an alarm time $\tau$ which depends only on the past and the present observations of the process $N$, and detects the disorder time $\theta$ as soon as possible.

More precisely, we shall assume that $\theta$ and $\Lambda$ are random variables on some probability space $(\Omega, \mathcal{H}, \mathbb{P})$, on which the process $N$ is also defined; the variables $\theta$, $\Lambda$ are independent of each other and of the process $N$. An alarm time is a stopping time $\tau$ of the history of the process $N$. We shall try to choose such a stopping time so as to minimize the Bayes risk

$$(1.1) \qquad \mathbb{P}\{\tau < \theta\} + c\mathbb{E}(\tau - \theta)^+,$$

namely, the sum of the frequency $\mathbb{P}\{\tau < \theta\}$ of the false alarms and the expected cost $c\mathbb{E}(\tau - \theta)^+$ of the detection delay.

We shall assume that the post-disorder arrival rate $\Lambda$ has some general prior distribution $\nu(\cdot)$. Similarly, the disorder time $\theta$ will be assumed to have an exponential distribution of the form

$$(1.2) \qquad \mathbb{P}\{\theta = 0\} = \pi \quad \text{and} \quad \mathbb{P}\{\theta > t | \theta > 0\} = e^{-\lambda t}, \qquad t \geq 0,$$

for some $\pi \in [0, 1)$ and $\lambda > 0$. The Poisson disorder problem with a *known* post-disorder rate (namely, $\Lambda$ equals a known constant with probability 1) was studied first by Galchuk and Rozovskii [10] and was solved completely by Peskir and Shiryaev [15]. In the meantime, Davis [8] noticed that several forms of Bayes risks, including (1.1), admit similar solutions. He called this class of problems *standard Poisson disorder problems*, and found a partial solution. Recently, Bayraktar and Dayanik [1] solved the Poisson disorder problem when the detection delay is penalized exponentially. Bayraktar, Dayanik and Karatzas [3] showed that the exponential detection delay penalty in fact leads to another variant of standard Poisson disorder problems if the "standards" suggested by Davis are restated under a suitable reference probability measure. It was also shown [3] that use of a suitable reference probability measure reduces the dimension of the Markovian sufficient statistic for the detection problem, and the solution of the standard Poisson disorder problem was described fully.

We believe that *unknown and unobservable* post-disorder arrival rate $\Lambda$ captures quite well real-life applications of change-point detection theory. Before the onset of the new regime, past experience may help us at most to fit an a priori distribution $\nu(\cdot)$ on the likely values of the new arrival rate of $N$ after the disorder. Even after the disorder happens, we may not have enough observations to get a reliable statistical estimate of the post-disorder rate. Indeed, since a good alarm is expected to sound as soon as the disorder happens, we may have very few observations of $N$ sampled from the new regime since the disorder. The quickest detection of an unknown and unobservable shift in the drift of a Wiener process has been tackled by Beibel [4] and Beibel and Lerche [5]. However, we are unaware of any work pertaining to Poisson processes.

Let us highlight our approach to the problem and our main results. We show that the most general such detection problem is equivalent, under a reference probability measure, to a discounted optimal stopping problem with a running cost for an infinite-dimensional Markovian sufficient statistic. However, the dimension becomes finite as soon as the prior probability distribution $\nu(\cdot)$ of the post-disorder arrival rate $\Lambda$ charges only a finite number of atoms. This class of problems is of considerable interest since, in many applications, we have typically an empirical distribution of the post-disorder arrival rate, constructed either from finite past data or from expert opinions on the most significant likely values.

We then study in detail the case where the new arrival rate after the disorder is expected either to increase or to decrease by the same amount. The detection problem turns in this case into an *optimal stopping problem for a two-dimensional piecewise-deterministic Markov process*, driven by the same point process. We solve this optimal stopping problem fully by describing $\varepsilon$-optimal and optimal stopping times and identifying explicitly the nontrivial shape of the optimal continuation region.

The common approach to an optimal stopping problem for a continuous-time Markov process is to reformulate it as a free-boundary problem in terms of the infinitesimal generator of the process. The free-boundary problems sometimes turn out to be quite hard, even in one dimension; see, for example, [1, 10, 15]. Here, the infinitesimal operator gets complicated further, and becomes a singular partial differential-delay operator. Moreover, it is a nontrivial task, even in two dimensions, to guess the location, shape and smoothness of the free-boundary separating the continuation and stopping regions, as well as the behavior of the value function along the boundary.

Instead, we follow a direct approach and work with integral operators rather than differential operators. As in [9] and [11] we use a suitable single-jump operator to strip the jumps off the original two-dimensional piecewise-deterministic Markov process and turn the original optimal stopping problem into a sequence of optimal stopping problems for a deterministic process with continuous paths. Using direct arguments, we are able to infer from the properties of the single-jump operator the location and shape of the optimal continuation region, as well as the smoothness of the switching boundary and the value function.

The single-jump operator also suggests a straightforward numerical method for calculating the value function and the optimal continuation region. The deterministic process obtained after removing the jumps from the original

Markov process has two fundamentally different types of behavior. We tailor the naive numerical method to each case, by exploiting the behavior of the paths.

We also raise the question when the value function should be a classical solution of the relevant free-boundary problem. For a large range of configurations of parameters, both the value function and the boundary of the continuation region turn out to be continuously differentiable, and one may also choose to use finite-difference methods for differential-difference equations to solve the problem numerically. For a few other cases, we cannot qualify completely the degree of smoothness of the value function. Viscosity approaches or some other techniques of nonsmooth analysis are very likely to fill the gap, but we do not pursue this direction here. We report one concrete example on "partial" failure of the *smooth-fit principle*: in certain cases, the value function is continuously differentiable everywhere on the state space except on a proper subset of the connected and continuously differentiable optimal stopping boundary.

This work is divided naturally in two parts. In Sections 2–6, we describe the problem, formulate a model, and develop an important approximation. In Sections 7–11, we use that approximation to develop the solution and study its properties. Long proofs are presented in the Appendix.

**2. Problem description.** Let $N = \{N_t; t \geq 0\}$ be a homogeneous Poisson process with some rate $\mu > 0$ on a fixed probability space $(\Omega, \mathcal{H}, \mathbb{P}_0)$, which also supports two random variables $\theta$ and $\Lambda$ independent of each other and of the process $N$. We shall denote by $\nu(\cdot)$ the distribution of the random variable $\Lambda$, assume that

$$(2.1) \quad m^{(k)} \triangleq \int_{\mathbb{R}} (v - \mu)^k \nu(dv), \qquad k \in \mathbb{N}_0, \qquad \text{are well defined and finite,}$$

and that

$$(2.2) \quad \mathbb{P}_0\{\theta = 0\} = \pi \quad \text{and} \quad \mathbb{P}_0\{\theta > t\} = (1 - \pi)e^{-\lambda t}, \qquad t \geq 0,$$

hold for some constants $\lambda > 0$ and $\pi \in [0, 1)$.

Let us denote by $\mathbb{F} = \{\mathcal{F}_t\}_{t \geq 0}$ the right-continuous enlargement with $\mathbb{P}_0$-null sets of the natural filtration $\sigma(N_s; 0 \leq s \leq t)$ of $N$. We also define a larger filtration $\mathbb{G} = \{\mathcal{G}_t\}_{t \geq 0}$ by setting $\mathcal{G}_t \triangleq \mathcal{F}_t \vee \sigma\{\theta, \Lambda\}$, $t \geq 0$. The $\mathbb{G}$-adapted, right-continuous (hence, $\mathbb{G}$-progressively measurable) process

$$(2.3) \quad h(t) \triangleq \mu \mathbf{1}_{\{t < \theta\}} + \Lambda \mathbf{1}_{\{t \geq \theta\}}, \qquad t \geq 0,$$

induces the $(\mathbb{P}_0, \mathbb{G})$-martingale (see [6], pages 165–168)

$$(2.4) \quad Z_t \triangleq \exp\left\{\int_0^t \log\left(\frac{h(s-)}{\mu}\right) dN_s - \int_0^t (h(s) - \mu)\, ds\right\}, \qquad t \geq 0.$$

This martingale defines a new probability measure $\mathbb{P}$ on every $(\Omega, \mathcal{G}_t)$ by

$$(2.5) \qquad \left.\frac{d\mathbb{P}}{d\mathbb{P}_0}\right|_{\mathcal{G}_t} = Z_t, \qquad t \geq 0.$$

Since $\mathbb{P}$ and $\mathbb{P}_0$ coincide on $\mathcal{G}_0 = \sigma\{\theta, \Lambda\}$, the random variables $\theta$ and $\Lambda$ are independent and have the same distributions under both $\mathbb{P}$ and $\mathbb{P}_0$.

Under the new probability measure $\mathbb{P}$ the counting process $N$ has $\mathbb{G}$-progressively measurable intensity given by $h(\cdot)$ of (2.3), namely $N_t - \int_0^t h(s)\, ds$, $t \geq 0$, is a $(\mathbb{P}, \mathbb{G})$-martingale. In other words, the $\mathbb{G}$-adapted process $N$ is a Poisson process whose rate changes at time $\theta$ from $\mu$ to $\Lambda$.

In the Poisson disorder problem, only the process $N$ is observable, and our objective is to detect the disorder time $\theta$ as quickly as possible. More precisely, we want to find an $\mathbb{F}$-stopping time $\tau$ that minimizes the *Bayes risk*

$$(2.6) \qquad R_\tau(\pi) \triangleq \mathbb{P}\{\tau < \theta\} + c\mathbb{E}(\tau - \theta)^+,$$

where $c > 0$ is a constant, and the expectation $\mathbb{E}$ is taken under the probability measure $\mathbb{P}$. Hence, we are interested in an alarm time $\tau$ which is adapted to the history of the process $N$, and minimizes the trade-off between the frequency of false alarms $\mathbb{P}\{\tau < \theta\}$ and the expected time of delay $\mathbb{E}(\tau - \theta)^+$ between the alarm time and the unobservable disorder time.

In the next section we shall formulate the quickest detection problem as a problem of optimal stopping for a suitable Markov process.

**3. Sufficient statistics for the adaptive Poisson disorder problem.** Let $\mathcal{S}$ be the collection of all $\mathbb{F}$-stopping times, and introduce the $\mathbb{F}$-adapted processes

$$(3.1) \qquad \begin{aligned} \Pi_t &\triangleq \mathbb{P}\{\theta \leq t | \mathcal{F}_t\} \quad \text{and} \\ \Phi_t^{(k)} &\triangleq \frac{\mathbb{E}[(\Lambda - \mu)^k \mathbf{1}_{\{\theta \leq t\}} | \mathcal{F}_t]}{1 - \Pi_t}, \qquad k \in \mathbb{N}_0, t \geq 0. \end{aligned}$$

Since $\Lambda$ has the same distribution $\nu(\cdot)$ under $\mathbb{P}$ and $\mathbb{P}_0$, each $\Phi^{(k)}$, $k \in \mathbb{N}_0$ is well defined by (2.1). The process $\Pi = \{\Pi_t, t \geq 0\}$ tracks the likelihood that a change in the intensity of $N$ has already occurred, given past and present observations of the process. Each $\Phi^{(k)} = \{\Phi_t^{(k)}, t \geq 0\}$, $k \in \mathbb{N}$, may be regarded as a (weighted) *odds-ratio process*.

Our first lemma below shows that the minimum Bayes risk can be found by solving a discounted optimal stopping problem, with discount rate $\lambda$ and running cost function $f(x) = x - \lambda/c$ for the $\mathbb{F}$-adapted process $\Phi^{(0)}$. By Lemma 3.2, the observation process $X$ and the sufficient statistic $\{\Phi^{(k)}\}_{k \geq 0}$ jump exactly at the same times and evolve deterministically between jumps;

therefore, their natural filtrations and the collection of their stopping times are the same.

The calculations are considerably easier when the process $\Phi^{(0)}$ has the Markov property. Unfortunately, this is not true in general. However, the explicit dynamics of $\Phi^{(0)}$ in Lemma 3.2 reveal that the infinite-dimensional sequence $\{\Phi^{(k)}\}_{k \in \mathbb{N}_0}$ of the processes in (3.1) is always a Markovian sufficient statistic for the quickest detection problem. The same result also suggests sufficient conditions for the existence of a *finite-dimensional* Markovian sufficient statistic, a case amenable to concrete analysis.

LEMMA 3.1. *The Bayes risk in* (2.6) *equals*

$$(3.2) \quad R_\tau(\pi) = 1 - \pi + c(1-\pi)\mathbb{E}_0\left[\int_0^\tau e^{-\lambda t}\left(\Phi_t^{(0)} - \frac{\lambda}{c}\right) dt\right], \qquad \tau \in \mathcal{S},$$

*where the expectation* $\mathbb{E}_0$ *is taken under the (reference) probability measure* $\mathbb{P}_0$.

The proof is very similar to that of Proposition 2.1 in [3]. Note that every $Z_t$ in (2.4) can be written as

$$(3.3) \qquad Z_t = \mathbf{1}_{\{t < \theta\}} + \frac{L_t}{L_\theta}\mathbf{1}_{\{t \geq \theta\}}$$

in terms of the *likelihood ratio process*

$$(3.4) \qquad L_t \triangleq \left(\frac{\Lambda}{\mu}\right)^{N_t} e^{-(\Lambda-\mu)t}, \qquad t \geq 0.$$

Then the generalized Bayes theorem (see, e.g., [13], Section 7.9) and (3.3) imply

$$(3.5) \qquad 1 - \Pi_t = \frac{\mathbb{E}_0[Z_t \mathbf{1}_{\{\theta > t\}}|\mathcal{F}_t]}{\mathbb{E}_0[Z_t|\mathcal{F}_t]} = \frac{\mathbb{P}_0\{\theta > t|\mathcal{F}_t\}}{\mathbb{E}_0[Z_t|\mathcal{F}_t]} = \frac{(1-\pi)e^{-\lambda t}}{\mathbb{E}_0[Z_t|\mathcal{F}_t]},$$

since $\theta$ is independent of the process $N$ under $\mathbb{P}_0$ and has the distribution (2.2).

LEMMA 3.2. *Let* $m^{(k)}$, $k \in \mathbb{N}_0$, *be defined as in* (2.1). *Then every* $\Phi^{(k)}$, $k \in \mathbb{N}_0$, *in* (3.1) *satisfies the equation*

$$(3.6) \quad \begin{aligned} d\Phi_t^{(k)} &= \lambda(m^{(k)} + \Phi_t^{(k)})\, dt + \frac{1}{\mu}\Phi_{t-}^{(k+1)}(dN_t - \mu\, dt), \qquad t > 0, \\ \Phi_0^{(k)} &= \frac{\pi}{1-\pi}m^{(k)}. \end{aligned}$$

PROOF. For every $k \in \mathbb{N}_0$, let us introduce the function

$$(3.7) \quad F^{(k)}(t,x) \triangleq \int \left(\frac{v}{\mu}\right)^x (v-\mu)^k e^{-(v-\mu)t} \nu(dv), \qquad t \in \mathbb{R}_+, x \in \mathbb{R}.$$

The generalized Bayes theorem, (3.5), and the independence of the random variables $\theta$, $\Lambda$ and the process $N$ under $\mathbb{P}_0$ imply that we have

$$\begin{aligned}
\Phi_t^{(k)} &= \frac{\mathbb{E}_0[(\Lambda - \mu)^k Z_t \mathbf{1}_{\{\theta \leq t\}} | \mathcal{F}_t]}{(1 - \Pi_t)\mathbb{E}_0[Z_t | \mathcal{F}_t]} \\
(3.8) \qquad &= \frac{\pi e^{\lambda t}}{1-\pi} F^{(k)}(t, N_t) + \lambda \int_0^t e^{\lambda(t-s)} F^{(k)}(t-s, N_t - N_s) \, ds \\
&= U_t^{(k)} + V_t^{(k)}
\end{aligned}$$

for every $k \in \geq 0$ and $t \in \mathbb{R}_+$, where we have set

$$(3.9) \quad \begin{aligned} U_t^{(k)} &\triangleq \frac{\pi e^{\lambda t}}{1-\pi} F^{(k)}(t, N_t) \quad \text{and} \\ V_t^{(k)} &\triangleq \lambda \int_0^t e^{\lambda(t-s)} F^{(k)}(t-s, N_t - N_s) \, ds. \end{aligned}$$

Every $F^{(k)}(\cdot, \cdot)$, $k \in \mathbb{N}_0$, in (3.7) is continuously differentiable, and

$$(3.10) \quad \frac{\partial}{\partial t} F^{(k)}(t,x) = -F^{(k+1)}(t,x), \qquad t > 0, x \in \mathbb{R}, k \in \mathbb{N}_0.$$

The change of variable formula for jump processes gives

$$\begin{aligned}
F^{(k)}(t, N_t) &= F^{(k)}(0,0) + \int_0^t \frac{\partial F^{(k)}}{\partial t}(s, N_s) \, ds + \int_0^t \frac{\partial F^{(k)}}{\partial x}(s, N_{s-}) \, dN_s \\
&\quad + \sum_{0 < s \leq t} \left[ F^{(k)}(s, N_s) - F^{(k)}(s, N_{s-}) - \frac{\partial F^{(k)}}{\partial x}(s, N_{s-})\Delta N_s \right] \\
(3.11) \qquad &= m^{(k)} - \int_0^t F^{(k+1)}(s, N_s) \, ds \\
&\quad + \sum_{0 < s \leq t} [F^{(k)}(s, N_s) - F^{(k)}(s, N_{s-})],
\end{aligned}$$

where $\Delta N_s \triangleq N_s - N_{s-} \in \{0,1\}$ for $s > 0$, and the last equality follows from (3.10),

$$F^{(k)}(0,0) = m^{(k)} \quad \text{and} \quad \int_0^t \frac{\partial F^{(k)}}{\partial x}(s, N_{s-}) \, dN_s = \sum_{0 < s \leq t} \frac{\partial F^{(k)}}{\partial x}(s, N_{s-})\Delta N_s$$

for every integer $k \geq 0$. However, $F^{(k)}(s, N_s) - F^{(k)}(s, N_{s-})$ is equal to

$$\int \left(\frac{v}{\mu}\right)^{N_{s-}+\Delta N_s} (v-\mu)^k e^{-(v-\mu)s} \nu(dv) - \int \left(\frac{v}{\mu}\right)^{N_{s-}} (v-\mu)^k e^{-(v-\mu)s} \nu(dv)$$

$$= \frac{\Delta N_s}{\mu} \int \left(\frac{v}{\mu}\right)^{N_{s-}} (v-\mu)^{k+1} e^{-(v-\mu)t} \nu(dv) = \frac{1}{\mu} F^{(k+1)}(s, N_{s-}) \Delta N_s,$$

since $[(v/\mu)^{\Delta N_s} - 1] = (\Delta N_s/\mu)(v-\mu)$. This equation and (3.11) imply

$$F^{(k)}(t, N_t) = m^{(k)} - \int_0^t F^{(k+1)}(s, N_s)\, ds + \sum_{0 < s \leq t} \frac{1}{\mu} F^{(k+1)}(s, N_{s-}) \Delta N_s$$

(3.12)

$$= m^{(k)} + \int_0^t \frac{1}{\mu} F^{(k+1)}(s, N_{s-})(dN_s - \mu\, ds), \qquad t \in \mathbb{R}_+, k \in \mathbb{N}_0.$$

This identity will help us derive the dynamics of $U^{(k)}$ and $V^{(k)}$ in (3.9). Note that

$$d\left(\frac{1-\pi}{\pi} U_t^{(k)}\right) = d(e^{\lambda t} F^{(k)}(t, N_t)) = e^{\lambda t} F^{(k)}(t, N_t) \lambda\, dt + e^{\lambda t} dF^{(k)}(t, N_t)$$

$$= \lambda \frac{1-\pi}{\pi} U_t^{(k)}\, dt + \frac{e^{\lambda t}}{\mu} F^{(k+1)}(t, N_{t-})(dN_t - \mu\, dt).$$

Therefore,

(3.13)

$$dU_t^{(k)} = \lambda U_t^{(k)} + \frac{1}{\mu} U_t^{(k+1)}(dN_t - \mu\, dt), \qquad t > 0,$$

$$U_0^{(k)} = \frac{\pi}{1-\pi} m^{(k)}.$$

The derivation of the dynamics of $V^{(k)}$ is trickier. For every fixed $s \in [0, t)$, let us define $N_u^{(s)} \triangleq N_{s+u} - N_s$, $0 \leq u \leq t - s$. This is also a Poisson process under $\mathbb{P}_0$. As in (3.12),

$$F^{(k)}(t-s, N_{t-s}^{(s)}) = m^{(k)} + \int_0^{t-s} \frac{1}{\mu} F^{(k+1)}(u, N_{u-}^{(s)})(dN_u^{(s)} - \mu\, du).$$

Changing the variable of integration and substituting $N_\bullet^{(s)} = N_{s+\bullet} - N_s$ into this equality gives

$$F^{(k)}(t-s, N_t - N_s) = m^{(k)} + \frac{1}{\mu} \int_s^t F^{(k+1)}(v-s, N_{v-} - N_s)(dN_v - \mu\, dv).$$

Let us plug this identity into $V_t^{(k)}$ in (3.9), multiply both sides by $e^{-\lambda t}$ and change the order of integration. Then

$$e^{-\lambda t} V_t^{(k)} = \int_0^t \lambda e^{-\lambda s} \left(m^{(k)} + \frac{1}{\mu} \int_0^t F^{(k+1)}(v-s, N_{v-} - N_s)(dN_v - \mu\, dv)\right) ds$$

$$= m^{(k)} \int_0^t \lambda e^{-\lambda s} \, ds$$

$$+ \frac{\lambda}{\mu} \int_0^t \left( \int_0^v e^{-\lambda s} F^{(k+1)}(v-s, N_v - N_s) \, ds \right) (dN_v - \mu \, dv)$$

$$= m^{(k)} \int_0^t \lambda e^{-\lambda s} \, ds + \frac{1}{\mu} \int_0^t e^{-\lambda v} V_v^{(k+1)} (dN_v - \mu \, dv).$$

Differentiating both sides and rearranging terms, we obtain

(3.14)
$$dV_t^{(k)} = \lambda(m^{(k)} + V_t^{(k)}) \, dt + \frac{1}{\mu} V_t^{(k+1)}(dN_t - \mu \, dt), \qquad t > 0,$$

$$V_0^{(k)} = 0.$$

Adding (3.13) and (3.14) as in (3.8) gives the dynamics (3.6) of the process $\Phi^{(k)}$. $\square$

Lemma 3.2 shows that the process $\Phi^{(0)}$ does not have the Markov property in general. This is because, as (3.6) shows, $\Phi^{(0)}$ depends on $\Phi^{(1)}$, then $\Phi^{(1)}$ depends on $\Phi^{(2)}$, and so on ad infinitum. However, a finite-dimensional Markovian sufficient statistic emerges if the system of stochastic differential equations in (3.6) is *closeable*, namely, if the process $\Phi^{(k)}$ can be expressed in terms of the processes $\Phi^{(0)}, \ldots, \Phi^{(k-1)}$, for some $k \in \mathbb{N}_0$. Our next corollary shows that this is true if $\Lambda$ takes finitely many distinct values.

COROLLARY 3.3. *Suppose that* $\nu(\{\lambda_1, \ldots, \lambda_k\}) = 1$ *for some positive numbers* $\lambda_1, \ldots, \lambda_k$. *Consider the polynomial*

$$p(v) \triangleq \prod_{i=1}^k (v - \lambda_i + \mu) \equiv v^k + \sum_{i=0}^{k-1} c_i v^i, \qquad v \in \mathbb{R},$$

*for suitable real numbers* $c_0, \ldots, c_{k-1}$. *Then* $\{\Phi^{(0)}, \Phi^{(1)}, \ldots, \Phi^{(k-1)}\}$ *is a $k$-dimensional sufficient Markov statistic, with* $\Phi^{(k)} = -\sum_{i=0}^{k-1} c_i \Phi^{(i)}$.

PROOF. Under the hypothesis, the random variable $p(\Lambda - \mu) = (\Lambda - \mu)^k + \sum_{i=0}^{k-1} c_i (\Lambda - \mu)^i$ is equal to zero almost surely. Therefore, (3.1) implies

$$\Phi_t^{(k)} + \sum_{i=0}^{k-1} c_i \Phi_t^{(i)} = \frac{\mathbb{E}[p(\Lambda - \mu) \mathbf{1}_{\{\theta \le t\}} | \mathcal{F}_t]}{1 - \Pi_t} = 0, \qquad \mathbb{P}\text{-a.s., for every } t \ge 0.$$

The process on the left-hand side has right-continuous sample paths, by (3.6). Therefore, $\Phi_t^{(k)} + \sum_{i=0}^{k-1} c_i \Phi_t^{(i)} = 0$ for all $t \in \mathbb{R}_+$ almost surely, that is, the process $\Phi^{(k)}$ is a linear combination of the processes $\Phi^{(0)}, \ldots, \Phi^{(k-1)}$. $\square$

In the remainder of the paper we shall study the case where the arrival rate of the observations after the disorder has a Bernoulli prior distribution.

**4. Poisson disorder problem with a Bernoulli post-disorder arrival rate.** We shall assume henceforth $\mu > 1$ and that the random variable $\Lambda$ has Bernoulli distribution

$$\nu(\{\mu - 1, \mu + 1\}) = 1. \tag{4.1}$$

Namely, the rate of the Poisson process $N$ is expected to increase or decrease by one unit after the disorder. Corollary 3.3 implies that $\Phi^{(2)} = \Phi^{(0)}$, and the sufficient statistic $(\Phi^{(0)}, \Phi^{(1)})$ is a Markov process. According to Lemma 3.2, the pair satisfies

$$d\Phi_t^{(0)} = \lambda(1 + \Phi_t^{(0)})\,dt + \frac{1}{\mu}\Phi_{t-}^{(1)}(dN_t - \mu\,dt), \qquad \Phi_0^{(0)} = \frac{\pi}{1-\pi}, \tag{4.2}$$

$$d\Phi_t^{(1)} = \lambda(m + \Phi_t^{(1)})\,dt + \frac{1}{\mu}\Phi_{t-}^{(0)}(dN_t - \mu\,dt), \qquad \Phi_0^{(1)} = \frac{\pi}{1-\pi}m, \tag{4.3}$$

where, as in (2.1), we set

$$m \equiv m^{(1)} = \mathbb{E}_0[\Lambda - \mu] = \mathbb{P}\{\Lambda = \mu + 1\} - \mathbb{P}\{\Lambda = \mu - 1\}. \tag{4.4}$$

The dynamics of the processes $\Phi^{(0)}$ and $\Phi^{(1)}$ in (4.2) and (4.3) are interdependent. However, if we define a new process

$$\widetilde{\boldsymbol{\Phi}} \equiv \begin{bmatrix} \widetilde{\Phi}^{(0)} \\ \widetilde{\Phi}^{(1)} \end{bmatrix} \triangleq \frac{1}{\sqrt{2}} \begin{bmatrix} \Phi^{(0)} - \Phi^{(1)} \\ \Phi^{(0)} + \Phi^{(1)} \end{bmatrix}, \tag{4.5}$$

then each of the new processes $\widetilde{\Phi}^{(0)}$ and $\widetilde{\Phi}^{(1)}$ is autonomous:

$$d\widetilde{\Phi}_t^{(0)} = \left[(\lambda + 1)\widetilde{\Phi}_t^{(0)} + \frac{\lambda(1-m)}{\sqrt{2}}\right]dt - \frac{1}{\mu}\widetilde{\Phi}_{t-}^{(0)}\,dN_t,$$

$$\widetilde{\Phi}_0^{(0)} = \frac{(1-m)\pi}{\sqrt{2}(1-\pi)},$$

(4.6)

$$d\widetilde{\Phi}_t^{(1)} = \left[(\lambda - 1)\widetilde{\Phi}_t^{(1)} + \frac{\lambda(1+m)}{\sqrt{2}}\right]dt + \frac{1}{\mu}\widetilde{\Phi}_{t-}^{(1)}\,dN_t,$$

$$\widetilde{\Phi}_0^{(1)} = \frac{(1+m)\pi}{\sqrt{2}(1-\pi)}.$$

The new coordinates $\widetilde{\Phi}^{(0)}$ and $\widetilde{\Phi}^{(1)}$ are in fact the conditional *odds-ratio* processes as in

$$\widetilde{\Phi}_t^{(0)} = \sqrt{2} \cdot \frac{\mathbb{P}\{\Lambda = \mu - 1, \theta \leq t | \mathcal{F}_t\}}{\mathbb{P}\{\theta > t | \mathcal{F}_t\}}$$

and

$$\widetilde{\Phi}_t^{(1)} = \sqrt{2} \cdot \frac{\mathbb{P}\{\Lambda = \mu + 1, \theta \leq t | \mathcal{F}_t\}}{\mathbb{P}\{\theta > t | \mathcal{F}_t\}}.$$

Therefore, both $\widetilde{\Phi}^{(0)}$ and $\widetilde{\Phi}^{(1)}$ are nonnegative processes.

Note that $m \in [-1, 1]$ in (4.4). The cases $m = \pm 1$ degenerate to Poisson disorder problems with *known* post-disorder rates, and were studied by Bayraktar, Dayanik and Karatzas [3]. Therefore, we will assume that $m \in (-1, 1)$ in the remainder.

REMARK 4.1. For every $\phi_0 \in \mathbb{R}$ and $\phi_1 \in \mathbb{R}$, let us denote by $x(t, \phi_0)$, $t \in \mathbb{R}$, and $y(t, \phi_1)$, $t \in \mathbb{R}$, the solutions of the differential equations

$$
\text{(4.7)} \quad \begin{aligned}
\frac{d}{dt} x(t, \phi_0) &= (\lambda + 1) x(t, \phi_0) + \frac{\lambda(1-m)}{\sqrt{2}}, & x(0, \phi_0) &= \phi_0, \\
\frac{d}{dt} y(t, \phi_1) &= (\lambda - 1) y(t, \phi_1) + \frac{\lambda(1+m)}{\sqrt{2}}, & y(0, \phi_1) &= \phi_1,
\end{aligned}
$$

respectively. These solutions are given by

$$
x(t, \phi_0) = -\frac{\lambda(1-m)}{\sqrt{2}(\lambda+1)} + e^{(\lambda+1)t} \left[ \phi_0 + \frac{\lambda(1-m)}{\sqrt{2}(\lambda+1)} \right], \qquad t \in \mathbb{R},
$$

$$
\text{(4.8)} \quad y(t, \phi_1) = \begin{cases} -\dfrac{\lambda(1+m)}{\sqrt{2}(\lambda-1)} \\ \quad + e^{(\lambda-1)t} \left[ \phi_1 + \dfrac{\lambda(1+m)}{\sqrt{2}(\lambda-1)} \right], & \lambda \neq 1, \quad t \in \mathbb{R}. \\ \phi_1 + \dfrac{1+m}{\sqrt{2}} t, & \lambda = 1, \end{cases}
$$

Both $x(\cdot, \phi_0)$ and $y(\cdot, \phi_1)$ have the semigroup property, that is, for every $t \in \mathbb{R}$ and $s \in \mathbb{R}$

$$\text{(4.9)} \quad x(t+s, \phi_0) = x(s, x(t, \phi_0)) \quad \text{and} \quad y(t+s, \phi_1) = y(s, y(t, \phi_1)).$$

Note from (4.6) and (4.7) that

$$
\text{(4.10)} \quad \widetilde{\Phi}_t^{(0)} = x(t - \sigma_n, \widetilde{\Phi}_{\sigma_n}^{(0)}) \quad \text{and} \quad \widetilde{\Phi}_t^{(1)} = y(t - \sigma_n, \widetilde{\Phi}_{\sigma_n}^{(1)}),
$$
$$\sigma_n \leq t < \sigma_{n+1}, n \in \mathbb{N}_0.$$

4.1. *An optimal stopping problem for the quickest detection of the Poisson disorder* In terms of the new sufficient statistics $\widetilde{\Phi}^{(1)}$ and $\widetilde{\Phi}^{(0)}$ in (4.5), (4.6), the Bayes risk of (2.6), (3.2) can be rewritten as

$$R_\tau(\pi) = 1 - \pi + \frac{c(1-\pi)}{\sqrt{2}} \cdot \mathbb{E}_0 \left[ \int_0^\tau e^{-\lambda t} \left( \widetilde{\Phi}_t^{(0)} + \widetilde{\Phi}_t^{(1)} - \frac{\lambda}{c} \sqrt{2} \right) dt \right], \qquad \tau \in \mathcal{S}.$$

Therefore, the *minimum Bayes risk* $U(\pi) \triangleq \inf_{\tau \in \mathcal{S}} R_\tau(\pi)$, $\pi \in [0,1)$, is given by

$$U(\pi) = 1 - \pi + \frac{c(1-\pi)}{\sqrt{2}} \cdot V\left(\frac{(1-m)\pi}{\sqrt{2}(1-\pi)}, \frac{(1+m)\pi}{\sqrt{2}(1-\pi)}\right), \quad (4.11)$$
$$\pi \in [0,1),$$

where $m$ is as in (4.4), the function $V(\cdot, \cdot)$ is the value function of the optimal stopping problem

$$V(\phi_0, \phi_1) \triangleq \inf_{\tau \in \mathcal{S}} \mathbb{E}_0^{\phi_0, \phi_1}\left[\int_0^\tau e^{-\lambda t} g(\widetilde{\Phi}_t^{(0)}, \widetilde{\Phi}_t^{(1)})\, dt\right],$$
$$(4.12)$$
$$g(\phi_0, \phi_1) \triangleq \phi_0 + \phi_1 - \frac{\lambda}{c}\sqrt{2}, \qquad (\phi_0, \phi_1) \in \mathbb{R}_+^2,$$

and $\mathbb{E}_0^{\phi_0, \phi_1}$ is the conditional $\mathbb{P}_0$-expectation given that $\widetilde{\Phi}_0^{(0)} = \phi_0$ and $\widetilde{\Phi}_0^{(1)} = \phi_1$.

It is clear from (4.12) that it is never optimal to stop before the process $\widetilde{\mathbf{\Phi}}$ leaves the region

$$(4.13) \qquad \mathbf{C}_0 \triangleq \left\{(\phi_0, \phi_1) \in \mathbb{R}_+^2 : \phi_0 + \phi_1 < \frac{\lambda}{c}\sqrt{2}\right\}.$$

In the next subsection we shall discuss the pathwise behavior of the process $\widetilde{\mathbf{\Phi}}$; this will give insight into the solution of the optimal stopping problem in (4.12).

4.2. *The sample-paths of the sufficient-statistic process* $\widetilde{\mathbf{\Phi}} = (\widetilde{\Phi}^{(0)}, \widetilde{\Phi}^{(1)})$.
The process $\widetilde{\Phi}^{(0)}$ jumps downward and increases between jumps; see (4.6). On the other hand, the process $\widetilde{\Phi}^{(1)}$ jumps upward, and its behavior between jumps depends on the sign of $1 - \lambda$. If $\lambda \geq 1$, then the process $\widetilde{\Phi}^{(1)}$ increases between jumps. If $0 < \lambda < 1$, then $\widetilde{\Phi}^{(1)}$ reverts to the (positive) "mean level"

$$(4.14) \qquad \phi_d \triangleq \frac{\lambda(1+m)}{(1-\lambda)\sqrt{2}}$$

between jumps; it never visits $\phi_d$ unless it starts there; and in this latter case, it stays at $\phi_d$ until the first jump and never returns to $\phi_d$ (i.e., $\phi_d > 0$ is an entrance boundary for $\widetilde{\Phi}^{(1)}$). Finally, $\phi_d$ and $1 - \lambda \neq 0$ have the same signs.

As for the solution of the optimal stopping problem in (4.12), it is worth waiting if the process $\widetilde{\mathbf{\Phi}}$ is in the region $\mathbf{C}_0$ of (4.13), or is likely to return to $\mathbf{C}_0$ shortly. The sample-paths of the process $\widetilde{\mathbf{\Phi}}$ are deterministic between jumps, and tend toward, or away from, the region $\mathbf{C}_0$. These two cases are described separately below. In both cases, however, the process $\widetilde{\mathbf{\Phi}}$ jumps in the same direction relative to its position before the jump. A jump at $(\phi_0, \phi_1)$ is an instantaneous displacement $(1/\mu)[-\phi_0 \ \phi_1]^{\mathrm{T}}$ in $\widetilde{\mathbf{\Phi}}$. Therefore, the jump direction is away from (resp. toward) the region $\mathbf{C}_0$ if $\phi_0 < \phi_1$ (resp. $\phi_0 > \phi_1$). Along a quarter of a circle in Figure 1(a), the directions of jumps at an equal distance from the origin are illustrated by the arrows. Note also that, along any fixed half-ray in $\mathbb{R}_+^2$, the jump direction [namely, the angle $\alpha$ in Figure 1(a)] does not change, but the size of the jump does.

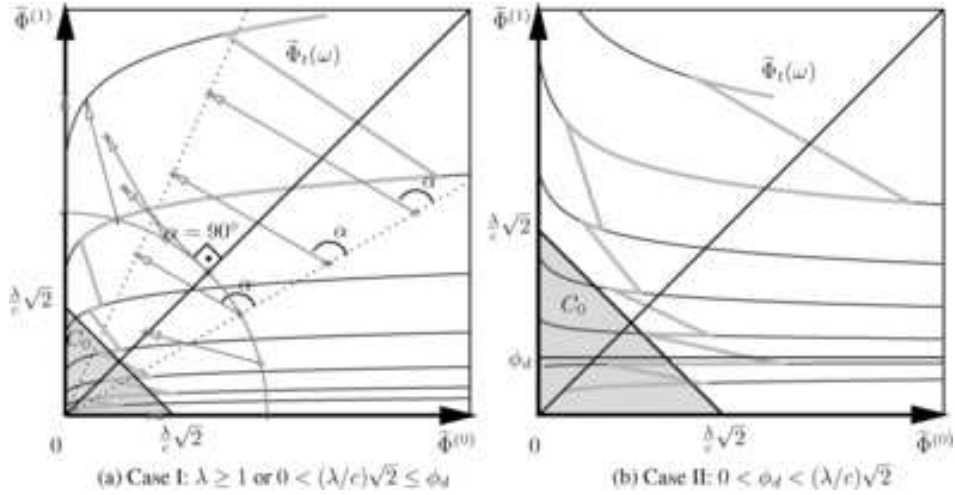

Fig. 1. *The sample-paths of $\widetilde{\mathbf{\Phi}}$.*

4.3. *Case* I: *A "large" disorder arrival rate.* Suppose that $\lambda \geq 1$ or $0 < (\lambda/c)\sqrt{2} \leq \phi_d$. Equivalently, $\lambda \geq [1 - (1+m)(c/2)]^+$ is "large." Between jumps, the process $\widetilde{\mathbf{\Phi}}$ gets farther away from the region $\mathbf{C}_0$. It may return to $\mathbf{C}_0$ by jumps only, and only if the jump originates in the region $L \triangleq \{(\phi_0, \phi_1) : \phi_0 > \phi_1\}$; see Figure 1(a). But, if $\widetilde{\Phi}^{(1)}$ reaches at or above $(\lambda/c)\sqrt{2}$, then $\widetilde{\mathbf{\Phi}}$ will never return to $\mathbf{C}_0$.

4.4. *Case* II: *A "small" disorder arrival rate.* Now suppose that $0 < \phi_d < (\lambda/c)\sqrt{2}$. Equivalently, $0 < \lambda < 1 - (1+m)(c/2)$ is "small." If the process $\widetilde{\mathbf{\Phi}}$ finds itself in a very close neighborhood of the upper-left corner of the triangular region $\mathbf{C}_0$, then it will drift into $\mathbf{C}_0$ before the next jump with positive probability. Otherwise, the behavior of the sample-paths of $\widetilde{\mathbf{\Phi}}$ relative to $\mathbf{C}_0$ is very similar to that in Case I; see Figure 1(b).

**5. A family of related optimal stopping problems.** Let us introduce for every $n \in \mathbb{N}$, the optimal stopping problem

$$V_n(\phi_0, \phi_1) \triangleq \inf_{\tau \in \mathcal{S}} \mathbb{E}_0^{\phi_0, \phi_1}\left[\int_0^{\tau \wedge \sigma_n} e^{-\lambda t} g(\widetilde{\Phi}_t^{(0)}, \widetilde{\Phi}_t^{(1)}) \, dt\right],$$

(5.1)
$$(\phi_0, \phi_1) \in \mathbb{R}_+^2,$$

obtained from (4.12) by stopping the process $\widetilde{\boldsymbol{\Phi}}$ at the $n$th jump time $\sigma_n$ of the process $N$. Since $g(\cdot, \cdot)$ in (4.12) is bounded from below by the constant $-(\lambda/c)\sqrt{2}$, the expectation in (5.1) is well defined for every stopping time $\tau \in \mathcal{S}$. In fact, $-\sqrt{2}/c \leq V_n \leq 0$ for every $n \in \mathbb{N}$. Since the sequence $(\sigma_n)_{n \geq 1}$ of jump times of the process $N$ is increasing almost surely, the sequence $(V_n)_{n \geq 1}$ is decreasing. Therefore, $\lim_{n \to \infty} V_n$ exists everywhere. It is also obvious that $V_n \geq V$, $n \in \mathbb{N}$.

PROPOSITION 5.1. *As $n \to \infty$, the sequence $V_n(\phi_0, \phi_1)$ converges to $V(\phi_0, \phi_1)$ uniformly in $(\phi_0, \phi_1) \in \mathbb{R}_+^2$. In fact, for every $n \in \mathbb{N}$ and $(\phi_0, \phi_1) \in \mathbb{R}_+^2$, we have*

(5.2)
$$\frac{\sqrt{2}}{c} \cdot \left(\frac{\mu}{\lambda + \mu}\right)^n \geq V_n(\phi_0, \phi_1) - V(\phi_0, \phi_1) \geq 0.$$

PROOF. Fix $(\phi_0, \phi_1) \in \mathbb{R}_+^2$. For $\tau \in \mathcal{S}$, $n \in \mathbb{N}$, we express $\mathbb{E}_0^{\phi_0, \phi_1}[\int_0^\tau e^{-\lambda s} \times g(\widetilde{\boldsymbol{\Phi}}_s)\, ds]$ as

$$\mathbb{E}_0^{\phi_0,\phi_1}\left[\int_0^{\tau \wedge \sigma_n} e^{-\lambda s} g(\widetilde{\boldsymbol{\Phi}}_s)\, ds\right] + \mathbb{E}_0^{\phi_0,\phi_1}\left[\mathbf{1}_{\{\tau \geq \sigma_n\}} \int_{\sigma_n}^\tau e^{-\lambda s} g(\widetilde{\boldsymbol{\Phi}}_s)\, ds\right]$$

$$\geq \mathbb{E}_0^{\phi_0,\phi_1}\left[\int_0^{\tau \wedge \sigma_n} e^{-\lambda s} g(\widetilde{\boldsymbol{\Phi}}_s)\, ds\right] - \frac{\lambda}{c}\sqrt{2} \cdot \mathbb{E}_0^{\phi_0,\phi_1}\left[\mathbf{1}_{\{\tau \geq \sigma_n\}} \int_{\sigma_n}^\tau e^{-\lambda s}\, ds\right]$$

$$\geq V_n(\phi_0, \phi_1) - \frac{\sqrt{2}}{c} \cdot \left(\frac{\mu}{\lambda + \mu}\right)^n.$$

We have used the bound $g(\phi_0, \phi_1) \geq -(\lambda/c)\sqrt{2}$ from (4.12), as well as the fact that $N$ is a Poisson process with rate $\mu$ under $\mathbb{P}_0$, and $\sigma_n$ is the $n$th jump time of $N$. Taking the infimum over $\tau \in \mathcal{S}$ gives the first inequality in (5.2). $\square$

We shall try to calculate now the functions $V_n(\cdot)$ of (5.1), following a method of Gugerli [11] and Davis [9]. Let us start by defining on the collection of bounded Borel functions $w: \mathbb{R}_+^2 \mapsto \mathbb{R}$ the operators

(5.3)
$$Jw(t, \phi_0, \phi_1) \triangleq \mathbb{E}_0^{\phi_0,\phi_1}\left[\int_0^{t \wedge \sigma_1} e^{-\lambda u} g(\widetilde{\Phi}_u^{(0)}, \widetilde{\Phi}_u^{(1)})\, du \right.$$
$$\left. + \mathbf{1}_{\{t \geq \sigma_1\}} e^{-\lambda \sigma_1} w(\widetilde{\Phi}_{\sigma_1}^{(0)}, \widetilde{\Phi}_{\sigma_1}^{(1)})\right],$$

(5.4) $\quad J_t w(\phi_0, \phi_1) \triangleq \inf_{u \in [t, \infty]} Jw(u, \phi_0, \phi_1) \quad$ for every $t \in [0, \infty]$.

The special structure of the stopping times of jump processes (see Lemma A.1 below) implies

$$
\begin{aligned}
J_0 w(\phi_0, \phi_1) = \inf_{\tau \in \mathcal{S}} \mathbb{E}_0^{\phi_0, \phi_1} \bigg[ &\int_0^{\tau \wedge \sigma_1} e^{-\lambda t} g(\widetilde{\Phi}_t^{(0)}, \widetilde{\Phi}_t^{(1)}) \, dt \\
&+ \mathbf{1}_{\{\tau \geq \sigma_1\}} e^{-\lambda \sigma_1} w(\widetilde{\Phi}_{\sigma_1}^{(0)}, \widetilde{\Phi}_{\sigma_1}^{(1)}) \bigg].
\end{aligned}
\tag{5.5}
$$

By relying on the strong Markov property of the process $N$ at its first jump time $\sigma_1$, one expects that the value function $V$ of (4.12) satisfies the equation $V = J_0 V$. Below, we show that this is indeed the case. In fact, if we define $v_n : \mathbb{R}_+^2 \mapsto \mathbb{R}$, $n \in \mathbb{N}_0$ sequentially by

$$
v_0 \equiv 0 \quad \text{and} \quad v_n \triangleq J_0 v_{n-1} \quad \forall n \in \mathbb{N}, \tag{5.6}
$$

then every $v_n$ is bounded and identical to $V_n$ of (5.1), whereas $\lim_{n \to \infty} v_n$ exists and is equal to the value function $V$ in (4.12).

Under $\mathbb{P}_0$, the first jump time $\sigma_1$ of the process $N$ has exponential distribution with rate $\mu$. Using the Fubini theorem and (4.10), we can write (5.3) as

$$
Jw(t, \phi_0, \phi_1) = \int_0^t e^{-(\lambda + \mu)u} (g + \mu \cdot w \circ S)(x(u, \phi_0), y(u, \phi_1)) \, du \tag{5.7}
$$

for every $t \in [0, \infty]$, where $x(\cdot, \phi_0)$ and $y(\cdot, \phi_1)$ are the solutions (4.8) of the ordinary differential equations in (4.7), and $S : \mathbb{R}_+^2 \mapsto \mathbb{R}_+^2$ is the linear mapping

$$
S(\phi_0, \phi_1) \triangleq \left( \left(1 - \frac{1}{\mu}\right) \phi_0, \left(1 + \frac{1}{\mu}\right) \phi_1 \right). \tag{5.8}
$$

REMARK 5.2. Using $\mu > 1$ and the explicit forms of $x(u, \phi_0)$ and $y(u, \phi_1)$ in (4.8), it is easy to check that the integrand in (5.7) is absolutely integrable on $\mathbb{R}_+$. Therefore,

$$
\lim_{t \to \infty} Jw(t, \phi_0, \phi_1) = Jw(\infty, \phi_0, \phi_1) < \infty,
$$

and the mapping $t \mapsto Jw(t, \phi_0, \phi_1) : [0, \infty] \mapsto \mathbb{R}$ is continuous. The infimum $J_t w(\phi_0, \phi_1)$ in (5.4) is attained for every $t \in [0, \infty]$.

LEMMA 5.3. *For every bounded Borel function $w : \mathbb{R}_+^2 \mapsto \mathbb{R}$, the mapping $J_0 w$ is bounded. If we define $\|w\| \triangleq \sup_{(\phi_0, \phi_1) \in \mathbb{R}_+^2} |w(\phi_0, \phi_1)| < \infty$, then*

$$
-\left( \frac{\lambda}{\lambda + \mu} \cdot \frac{\sqrt{2}}{c} + \frac{\mu}{\lambda + \mu} \cdot \|w\| \right) \leq J_0 w(\phi_0, \phi_1) \leq 0, \tag{5.9}
$$

$$(\phi_0, \phi_1) \in \mathbb{R}_+^2.$$

*If the function $w(\phi_0, \phi_1)$ is concave, then so is $J_0 w(\phi_0, \phi_1)$. If $w_1 \leq w_2$ are real-valued and bounded Borel functions defined on $\mathbb{R}_+^2$, then $J_0 w_1 \leq J_0 w_2$.*

COROLLARY 5.4.   *Every $v_n$, $n \in \mathbb{N}_0$, in (5.6) is bounded and concave, and $-\sqrt{2}/c \leq \cdots \leq v_n \leq v_{n-1} \leq v_1 \leq v_0 \equiv 0$. The limit*

$$(5.10) \qquad v(\phi_0, \phi_1) \triangleq \lim_{n \to \infty} v_n(\phi_0, \phi_1), \qquad (\phi_0, \phi_1) \in \mathbb{R}_+^2,$$

*exists, and is also bounded and concave.*

*Both $v_n : \mathbb{R}_+^2 \mapsto \mathbb{R}$, $n \in \mathbb{N}$, and $v : \mathbb{R}_+^2 \mapsto \mathbb{R}$ are continuous, increasing in each of their arguments, and their left and right partial derivatives are bounded on every compact subset of $\mathbb{R}_+^2$.*

PROPOSITION 5.5.   *For every $n \in \mathbb{N}$, the functions $v_n$ of (5.6) and $V_n$ of (5.1) coincide. For every $\varepsilon \geq 0$, let for every $n = 0, 1, \ldots, (\phi_0, \phi_1) \in \mathbb{R}_+^2$,*

$$r_n^\varepsilon(\phi_0, \phi_1) \triangleq \inf\{s \in (0, \infty] : Jv_n(s, \phi_0, \phi_1) \leq J_0 v_n(\phi_0, \phi_1) + \varepsilon\},$$

$$S_1^\varepsilon \triangleq r_0^\varepsilon(\widetilde{\boldsymbol{\Phi}}_0) \wedge \sigma_1 \quad \text{and}$$

$$S_{n+2}^\varepsilon \triangleq \begin{cases} r_{n+1}^{\varepsilon/2}(\widetilde{\boldsymbol{\Phi}}_0), & \text{if } \sigma_1 > r_{n+1}^{\varepsilon/2}(\widetilde{\boldsymbol{\Phi}}_0), \\ \sigma_1 + S_{n+1}^{\varepsilon/2} \circ \theta_{\sigma_1}, & \text{if } \sigma_1 \leq r_{n+1}^{\varepsilon/2}(\widetilde{\boldsymbol{\Phi}}_0), \end{cases}$$

*where $\theta_s$ is the shift-operator on $\Omega : N_t \circ \theta_s = N_{s+t}$. Then*

$$(5.11) \quad \mathbb{E}_0^{\phi_0, \phi_1}\left[\int_0^{S_n^\varepsilon} e^{-\lambda t} g(\widetilde{\boldsymbol{\Phi}}_t)\, dt\right] \leq v_n(\phi_0, \phi_1) + \varepsilon, \qquad n = 1, 2, \ldots, \varepsilon \geq 0.$$

PROPOSITION 5.6.   *We have $v(\phi_0, \phi_1) = V(\phi_0, \phi_1)$ for every $(\phi_0, \phi_1) \in \mathbb{R}_+^2$. Moreover, $V$ is the largest nonpositive solution $U$ of the equation $U = J_0 U$.*

LEMMA 5.7.   *Let $w : \mathbb{R}_+^2 \mapsto \mathbb{R}$ be a bounded function. For every $t \in \mathbb{R}_+$ and $(\phi_0, \phi_1) \in \mathbb{R}_+^2$,*

$$(5.12) \quad J_t w(\phi_0, \phi_1) = Jw(t, \phi_0, \phi_1) + e^{-(\lambda + \mu)t} J_0 w(x(t, \phi_0), y(t, \phi_1)).$$

COROLLARY 5.8.   *Let*

$$(5.13) \quad r_n(\phi_0, \phi_1) = \inf\{s \in (0, \infty] : Jv_n(s, (\phi_0, \phi_1)) = J_0 v_n(\phi_0, \phi_1)\}$$

*be the same as $r_n^\varepsilon(\phi_0, \phi_1)$ in Proposition 5.5 with $\varepsilon = 0$. Then*

$$(5.14) \quad r_n(\phi_0, \phi_1) = \inf\{t > 0 : v_{n+1}(x(t, \phi_0), y(t, \phi_1)) = 0\} \qquad (\inf \varnothing \equiv \infty).$$

PROOF. Let us fix $(\phi_0, \phi_1) \in \mathbb{R}_+^2$, and denote $r_n(\phi_0, \phi_1)$ by $r_n$. By Remark 5.2, we have $Jv_n(r_n, \phi_0, \phi_1) = J_0 v_n(\phi_0, \phi_1) = J_{r_n} v_n(\phi_0, \phi_1)$.

Suppose first that $r_n < \infty$. Since $J_0 v_n = v_{n+1}$, taking $t = r_n$ and $w = v_n$ in (5.12) implies that $Jv_n(r_n, \phi_0, \phi_1)$ equals

$$J_{r_n} v_n(\phi_0, \phi_1) = Jv_n(r_n, \phi_0, \phi_1) + e^{-(\lambda + \mu) r_n} v_{n+1}(x(r_n, \phi_0), y(r_n, \phi_1)).$$

Therefore, $v_{n+1}(x(r_n, \phi_0), y(r_n, \phi_1)) = 0$.

If $0 < t < r_n$, then $J_t v_n(t, \phi_0, \phi_1) > J_0 v_n(\phi_0, \phi_1) = J_{r_n} v_n(\phi_0, \phi_1) = J_t v_n(\phi_0, \phi_1)$ since $u \mapsto J_u v_n(\phi_0, \phi_1)$ is nondecreasing. Taking $t \in (0, r_n)$ and $w = v_n$ in (5.12) implies

$$J_0 v_n(\phi_0, \phi_1) = J_t v_n(\phi_0, \phi_1)$$
$$= J v_n(t, \phi_0, \phi_1) + e^{-(\lambda+\mu)t} v_{n+1}(x(t, \phi_0), y(t, \phi_1)).$$

Therefore, $v_{n+1}(x(t, \phi_0), y(t, \phi_1)) < 0$ for every $t \in (0, r_n)$, and (5.14) follows.

Suppose now that $r_n = \infty$. Then we have $v_{n+1}(x(t, \phi_0), y(t, \phi_1)) < 0$ for every $t \in (0, \infty)$ by the same argument in the last paragraph above. Hence, $\{t > 0 : v_{n+1}(x(t, \phi_0), y(t, \phi_1)) = 0\} = \varnothing$, and (5.14) still holds. □

REMARK 5.9. For every $t \in [0, r_n(\phi_0, \phi_1)]$, we have $J_t v_n(\phi_0, \phi_1) = J_0 v_n(\phi_0, \phi_1) = v_{n+1}(\phi_0, \phi_1)$. Then substituting $w(\cdot, \cdot) = v_n(\cdot, \cdot)$ in (5.12) gives the *dynamic programming equation* for the family $\{v_k(\cdot, \cdot)\}_{k \in \mathbb{N}_0}$: for every $(\phi_0, \phi_1) \in \mathbb{R}_+^2$ and $n \in \mathbb{N}_0$,

(5.15)
$$v_{n+1}(\phi_0, \phi_1) = J v_n(t, \phi_0, \phi_1) + e^{-(\lambda+\mu)t} v_{n+1}(x(t, \phi_0), y(t, \phi_1)),$$
$$t \in [0, r_n(\phi_0, \phi_1)].$$

REMARK 5.10 (Dynamic programming equation). Since $V(\cdot, \cdot)$ is bounded, and $V = J_0 V$ by Proposition 5.6, Lemma 5.7 gives

(5.16) $$J_t V(\phi_0, \phi_1) = JV(t, \phi_0, \phi_1) + e^{-(\lambda+\mu)t} V(x(t, \phi_0), y(t, \phi_1)),$$
$$t \in \mathbb{R}_+,$$

for every $(\phi_0, \phi_1) \in \mathbb{R}_+^2$; and if we define

(5.17) $$r(\phi_0, \phi_1) \triangleq \inf\{t > 0 : JV(t, \phi_0, \phi_1) = J_0 V(\phi_0, \phi_1)\},$$
$$(\phi_0, \phi_1) \in \mathbb{R}_+^2,$$

then arguments similar to those in the proof of Corollary 5.8, and (5.16), give

(5.18) $r(\phi_0, \phi_1) = \inf\{t > 0 : V(x(t, \phi_0), y(t, \phi_1)) = 0\}, \quad (\phi_0, \phi_1) \in \mathbb{R}_+^2,$

as well as the dynamic programming equation

(5.19)
$$V(\phi_0, \phi_1) = JV(t, \phi_0, \phi_1) + e^{-(\lambda+\mu)t} V(x(t, \phi_0), y(t, \phi_1)),$$
$$t \in [0, r(\phi_0, \phi_1)],$$

for the function $V(\cdot, \cdot)$ of (4.12). Because $t \mapsto Jw(t, (\phi_0, \phi_1))$ and $t \mapsto J_t w(\phi_0, \phi_1)$ are continuous for every bounded $w : \mathbb{R}_+^2 \mapsto \mathbb{R}$ [see, e.g., (5.7)], the identity (5.16) implies that $t \mapsto V(x(t, \phi_0), y(t, \phi_1))$ is continuous. Therefore, every realization of $t \mapsto V(\widetilde{\mathbf{\Phi}}_t)$ is right-continuous and has left limits.

Let us define the $\mathbb{F}$-stopping times

(5.20) $$U_\varepsilon \triangleq \inf\{t \geq 0 : V(\widetilde{\boldsymbol{\Phi}}_t) \geq -\varepsilon\}, \qquad \varepsilon \geq 0.$$

By Remark 5.10, we have

(5.21) $$V(\widetilde{\boldsymbol{\Phi}}_{U_\varepsilon}) \geq -\varepsilon \text{ on the event } \{U_\varepsilon < \infty\}.$$

PROPOSITION 5.11. *Let $M_t \triangleq e^{-\lambda t} V(\widetilde{\boldsymbol{\Phi}}_t) + \int_0^t e^{-\lambda s} g(\widetilde{\boldsymbol{\Phi}}_s)\,ds$, $t \geq 0$. For every $n \in \mathbb{N}$, $\varepsilon \geq 0$ and $(\phi_0, \phi_1) \in \mathbb{R}_+^2$, we have $\mathbb{E}_0^{\phi_0,\phi_1}[M_0] = \mathbb{E}_0^{\phi_0,\phi_1}[M_{U_\varepsilon \wedge \sigma_n}]$, that is,*

(5.22) $$V(\phi_0, \phi_1) = \mathbb{E}_0^{\phi_0,\phi_1}\left[e^{-\lambda(U_\varepsilon \wedge \sigma_n)} V(\widetilde{\boldsymbol{\Phi}}_{U_\varepsilon \wedge \sigma_n}) + \int_0^{U_\varepsilon \wedge \sigma_n} e^{-\lambda s} g(\widetilde{\boldsymbol{\Phi}}_s)\,ds\right].$$

PROPOSITION 5.12. *For every $\varepsilon \geq 0$, the stopping time $U_\varepsilon$ in (5.20) is $\varepsilon$-optimal for the problem (4.12), that is,*

$$\mathbb{E}_0^{\phi_0,\phi_1}\left[\int_0^{U_\varepsilon} e^{-\lambda s} g(\widetilde{\boldsymbol{\Phi}}_s)\,ds\right] \leq V(\phi_0, \phi_1) + \varepsilon \qquad \text{for every } (\phi_0, \phi_1) \in \mathbb{R}_+^2.$$

**6. A bound on the alarm time.** We shall show that the optimal continuation region $\mathbf{C} = \{(\phi_0, \phi_1) \in \mathbb{R}_+^2 : V(\phi_0, \phi_1) < 0\}$ is contained in some set

(6.1) $$D = \{(\phi_0, \phi_1) \in \mathbb{R}_+^2 : \phi_0 + \phi_1 < \xi^*\}$$
$$\text{for a suitable } \xi^* \in \left[\frac{\lambda + \mu}{c}\sqrt{2}, \infty\right).$$

Therefore, the region $\mathbf{C}$ has compact closure; this will be very useful in proving in the next section that $\mathbf{C}$ has a boundary which is strictly decreasing and convex.

Recall from Section 4.1 that it is not optimal to stop before the process $\widetilde{\boldsymbol{\Phi}}$ leaves the region $\mathbf{C}_0$ in (4.13). Thus, the optimal stopping time $U_0$ of Proposition 5.12 is bounded from below and above as in

(6.2) $$\tau_{C_0} \triangleq \inf\left\{t \geq 0 : \widetilde{\Phi}_t^{(0)} + \widetilde{\Phi}_t^{(1)} \geq \frac{\lambda}{c}\sqrt{2}\right\} \leq U_0 \leq \tau_D$$
$$\triangleq \inf\{t \geq 0 : \widetilde{\Phi}_t^{(0)} + \widetilde{\Phi}_t^{(1)} \geq \xi^*\}$$

in terms of the exit times $\tau_{C_0}$ and $\tau_D$ of the process $\widetilde{\boldsymbol{\Phi}}$ from the regions $C_0$ and $D$, respectively. The constant threshold $\xi^*$ in (6.1) is essentially determined by the number $(\lambda + \mu)\sqrt{2}/c$ [see (6.5), (6.9) and (6.11)], and our calculations below suggest that they are close. Therefore, the bounds in (6.2) may prove useful in practice. The difference $[(\lambda+\mu)/c]\sqrt{2} - (\lambda/c)\sqrt{2} = (\mu/c)\sqrt{2}$ between the thresholds that determine the latest and the earliest alarm times is also meaningful. It increases as $\mu/c$ increases: waiting longer is encouraged if the new information arrives at a rate higher than the cost for detection delay per unit time when the disorder has already happened.

Finally, we prove in Lemma 6.1 that $\tau_D$ in (6.2) has finite expectation. Therefore,

$$\mathbb{E}_0^{\phi_0,\phi_1}[U_0] \leq \mathbb{E}_0^{\phi_0,\phi_1}[\tau_D] < \infty \qquad \text{for every } (\phi_0,\phi_1) \in \mathbb{R}_+^2.$$

Let $\tau \in \mathcal{S}$ be any $\mathbb{F}$-stopping time. By Lemma A.1, there is a constant $t \geq 0$ such that $\tau \wedge \sigma_1 = t \wedge \sigma_1$ almost surely. Therefore

$$\mathbb{E}_0^{\phi_0,\phi_1}\left[\int_0^\tau e^{-\lambda s} g(\widetilde{\boldsymbol{\Phi}}_s)\, ds\right]$$

$$= \mathbb{E}_0^{\phi_0,\phi_1}\left[\int_0^{\tau \wedge \sigma_1} e^{-\lambda s} g(\widetilde{\boldsymbol{\Phi}}_s)\, ds\right] + \mathbb{E}_0^{\phi_0,\phi_1}\left[\mathbf{1}_{\{\tau \geq \sigma_1\}} \int_{\sigma_1}^\tau e^{-\lambda s} g(\widetilde{\boldsymbol{\Phi}}_s)\, ds\right]$$

(6.3)
$$\geq \mathbb{E}_0^{\phi_0,\phi_1}\left[\int_0^t \mathbf{1}_{\{s \leq \sigma_1\}} e^{-\lambda s} g(x(s,\phi_0), y(s,\phi_1))\, ds\right]$$

$$- \frac{\sqrt{2}}{c} \cdot \mathbb{E}_0^{\phi_0,\phi_1}[\mathbf{1}_{\{t \geq \sigma_1\}} e^{-\lambda \sigma_1}]$$

$$= \int_0^t e^{-(\lambda+\mu)s}\left[g(x(s,\phi_0), y(s,\phi_1)) - \frac{\mu}{c}\sqrt{2}\right] ds.$$

The inequality follows from $g(\phi_0,\phi_1) \geq g(0,0) = -(\lambda/c)\sqrt{2}$; see (4.12). The functions $x(\cdot,\phi_0)$ and $y(\cdot,\phi_1)$ are the solutions of (4.7) (see Remark 4.1), and $\sigma_1$ has exponential distribution with rate $\mu$ under $\mathbb{P}_0$. Clearly, if for $0 < s < \infty$ we have

$$(6.4) \quad 0 < g(x(s,\phi_0), y(s,\phi_1)) - \frac{\mu}{c}\sqrt{2} = x(s,\phi_0) + y(s,\phi_1) - \frac{\lambda+\mu}{c}\sqrt{2},$$

then (6.3) implies that $\mathbb{E}_0^{\phi_0,\phi_1}[\int_0^\tau e^{-\lambda s} g(\widetilde{\boldsymbol{\Phi}}_s)\, ds] > 0$ for every $\mathbb{F}$-stopping time $\tau \neq 0$ almost surely (since the filtration $\mathbb{F}$ is right-continuous, the probability of $\{\tau \geq 0\} \in \mathcal{F}_0$ equals zero or 1). Thus, "*stopping immediately*" *is optimal at every* $(\phi_0,\phi_1)$ *for which* (6.4) *holds.*

If $\lambda \geq 1$, then $s \mapsto x(s,\phi_0)$ and $s \mapsto y(s,\phi_1)$ are increasing for every $(\phi_0,\phi_1) \in \mathbb{R}_+^2$; see (4.7) and Figure 2(a). Therefore, $x(s,\phi_0) + y(s,\phi_1) > x(0,\phi_0) + y(0,\phi_1) = \phi_0 + \phi_1$ for every $0 < s < \infty$. Hence, (6.4) holds, and therefore it is optimal to stop immediately outside the region

$$(6.5) \qquad D_1 \triangleq \left\{(\phi_0,\phi_1) \in \mathbb{R}_+^2 : \phi_0 + \phi_1 < \frac{\lambda+\mu}{c}\sqrt{2}\right\} \qquad \text{if } \lambda \geq 1.$$

Suppose now that $0 < \lambda < 1$; equivalently, $\phi_d$ of (4.14) is positive. Then $s \mapsto x(s,\phi_0)$ is increasing for every $\phi_0 \in \mathbb{R}_+$. For $\phi_1 = \phi_d$, the derivative $dy(s,\phi_d)/ds$ in (4.7) vanishes for every $0 < s < \infty$. The mapping $s \mapsto y(s,\phi_1)$ is increasing if $\phi_1 \in [0,\phi_d)$, decreasing if $\phi_1 \in (\phi_d, \infty)$, and $\phi(s,\phi_d) = \phi_d$ for every $0 \leq s < \infty$; see (4.7) and Figures 2(b),(c). The derivative

$$(6.6) \quad \frac{d}{dt}[x(t,\phi_0) + y(t,\phi_1)] = (\lambda+1)x(t,\phi_0) + (\lambda-1)y(t,\phi_1) + \lambda\sqrt{2}$$

of the right-hand side of (6.4) [see also (4.7)] vanishes if $s \mapsto (x(s,\phi_0), y(s,\phi_1))$ meets at $s = t$ the line

$$(6.7) \quad \ell : (\lambda+1)x + (\lambda-1)y + \lambda\sqrt{2} = 0 \quad \text{or} \quad y = \frac{1+\lambda}{1-\lambda}x + \frac{\lambda}{1-\lambda}\sqrt{2}.$$

Since $m \in (-1,1)$, the "mean-level" $\phi_d$ in (4.14) and the $y$-intercept of the line $\ell$ in (6.7) are related as in

$$\phi_d = \frac{\lambda}{1-\lambda} \cdot \frac{1+m}{\sqrt{2}} < \frac{\lambda}{1-\lambda}\sqrt{2}.$$

Because $\ell$ is increasing, this relationship implies that the line $\ell$ in contained in $\mathbb{R}_+ \times (\phi_d, \infty)$ [see Figure 2(b),(c)]. However, every curve $t \mapsto (x(t,\phi_0), y(t,\phi_1))$ starting at some $(\phi_0, \phi_1)$ in $\mathbb{R}_+ \times (\phi_d, \infty)$ is "decreasing," and the derivative in (6.6) is increasing. Therefore, any curve $t \mapsto (x(t,\phi_0), y(t,\phi_1))$, $(\phi_0, \phi_1) \in \mathbb{R}_+^2$ may meet $\ell$ at most once, and

(6.8) if $t \mapsto (x(t,\phi_0), y(t,\phi_1))$ meets the line $\ell$ at $t_\ell = t_\ell(\phi_0, \phi_1)$, then $t \mapsto x(t,\phi_0) + y(t,\phi_1)$ is decreasing (resp., increasing) on $[0, t_\ell]$ (resp., on $[t_\ell, \infty)$). Otherwise, $t \mapsto x(t,\phi_0) + y(t,\phi_1)$ is increasing on $[0, \infty)$.

Consider now the first of two possible cases: the line $\ell$ does not meet $D_1$ of (6.5); that is, $\lambda/(1-\lambda) \geq (\lambda+\mu)/c$, as in Figure 2(b). Then $\phi_0 + \phi_1 \geq$ $(\lambda+\mu)\sqrt{2}/c$ for every $(\phi_0,\phi_1) \in \ell$. Therefore, (6.8) implies that (6.4) holds, that is, it is optimal to stop immediately, outside

$$(6.9) \quad D_1 = \left\{(\phi_0,\phi_1) \in \mathbb{R}_+^2 : \phi_0 + \phi_1 < \frac{\lambda+\mu}{c}\sqrt{2}\right\} \qquad \text{if } \frac{\lambda}{1-\lambda} \geq \frac{\lambda+\mu}{c}.$$

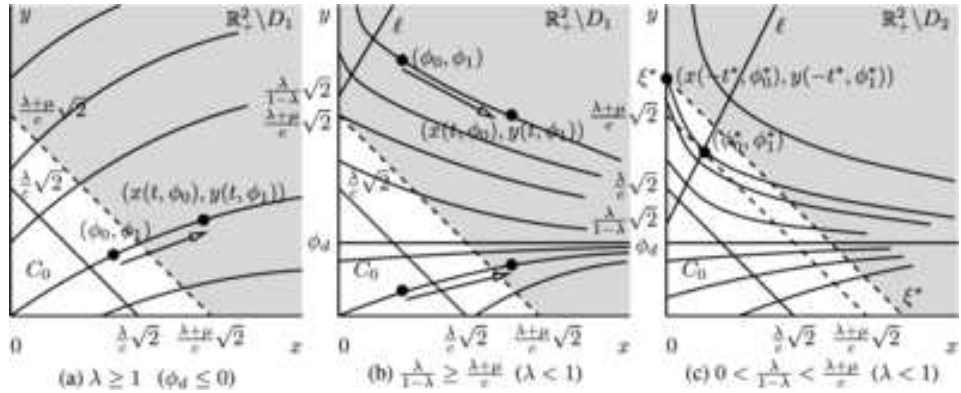

Fig. 2. *Region D.*

In the second case, the line $\ell$ of (6.7) meets the region $D_1$, that is, $0 < \lambda/(1-\lambda) < (\lambda+\mu)/c$; see Figure 2(c). Let us denote by $(\phi_0^*, \phi_1^*)$ the point at the intersection of the line $\ell$ and the boundary $x + y - (\lambda+\mu)\sqrt{2}/c = 0$ of the region $D_1$. By running the time "backward," we can find $\xi^*$ (and $t^*$) such that

$$(6.10) \qquad (0,\xi^*) = (x(-t^*,\phi_0^*), y(-t^*,\phi_1^*)).$$

Indeed, using (4.8), we can obtain first $t^* \geq 0$ by solving $0 = x(-t^*,\phi_0^*)$, and then $\xi^* \triangleq y(-t^*,\phi_1^*)$. By the semigroup property (4.9), we have

$$x(t^*,0) = x(t^*, x(-t^*,\phi_0^*)) = x(t^* + (-t^*), \phi_0^*) = x(0,\phi_0^*) = \phi_0^*,$$
$$y(t^*,\xi^*) = y(t^*, y(-t^*,\phi_1^*)) = y(t^* + (-t^*), \phi_1^*) = y(0,\phi_1^*) = \phi_1^*.$$

Hence, the curve $t \mapsto (x(t,0), y(t,\xi^*))$, $t \geq 0$, meets $\ell$ at $(\phi_0^*, \phi_1^*)$, and $t_\ell$ in (6.8) equals $t^*$; see Figure 2(c). Therefore, (6.8) implies that

$$x(t,0) + y(t,\xi^*) \geq x(t^*,0) + y(t^*,\xi^*) = \phi_0^* + \phi_1^* = \frac{\lambda+\mu}{c}\sqrt{2}, \qquad 0 \leq t < \infty.$$

In particular, $\xi^* = 0 + \xi^* = x(0,0) + y(0,\xi^*) \geq (\lambda+\mu)\sqrt{2}/c$. We are now ready to show that it is optimal to stop immediately outside the region

$$(6.11) \quad D_2 \triangleq \{(\phi_0,\phi_1) \in \mathbb{R}_+^2 : \phi_0 + \phi_1 < \xi^*\} \qquad \text{if } 0 < \frac{\lambda}{1-\lambda} < \frac{\lambda+\mu}{c},$$

where $\xi^*$ is as in (6.10). The curve $t \mapsto (x(t,0), y(t,\xi^*))$ divides $\mathbb{R}_+^2$ into two connected components each containing the region $D_1$ of (6.5) and

$$M \triangleq (\mathbb{R}_+^2 \setminus D_2) \cap \{(x,y) \in \mathbb{R}_+^2 : (\lambda+1)x + (\lambda-1)y + \lambda\sqrt{2} < 0\},$$

respectively [see (6.7)]. Every curve $t \mapsto (x(t,\phi_0), y(t,\phi_1))$, $t \geq 0$, starting at $(\phi_0,\phi_1)$ in $M$ will stay in the same component as $M$. Therefore, the curve intersects the line $\ell$ away from $D_1$, and (6.8) implies that (6.4) is satisfied for every $(\phi_0,\phi_1) \in M$.

For $(\phi_0,\phi_1) \in (\mathbb{R}_+^2 \setminus D_2) \cap \{(x,y) \in \mathbb{R}_+^2 : (\lambda+1)x + (\lambda-1)y + \lambda\sqrt{2} \geq 0\}$, the curve $t \mapsto (x(t,\phi_0), y(t,\phi_1))$, $t \geq 0$, does not meet $\ell$; therefore, $t \mapsto x(t,\phi_0) + y(t,\phi_1)$ increases by (6.8) and

$$x(t,\phi_0) + y(t,\phi_1) > x(0,\phi_0) + y(0,\phi_1)$$
$$= \phi_0 + \phi_1 \geq \xi^* \geq \frac{\lambda+\mu}{c}\sqrt{2}, \qquad 0 < s < \infty.$$

Thus, the sufficient condition (6.4) for the optimality of immediate stopping holds for every $(\phi_0,\phi_1) \in \mathbb{R}_+^2 \setminus D_2$.

LEMMA 6.1. *Let $\tau_D$ be the exit time of the process $\widetilde{\mathbf{\Phi}}$ from the region $D$ in (6.1). Then $\mathbb{E}_0^{\phi_0,\phi_1}[\tau_D]$ is finite for every $(\phi_0, \phi_1) \in \mathbb{R}_+^2$.*

PROOF. Let $f(\phi_0, \phi_1) \triangleq \phi_0 + \phi_1$, $(\phi_0, \phi_1) \in \mathbb{R}_+^2$. Using the explicit form of the infinitesimal generator $\widetilde{\mathcal{A}}$ of the process $\widetilde{\mathbf{\Phi}}$ in (A.4), we obtain

$$\widetilde{\mathcal{A}} f(\phi_0, \phi_1) = (\lambda + 1)\phi_0 + \frac{\lambda(1-m)}{\sqrt{2}} + (\lambda - 1)\phi_1 + \frac{\lambda(1+m)}{\sqrt{2}}$$

(6.12)
$$+ \mu \left[ \left(1 - \frac{1}{\mu}\right) \phi_0 + \left(1 + \frac{1}{\mu}\right) \phi_1 - (\phi_0 + \phi_1) \right]$$

$$= \lambda(\phi_0 + \phi_1 + \sqrt{2}) \geq \lambda\sqrt{2}$$

for every $(\phi_0, \phi_1) \in \mathbb{R}_+^2$. Since $f(\cdot, \cdot)$ is bounded on $\overline{D}$ of (6.1) and $\tau_D \wedge t$, $t \geq 0$, is a bounded $\mathbb{F}$-stopping time, (A.3) holds for $\tau = \tau_D \wedge t$. Then we have

$$\xi^* \left(1 + \frac{1}{\mu}\right) \geq \mathbb{E}_0^{\phi_0,\phi_1}[f(\widetilde{\mathbf{\Phi}}_{\tau_D \wedge t})]$$

(6.13)
$$= f(\phi_0, \phi_1) + \mathbb{E}_0^{\phi_0,\phi_1}\left[ \int_0^{\tau_D \wedge t} \widetilde{\mathcal{A}} f(\widetilde{\mathbf{\Phi}}_t)\, dt \right]$$

$$\geq \lambda\sqrt{2} \mathbb{E}_0^{\phi_0,\phi_1}[\tau_D \wedge t], \qquad t \geq 0.$$

The process $\widetilde{\mathbf{\Phi}}$ may leave the region $D$ in (6.1) continuously or by a jump. Since $f(S(\phi_0, \phi_1)) = (1 + 1/\mu)\phi_0 + (1 - 1/\mu)\phi_1 \leq (1 + 1/\mu)(\phi_0 + \phi_1) = (1 + 1/\mu)f(\phi_0, \phi_1) \leq (1 + 1/\mu)\xi^*$ for every $(\phi_0, \phi_1) \in D$, and this upper bound is larger than $\xi^*$, the first inequality in (6.13) follows. The second inequality is due to (6.12). Finally, the monotone convergence theorem and (6.13) imply that $\mathbb{E}_0^{\phi_0,\phi_1}[\tau_D]$ is finite. □

**7. The solution.** In Proposition 5.1, we showed that the function $V(\phi_0, \phi_1)$ of our original optimal stopping problem in (4.12) is approximated *uniformly* in $(\phi_0, \phi_1) \in \mathbb{R}_+^2$ by the decreasing sequence $\{V_n(\phi_0, \phi_1)\}_{n \in \mathbb{N}}$ of the value functions of the optimal stopping problems in (5.1). The value functions $V_n(\cdot, \cdot) = v_n(\cdot, \cdot)$, $n \in \mathbb{N}$ can be calculated sequentially by setting $v_0 \equiv 0$, and

(7.1)
$$v_{n+1}(\phi_0, \phi_1) = J_0 v_n(\phi_0, \phi_1) = \inf_{t \in [0,\infty]} J v_n(t, \phi_0, \phi_1),$$
$$(\phi_0, \phi_1) \in \mathbb{R}_+^2,$$

where the operator $J$ is defined in (5.3); see Proposition 5.5.

Finding the infimum in (7.1) is not as formidable as it may appear. By Proposition 5.5, the infimum in (7.1) is always attained [i.e., the case $\varepsilon = 0$ in (5.11)]. By Corollary 5.8, it is attained at the time, $r_n(\phi_0, \phi_1)$, the deterministic continuous curve $t \mapsto (x(t,\phi_0), y(t,\phi_1))$ in (4.7) exits from the set

$$\{(\phi_0, \phi_1) \in \mathbb{R}_+^2 : v_{n+1}(\phi_0, \phi_1) < 0\} \subseteq \{(\phi_0, \phi_1) \in \mathbb{R}_+^2 : v(\phi_0, \phi_1) < 0\} \subseteq D;$$

here $D$ is the triangular region in (6.1), and the last inclusion is proven in Section 6. Therefore, the search for the infimum in (7.1) can be confined for every $n \in \mathbb{N}$ to

$$(7.2) \quad Jv_n(t, \phi_0, \phi_1) = \int_0^t e^{-(\lambda+\mu)u} [g + \mu \cdot v_n \circ S](x(u,\phi_0), y(u,\phi_1))\, du,$$
$$t \in [0, \bar{r}(\phi_0, \phi_1)],$$

over the interval $t \in [0, \bar{r}(\phi_0, \phi_1)]$, where

$$\bar{r}(\phi_0, \phi_1) \triangleq \inf\{t \geq 0 : x(t, \phi_0) + y(t, \phi_1) \geq \xi^*\}, \qquad (\phi_0, \phi_1) \in \mathbb{R}_+^2,$$

is the (bounded) exit time of the curve $t \mapsto (x(t,\phi_0), y(t,\phi_1))$ from the region $D$ of (6.1).

Finally, the error in approximating $V(\cdot, \cdot)$ of (4.12) by $\{v_n(\cdot, \cdot)\}_{n \in \mathbb{N}}$ in (7.1) can be controlled. For every $\varepsilon > 0$,

$$(7.3) \quad \frac{\sqrt{2}}{c}\left(\frac{\mu}{\lambda+\mu}\right)^n < \varepsilon \quad \Longrightarrow \quad -\varepsilon \leq V(\phi_0, \phi_1) - v_n(\phi_0, \phi_1) \leq 0,$$
$$(\phi_0, \phi_1) \in \mathbb{R}_+^2,$$

by Propositions 5.1 and 5.5. The exponential rate of the uniform convergence of $\{v_n(\cdot, \cdot)\}_{n \in \mathbb{N}}$ to $V(\cdot, \cdot)$ on $\mathbb{R}_+^2$ in (7.3) may also reduce the computational burden by allowing relatively small number of iterations in (7.1).

In the remainder, we draw attention to certain special cases where the value function $V(\cdot, \cdot)$ can be calculated *gradually* at each iteration in (7.1); see Proposition 8.3. In the meantime, we shall give a precise geometric description of the stopping regions

(7.4) $\boldsymbol{\Gamma}_n \triangleq \{(\phi_0, \phi_1) \in \mathbb{R}_+^2 : v_n(\phi_0, \phi_1) = 0\}, \qquad \mathbf{C}_n \triangleq \mathbb{R}_+^2 \setminus \boldsymbol{\Gamma}_n, \qquad n \in \mathbb{N},$

(7.5) $\quad \boldsymbol{\Gamma} \triangleq \{(\phi_0, \phi_1) \in \mathbb{R}_+^2 : v(\phi_0, \phi_1) = 0\}, \qquad \mathbf{C} \triangleq \mathbb{R}_+^2 \setminus \boldsymbol{\Gamma},$

and describe the optimal stopping strategies.

**8. The structure of the stopping regions.** By Proposition 5.12, the set $\boldsymbol{\Gamma}$ is the *optimal stopping region* for the problem (4.12). Namely, stopping at the first hitting time $U_0 = \inf\{t \in \mathbb{R}_+ : \widetilde{\boldsymbol{\Phi}}_t \in \boldsymbol{\Gamma}\}$ of the process $\widetilde{\boldsymbol{\Phi}} = (\widetilde{\Phi}^{(0)}, \widetilde{\Phi}^{(1)})$ to the set $\boldsymbol{\Gamma}$ is optimal for (4.12).

Similarly, we shall call each set $\boldsymbol{\Gamma}_n$, $n \in \mathbb{N}$, a *stopping region* for the family of optimal stopping problems in (5.1). However, unlike the case above, we need the first $n$ stopping regions, $\boldsymbol{\Gamma}_1, \ldots, \boldsymbol{\Gamma}_n$, in order to describe an optimal stopping time for the problem of (5.1). Using Corollary 5.8, the optimal stopping time $S_n \equiv S_n^0$ in Proposition 5.5 for $V_n$ of (5.1) may be described as follows: *Stop if the process $\widetilde{\Phi}$ hits $\Gamma_n$ before $N$ jumps. If $N$ jumps before $\widetilde{\Phi}$ reaches $\Gamma_n$, then wait, and stop if $\widetilde{\Phi}$ hits $\Gamma_{n-1}$ before the next jump of $N$, and so on. If the rule is not met before the $(n-1)$st jump of $N$, then stop at the earliest of the hitting time of $\Gamma_1$ and the next jump time of $N$*. See Figure 3(b) for three realizations of the stopping time $S_2$.

We shall call each $\mathbf{C}_n \triangleq \mathbb{R}_+^2 \setminus \mathbf{\Gamma}_n$, $n \in \mathbb{N}$, a *continuation region* for the family of optimal stopping problems in (5.1), and $\mathbf{C} \triangleq \mathbb{R}_+^2 \setminus \mathbf{\Gamma}$ the *optimal continuation region* for (4.12). The stopping regions are related by

$$\mathbb{R}_+^2 \setminus D \subset \mathbf{\Gamma} \subset \cdots \subset \mathbf{\Gamma}_n \subset \mathbf{\Gamma}_{n-1} \subset \cdots \subset \mathbf{\Gamma}_1 \subset \mathbb{R}_+^2 \setminus \mathbf{C}_0 \quad \text{and}$$
(8.1)
$$\mathbf{\Gamma} = \bigcap_{n=1}^{\infty} \mathbf{\Gamma}_n,$$

since the sequence of nonpositive functions $\{v_n\}_{n \in \mathbb{N}}$ is decreasing, and $v = \lim_{n \to \infty} \downarrow v_n$ by Lemma 5.4. The sets $D$ and $\mathbf{C}_0$ are defined in (6.1) and (4.13), respectively. Since $v_n$, $n \in \mathbb{N}$ and $v$ are concave and continuous mappings from $\mathbb{R}_+^2$ into $(-\infty, 0]$ by Lemma 5.4, the stopping regions $\mathbf{\Gamma}_n$, $n \in \mathbb{N}$, and $\mathbf{\Gamma}$ are convex and closed. Let us define the functions $\gamma_n : \mathbb{R}_+ \mapsto \mathbb{R}_+$, $n \in \mathbb{N}$, and $\gamma : \mathbb{R}_+ \mapsto \mathbb{R}_+$ by [see also Figure 3(a)]

$$\gamma_n(x) \triangleq \inf\{y \in \mathbb{R}_+ : (x, y) \in \mathbf{\Gamma}_n\}, \qquad x \in \mathbb{R}_+,$$
$$\gamma(x) \triangleq \inf\{y \in \mathbb{R}_+ : (x, y) \in \mathbf{\Gamma}\}, \qquad x \in \mathbb{R}_+,$$

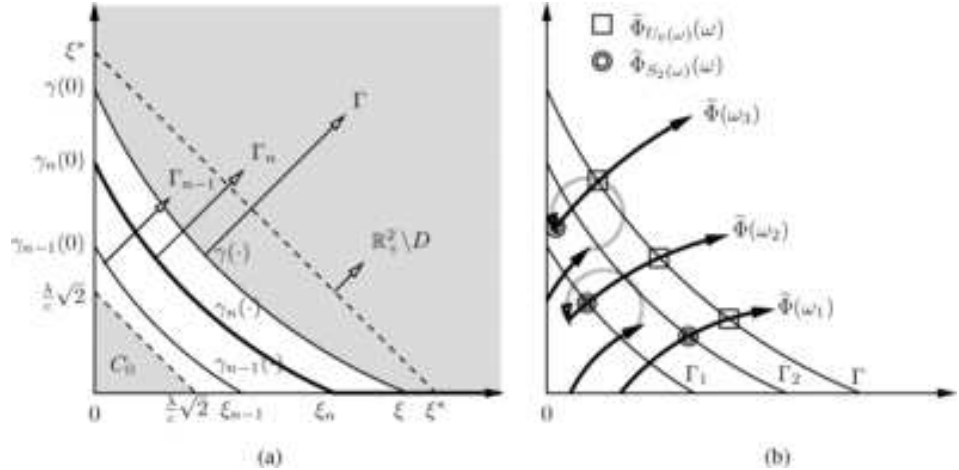

Fig. 3. (a) *The stopping regions (each arrow at the boundary of a region points toward the interior of that region), and* (b) *three sample-paths and the optimal stopping times $S_2$ and $U_0$ for the optimal stopping problems $V_2$ in* (5.1) *and $V$ in* (4.12), *respectively.*

and the numbers

$$\xi_n \triangleq \inf\{x \in \mathbb{R}_+ : \gamma_n(x) = 0\}, \qquad n \in \mathbb{N}, \quad \text{and} \quad \xi \triangleq \inf\{x \in \mathbb{R}_+ : \gamma(x) = 0\}.$$

Then the stopping regions $\mathbf{\Gamma}_n$, $n \in \mathbb{N}$, and $\mathbf{\Gamma}$ are the convex and closed epigraphs of the functions $\gamma_n(\cdot)$, $n \in \mathbb{N}$, and $\gamma(\cdot)$, respectively. Therefore, $\gamma_n(\cdot)$, $n \in \mathbb{N}$, and $\gamma(\cdot)$ are convex and continuous mappings from $\mathbb{R}_+$ into $\mathbb{R}_+$.

By the set-inclusions in (8.1), we have $(\lambda/c)\sqrt{2} \leq \xi_{n-1} \leq \xi_n \leq \xi \leq \xi^*$ for the same $\xi^* \in \mathbb{R}_+$ in the description (6.1) of the set $D$. Since $v_n$, $n \in \mathbb{N}$, and $v$ vanish on $\mathbb{R}_+ \setminus D = \{(\phi_0, \phi_1) \in \mathbb{R}_+^2 : \phi_0 + \phi_1 \geq \xi^*\}$ by (8.1), the functions $\gamma_n(\cdot)$, $n \in \mathbb{N}$, and $\gamma(\cdot)$ vanish on $[\xi^*, \infty)$. However, $\xi_n$ and $\xi$ are the *smallest* zeros of the continuous functions $\gamma_n(\cdot)$, $n \in \mathbb{N}$, and $\gamma(\cdot)$, respectively. Since both functions are also nonnegative and convex, the function $\gamma_n(\cdot)$, $n \in \mathbb{N}$ [resp. $\gamma(\cdot)$] is zero on $[\xi_n, \infty)$ (resp. on $[\xi, \infty)$) and *strictly* decreasing on $[0, \xi_n]$ (resp. on $[0, \xi]$). For future reference, we now summarize our results.

PROPOSITION 8.1. *There are decreasing, convex and continuous mappings $\gamma_n : \mathbb{R}_+ \mapsto \mathbb{R}_+$, $n \in \mathbb{N}$, and $\gamma : \mathbb{R}_+ \mapsto \mathbb{R}_+$ such that*

$$\mathbf{\Gamma}_n = \{(\phi_0, \phi_1) \in \mathbb{R}_+^2 : \phi_1 \geq \gamma_n(\phi_0)\}, \qquad n \in \mathbb{N},$$

$$\mathbf{\Gamma} = \{(\phi_0, \phi_1) \in \mathbb{R}_+^2 : \phi_1 \geq \gamma(\phi_0)\}.$$

*The sequence $\{\gamma_n(\phi_0)\}_{n \in \mathbb{N}}$ is increasing and $\gamma(\phi_0) = \lim \uparrow \gamma_n(\phi_0)$ for every $\phi_0 \in \mathbb{R}_+$. There are numbers*

$$(8.2) \qquad \frac{\lambda}{c}\sqrt{2} \leq \xi_1 \leq \cdots \leq \xi_{n-1} \leq \xi_n \leq \cdots \leq \xi < \xi^* < \infty$$

*such that $\gamma_n(\cdot)$, $n \in \mathbb{N}$ [resp., $\gamma(\cdot)$] is strictly decreasing on $[0, \xi_n]$, $n \in \mathbb{N}$ (resp., $[0, \xi]$) and vanishes on $[\xi_n, \infty)$, $n \in \mathbb{N}$ (resp., $[\xi, \infty)$). Moreover,*

$$(8.3) \quad \frac{\lambda}{c}\sqrt{2} \leq \gamma_1(0) \leq \cdots \leq \gamma_{n-1}(0) \leq \gamma_n(0) \leq \cdots \leq \gamma(0) < \xi^* < \infty.$$

*The number $\xi^*$ is the same as in the definition of the set $D$ in* (6.1).

NOTATION 8.2. Let $S : \mathbb{R}_+^2 \mapsto \mathbb{R}_+^2$ be the same linear map as in (5.8).

(N1) For any subset $R \subseteq \mathbb{R}_+^2$,

$$S^{-(n+1)}(R) \triangleq S^{-1}(S^{-n}(R)), \qquad n \in \mathbb{N},$$
$$S^{-1}(R) \triangleq \{(x, y) \in \mathbb{R}_+^2 : S(x, y) \in R\},$$
$$S^{n+1}(R) \triangleq S(S^n(R)), \qquad n \in \mathbb{N},$$
$$S(R) \triangleq \{S(x, y) \in \mathbb{R}_+^2 : (x, y) \in R\},$$

and $S^0(R) = S(S^{-1}(R)) = S^{-1}(S(R)) = R$.

(N2) For every singleton $\{(x,y)\} \subseteq \mathbb{R}_+^2$, we write

$$S^m(\{x,y\}) = S^m(x,y) = \left(\left(1 - \frac{1}{\mu}\right)^m x, \left(1 + \frac{1}{\mu}\right)^m y\right), \qquad m \in \mathbb{Z}.$$

(N3) For any function $g:\mathbb{R}_+ \mapsto \mathbb{R}_+$, we define the function $S^n[g]:\mathbb{R}_+ \mapsto \mathbb{R}_+$, $n \in \mathbb{Z}$ by

$$S^n[g](x) \triangleq \inf\{y \in \mathbb{R}_+ : (x,y) \in S^n(\mathrm{epi}(g))\}, \qquad x \in \mathbb{R}_+.$$

That is, $S^n[g]$ is the function whose epigraph is the set $S^n(\mathrm{epi}(g))$. Note that we use $S^n(\cdot)$ and $S^n[\cdot]$ to distinguish the sets and the functions.

(N4) For every subset $R$ of $\mathbb{R}_+^2$, we denote by $\mathrm{cl}(R)$ its *closure* in $\mathbb{R}_+^2$ and by $\mathrm{int}(R)$ its *interior*. We shall denote the *support* of a function $g:\mathbb{R}_+ \mapsto \mathbb{R}_+$ by

$$\mathrm{supp}(g) = \mathrm{cl}(\{x \in \mathbb{R}_+ : g(x) > 0\}).$$

The process $\widetilde{\Phi}$ jumps into the region $\Gamma$ [resp. $S^{-n}(\Gamma)$, $n \in \mathbb{N}$] if the process $N$ jumps while $\widetilde{\Phi}$ is in the region $S^{-1}(\Gamma)$ [resp. $S^{-(n+1)}(\Gamma)$, $n \in \mathbb{N}$]. Clearly, if the process $\widetilde{\Phi}$ can never leave the region $S^{-1}(\Gamma)$ before a jump, then the value functions $V(\cdot,\cdot)$ and $V_1(\cdot,\cdot)$ in (5.1) must coincide on the region $S^{-1}(\Gamma)$.

PROPOSITION 8.3. *Suppose that*

(8.4)
$$\forall n \in \mathbb{N} : (\phi_0, \phi_1) \in S^{-n}(\Gamma) \implies (x(t,\phi_0), y(t,\phi_1)) \in S^{-n}(\Gamma),$$
$$t \in [0,\infty),$$

*holds. Then for every $n \in \mathbb{N}$, we have*

(8.5)
$$V(\phi_0, \phi_1) = V_n(\phi_0, \phi_1) = V_{n+1}(\phi_0, \phi_1) = \cdots \qquad \forall (\phi_0, \phi_1) \in S^{-n}(\Gamma),$$
$$S^{-n}(\Gamma) \cap \Gamma = S^{-n}(\Gamma) \cap \Gamma_n = S^{-n}(\Gamma) \cap \Gamma_{n+1} = \cdots,$$
$$S^{-n}(\Gamma) \cap \mathbf{C} = S^{-n}(\Gamma) \cap \mathbf{C}_n = S^{-n}(\Gamma) \cap \mathbf{C}_{n+1} = \cdots.$$

Since $\Gamma$ and $\Gamma_n$ are convex and closed, and $S(\cdot,\cdot)$ is a linear mapping, the sets $S^{-n}(\Gamma)$ and $S^{-n}(\Gamma_n)$, $n \in \mathbb{N}$, are convex and closed. The sets $\Gamma$ and $\Gamma_n$, $n \in \mathbb{N}$, are the epigraphs of the continuous functions $\gamma(\cdot)$ and $\gamma_n(\cdot)$, $n \in \mathbb{N}$, in Proposition 8.3, respectively. Therefore,

(8.6)
$$S^{-n}(\Gamma) = \{(x,y) \in \mathbb{R}_+^2 : y \geq S^{-n}[\gamma](x)\},$$
$$S^{-n}(\Gamma_n) = \{(x,y) \in \mathbb{R}_+^2 : y \geq S^{-n}[\gamma_n](x)\},$$

are the epigraphs of the functions $S^{-n}[\gamma](\cdot)$ and $S^{-n}[\gamma_n](\cdot)$ for every $n \in \mathbb{N}_0$. These functions are decreasing, continuous and convex. In fact,

$$(8.7) \quad S^{-n}[\gamma](x) = \left(\frac{\mu}{\mu+1}\right)^n \gamma\left(\left(\frac{\mu-1}{\mu}\right)^n x\right), \qquad x \in \mathbb{R}_+, n \in \mathbb{Z},$$

and the function $S^{-n}[\gamma_n](\cdot)$ is obtained by replacing $\gamma$ with $\gamma_n$ in (8.7). The support of the functions $S^{-n}[\gamma](\cdot)$ and $S^{-n}[\gamma_n](\cdot)$ are

$$(8.8) \quad \begin{aligned} \operatorname{supp}(S^{-n}[\gamma]) &= \left[0, \left(\frac{\mu}{\mu-1}\right)^n \xi\right], \\ \operatorname{supp}(S^{-n}[\gamma_n]) &= \left[0, \left(\frac{\mu}{\mu-1}\right)^n \xi_n\right], \end{aligned}$$

respectively, for every $n \in \mathbb{Z}$. By Proposition 8.1, the functions $S^{-n}[\gamma](\cdot)$ and $S^{-n}[\gamma_n](\cdot)$ are strictly decreasing on their supports.

Since $S^{-n}[\gamma](0) = (\mu/(\mu+1))^n \gamma(0) < \gamma(0)$ and $S^{-n}[\gamma](\xi) > 0 = \gamma(\xi)$ for every $n \in \mathbb{N}$, the functions $S^{-n}[\gamma](\cdot)$ and $\gamma(\cdot)$ intersect, and

$$(8.9) \quad x_n(\gamma) \triangleq \min\{x \in \mathbb{R}_+ : S^{-n}[\gamma](x) = \gamma(x)\} \in (0, \infty), \qquad n \in \mathbb{N}.$$

COROLLARY 8.4. *Suppose that* (8.4) *holds. Then for every* $n = 1, 2, \ldots$ *and* $k \geq n$, *we have* $x_n \equiv x_n(\gamma) = x_n(\gamma_k)$, *and*

$$(8.10) \quad S^{-n}(\mathbf{\Gamma}) \cap \mathbf{C} \cap ([0, x_n] \times \mathbb{R}_+) = S^{-n}(\mathbf{\Gamma}_k) \cap \mathbf{C}_k \cap ([0, x_n] \times \mathbb{R}_+).$$

*Particularly, we have* $\gamma(x) = \gamma_n(x)$ *for every* $x \in [0, x_n]$, *and*

$$(8.11) \quad V(x,y) = V_n(x,y) \qquad \forall (x,y) \in S^{-n}(\mathbf{\Gamma}_n) \cap \mathbf{C}_n \cap ([0, x_n] \times \mathbb{R}_+),$$
$$n \in \mathbb{N}.$$

PROOF. Let us fix any $k \geq n \in \mathbb{N}$. Since the value functions $V(\cdot, \cdot)$ and $V_k(\cdot, \cdot)$ are equal on the region $S^{-n}(\mathbf{\Gamma})$ by Proposition 8.3, the boundaries of the regions $\mathbf{\Gamma}$ and $\mathbf{\Gamma}_k$ coincide in the region $S^{-n}(\mathbf{\Gamma})$. Particularly, we have

$$(8.12) \quad \gamma(x) = \gamma_k(x) \qquad \text{for every } x \in [0, x_n(\gamma)]$$

since $S^{-n}[\gamma](x) < \gamma(x)$ for every $x \in [0, x_n(\gamma))$. Therefore,

$$(8.13) \quad \begin{aligned} S^{-n}[\gamma](x) &= S^{-n}[\gamma_k](x) \\ &\text{for every } x \in \left[0, \left(\frac{\mu}{\mu-1}\right)^n x_n(\gamma)\right] \supset [0, x_n(\gamma)]. \end{aligned}$$

Now, (8.12) and (8.13) imply that $x_n(\gamma) = x_n(\gamma_k)$, and (8.6) implies that

$$S^{-n}(\mathbf{\Gamma}) \cap \mathbf{C} \cap ([0, x_n(\gamma)] \times \mathbb{R}_+) = S^{-n}(\mathbf{\Gamma}_k) \cap \mathbf{C}_k \cap ([0, x_n(\gamma)] \times \mathbb{R}_+).$$

Equality (8.11) follows immediately from Proposition 8.3. □

The identity in (8.11) suggests that, in a *finite number of iterations* of (7.1), we can find the restrictions of the value function $V(\cdot,\cdot)$ and the continuation region $\mathbf{C}$ to the set $\mathbb{R}_+ \times [B,\infty)$ for any $B > 0$, when the condition (8.4) holds:

*Step* A.1. Calculate the value function $v_1(0,y)$ for every $y \in [0,\xi^*]$, and determine $\gamma(0) = \gamma_1(0) = \inf\{y \in \mathbb{R}_+ : v_1(x,y) = 0\} \in (0,\xi^*)$; see (8.3), Corollary 8.4.

*Step* A.2. Given any $B > 0$, find the smallest $n \in \mathbb{N}$ such that

$$(8.14) \qquad B > \left(\frac{\mu}{\mu+1}\right)^n \gamma_1(0) = \left(\frac{\mu}{\mu+1}\right)^n \gamma_n(0) = S^{-n}[\gamma_n](0).$$

Because every $S^{-m}[\gamma_m](\cdot)$, $m \in \mathbb{N}$, is decreasing, this implies $\mathbb{R}_+ \times [B,\infty) \subset S^{-n}(\mathbf{\Gamma}_n)$; see (8.6). We also have $n \leq \min\{m \in \mathbb{N} : B > (\mu/(\mu+1))^m \xi^*\}$ since $\gamma_1(0) \in (0,\xi^*)$.

*Step* A.3. Calculate $v_n(\phi_0, \phi_1)$ for every $(\phi_0, \phi_1) \in \mathbb{R}_+ \setminus D$ by (7.1), where $D$ is as in (6.1). By (8.1), $D \subseteq \mathbf{\Gamma}_n$ and $v_n \equiv 0$ on $D$.

Then the value functions $V(\cdot,\cdot)$ and $v_n(\cdot,\cdot)$ are equal on $\mathbb{R}_+ \times [B,\infty)$ and $(\mathbb{R}_+ \times [B,\infty)) \cap \mathbf{C} = (\mathbb{R}_+ \times [B,\infty)) \cap \mathbf{C}_n$.

The next lemma implies that we can calculate the exact value function $V(\cdot,\cdot)$ under condition (8.4) on the set $\mathbb{R}_+ \times (0,\infty)$ along an increasing sequence of sets $\mathbb{R}_+ \times [B_n,\infty)$, and on $\mathbb{R}_+ \times \{0\}$ by the continuity of the function $V(\cdot,\cdot)$ on $\mathbb{R}_+^2$.

LEMMA 8.5. *Suppose that* (8.4) *holds. Let* $\xi^*$ *be the same number as in the definition of the region* $D$ *in* (6.1). *Then* $\lim_{n\to\infty} S^{-n}(\mathbf{\Gamma}) = \mathbb{R}_+ \times (0,\infty)$, *and*

$$(8.15) \qquad \mathbb{R}_+ \times \left[\left(\frac{\mu}{1+\mu}\right)^n \cdot \xi^*, +\infty\right) \subseteq S^{-n}(\mathbf{\Gamma}), \qquad n \in \mathbb{N}.$$

PROOF. Recall from (8.1) that $\mathbb{R}_+ \times [\xi^*, \infty) \subset \mathbb{R}_+^2 \setminus D \subset \mathbf{\Gamma}$. The rectangle on the left-hand side in (8.15) is the same set as $S^{-n}(\mathbb{R}_+ \times [\xi^*, \infty)) \subset S^{-n}(\mathbf{\Gamma})$. But, (8.15) implies that $\mathbb{R}_+ \times (0,\infty) \subseteq \liminf_{n\to\infty} S^{-n}(\mathbf{\Gamma})$.

On the other hand, for every $x \in \mathbb{R}_+$, there exists number $N(x)$ such that $S^n(x,0) = ((1-1/\mu)^n x, 0) \notin \mathbf{\Gamma}$, $n \geq N(x)$. Then $(x,0) \notin S^{-n}(\mathbf{\Gamma})$ for every $n \geq N(x)$. This implies that $\limsup_{n\to\infty} S^{-n}(\mathbf{\Gamma}) \subseteq \mathbb{R}_+ \times (0,\infty)$. □

REMARK 8.6. Every set $S^{-n}(\mathbf{\Gamma})$, $n \in \mathbb{N}$, is separated from its complement by the strictly decreasing, convex and continuous function $S^{-n}[\gamma](x)$, $x \in [0, (\mu/(\mu-1))^n \xi]$. Therefore, the condition (8.4) will be satisfied, for example, if the mappings $t \mapsto x(t,\phi_0)$, $t \in \mathbb{R}_+$, and $t \mapsto y(t,\phi_1)$, $t \in \mathbb{R}_+$, are increasing for every $(\phi_0,\phi_1) \in \mathbb{R}_+^2$. We have seen in Section 4.3 that this is always the case when $\lambda$ is "large."

Thus, if $\lambda$ is "large," then there is a sequence of sets $\mathbb{R}_+ \times [B_n, \infty)$, $n \in \mathbb{N}$, increasing to $\mathbb{R}_+ \times (0, \infty)$ in the limit, such that $V(\cdot, \cdot) = v_n(\cdot, \cdot)$ on $\mathbb{R}_+ \times [B_n, \infty)$ for every $n \in \mathbb{N}$. See also Section 10 below.

**9. The boundaries of the stopping regions.** We shall show that the optimization in (7.1) can be avoided in principle, and $v_1, v_2, \ldots$ can be calculated by integration.

Note that we obtain $Jv_n(t, \phi_0, \phi_1)$ in (7.1) by integrating the function $[g + \mu \cdot v_n \circ S](\cdot, \cdot)$ along the curve $u \mapsto (x(u, \phi_0), y(u, \phi_1))$ on $u \in [0, t]$; see (5.3). Therefore, the infimum in (7.1) is determined by the excursions of $u \mapsto (x(u, \phi_0), y(u, \phi_1))$, $u \in \mathbb{R}_+$, into the regions where the sign of the continuous mapping $[g + \mu \cdot v_n \circ S](\cdot, \cdot)$ is negative and positive.

LEMMA 9.1. *For every $n \in \mathbb{N}$, we have*

$$(9.1) \qquad A_n \triangleq \{(x, y) \in \mathbb{R}_+^2 : [g + \mu \cdot v_n \circ S](x, y) < 0\} \subseteq \mathbf{C}_{n+1}.$$

PROOF. Let $(\phi_0, \phi_1) \in A_n$. Since the function $u \mapsto [g + \mu \cdot v_n \circ S](x(u, \phi_0), y(u, \phi_1))$ is continuous, there exists some $t = t(\phi_0, \phi_1) > 0$ such that

$$Jv_n(t, \phi_0, \phi_1) = \int_0^t e^{-(\lambda + \mu)u}[g + \mu \cdot v_n \circ S](x(u, \phi_0), y(u, \phi_1))\, du < 0.$$

Hence, $v_{n+1}(\phi_0, \phi_1) = J_0 v_n(\phi_0, \phi_1) \leq Jv_n(t, \phi_0, \phi_1) < 0$, and $(\phi_0, \phi_1) \in \mathbf{C}_{n+1}$. □

For certain cases, the regions $A_n$ and $\mathbf{C}_{n+1}$ coincide, that is, the continuation region $\mathbf{C}_{n+1}$ for $v_{n+1}(\cdot, \cdot)$ can be found immediately when the value function $v_n(\cdot, \cdot)$ is available. Then $v_{n+1} \equiv 0$ on $\mathbf{\Gamma}_{n+1} = \mathbb{R}_+^2 \setminus \mathbf{C}_{n+1}$, and we calculate $v_{n+1}(\cdot, \cdot)$ on $\mathbf{C}_{n+1}$ by the integration

$$
\begin{aligned}
v_{n+1}(\phi_0, \phi_1) &= Jv_n(t, \phi_0, \phi_1)|_{t = r_n(\phi_0, \phi_1)} \\
(9.2) \qquad &= \int_0^{r_n(\phi_0, \phi_1)} e^{-(\lambda + \mu)u}[g + \mu v_n \circ S](x(u, \phi_0), y(u, \phi_1))\, du,
\end{aligned}
$$

$$(\phi_0, \phi_1) \in \mathbf{C}_{n+1},$$

of the function $[g + \mu \cdot v_n \circ S](\cdot, \cdot)$ over the curve $(x(\cdot, \phi_0), y(\cdot, \phi_1))$ until the exit time $r_n(\phi_0, \phi_1)$ [see (5.14)] of the continuous curve $u \mapsto (x(u, \phi_0), y(u, \phi_1))$, $u \in \mathbb{R}_+$ from the continuation region $\mathbf{C}_{n+1}$.

The region $A_n$ in (9.1) has properties very similar to those of the continuation region $\mathbf{C}_{n+1}$; compare Lemma 9.2 and Proposition 8.1. For example, both sets are separated from their complements by a strictly decreasing, convex and continuous function which stays flat on the $x$-axis for all large $x$ values.

For every $n \in \mathbb{N}$, let us define the function $a_n : \mathbb{R}_+ \mapsto \mathbb{R}_+$ by

$$
\begin{aligned}
(9.3) \quad a_n(x) &\triangleq \inf\{y \geq 0 : (x,y) \in \mathbb{R}_+^2 \setminus A_n\} \\
&= \inf\{y \geq 0 : [g + \mu v_n \circ S](x,y) \geq 0\}.
\end{aligned}
$$

The function $a_n(\cdot)$ is finite since, given any $x \in \mathbb{R}_+$, we have $[g + \mu \cdot v_n \circ S](x,y) > 0$ for every large $y \in \mathbb{R}_+$. Recall that the function $v_n(\cdot, \cdot)$ vanishes outside the bounded region $\mathbf{C}_n$. The linear mapping $S : \mathbb{R}_+^2 \mapsto \mathbb{R}_+^2$ in (5.8) is increasing in both $x$ and $y$. The affine mapping $g : \mathbb{R}_+^2 \mapsto \mathbb{R}$ in (4.12) is also increasing and grows unboundedly in both $x$ and $y$.

Similarly, given for large $x \in \mathbb{R}_+$, we have $[g + \mu \cdot v_n \circ S](x,y) \geq 0$, $\forall\, y \in \mathbb{R}_+$. Therefore, $a_n(x) = 0$ for $x \in [\alpha, \infty)$ for some $\alpha \geq 0$, and

$$
(9.4) \quad \alpha_n \triangleq \inf\{x \geq 0 : a_n(x) = 0\} \quad \text{is finite.}
$$

The set $\mathbb{R}_+^2 \setminus A_n = \{(x,y) \in \mathbb{R}_+^2 : [g + \mu \cdot v_n \circ S](x,y) \geq 0\}$ is convex and closed since $v_n(\cdot, \cdot)$ is concave and continuous, $S(\cdot, \cdot)$ is linear and $g(\cdot, \cdot)$ is affine. Because $\mathbb{R}_+^2 \setminus A_n$ is the epigraph of $a_n(\cdot)$, this implies that $a_n(\cdot)$ is a convex and continuous mapping from $\mathbb{R}_+$ into $\mathbb{R}_+$.

The function $a_n(\cdot)$ does not vanish identically on $\mathbb{R}_+$; in particular, $a_n(0) > 0$ since the continuous function $[g + \mu \cdot v_n \circ S](x,y)$ is strictly negative at $(x,y) = (0,0)$:

$$
[g + \mu \cdot v_n \circ S](0,0) = g(0,0) + \mu \cdot v_n(0,0) \leq g(0,0) = -\frac{\lambda}{c}\sqrt{2} < 0.
$$

Because $a_n(\cdot)$ is continuous, this implies that the number $\alpha_n$ in (9.4) is strictly positive. Since $a_n(\cdot)$ is convex and vanishes for every large $x \in \mathbb{R}_+$, it is *strictly decreasing* on $[0, \alpha_n)$, and equals zero on $[\alpha_n, \infty)$.

LEMMA 9.2. *For $n \geq 1$, there exist a number $\alpha_n \in (0, \infty)$ and a strictly decreasing, convex and continuous mapping $a_n : [0, \alpha_n] \mapsto \mathbb{R}_+$ such that $a_n(\alpha_n) = 0$, and*

$$
(9.5) \quad \{(x, a_n(x)); x \in [0, \alpha_n]\} = \{(x,y) \in \mathbb{R}_+^2; [g + \mu \cdot v_n \circ S](x,y) = 0\}.
$$

*Moreover, the continuous mapping $(x,y) \mapsto [g + \mu \cdot v_n \circ S](x,y)$, $n \geq 1$, is strictly increasing in each argument, and for every $n \geq 1$,*

$$
\begin{aligned}
(9.6) \quad &\{(x,y) \in [0, \alpha_n) \times \mathbb{R}_+; y < a_n(x)\} \\
&= \{(x,y) \in \mathbb{R}_+^2; [g + \mu \cdot v_n \circ S](x,y) < 0\} \equiv A_n.
\end{aligned}
$$

Next, we shall relate the regions $A_n$ in (9.1) and $\mathbf{C}_{n+1}$, and their boundaries $a_n(\cdot)$ and $\gamma_{n+1}(\cdot)$, respectively, for every $n \in \mathbb{N}$.

Using the characterization of the stopping regions $\mathbf{\Gamma}_n$, $n \in \mathbb{N}$, in Proposition 8.1 in terms of the switching curves $\gamma_n(\cdot)$, the exit time $r_n(\cdot, \cdot)$ in Corollary 5.8 can be expressed as

$$(9.7) \quad r_n(\phi_0, \phi_1) = \inf\{t > 0 : y(t, \phi_1) = \gamma_{n+1}(x(t, \phi_0))\}, \qquad (\phi_0, \phi_1) \in \mathbb{R}_+^2,$$

since the functions $x(\cdot, \phi_0)$ and $y(\cdot, \phi_1)$ in (4.8) are continuous. Because every $\gamma_{n+1}(\cdot)$, $n \in \mathbb{N}$, is bounded, the function $r_n(\cdot, \cdot)$ is real-valued. Thus

$$0 < r_n(\phi_0, \phi_1) < \infty \qquad \text{for every } (\phi_0, \phi_1) \in \mathbf{C}_{n+1}.$$

Therefore, the (smallest) minimizer $r_n(\phi_0, \phi_1)$ of the function $t \mapsto Jv_n(t, \phi_0, \phi_1)$ as in (5.13), is an interior point of $(0, \infty]$ for every $(\phi_0, \phi_1) \in \mathbf{C}_{n+1}$, and the derivative $\partial Jv_n(t, \phi_0, \phi_1)/\partial t$ vanishes at $t = r_n(\phi_0, \phi_1)$. Using (5.7) and (9.7) gives

$$\begin{aligned}
0 &= [g + \mu \cdot v_n \circ S](x(t, \phi_0), y(t, \phi_1))|_{t=r_n(\phi_0, \phi_1)} \\
(9.8) \quad &= [g + \mu \cdot v_n \circ S](x(t, \phi_0), \gamma_{n+1}(x(t, \phi_0)))|_{t=r_n(\phi_0, \phi_1)}, \\
&\qquad\qquad\qquad\qquad\qquad\qquad\qquad\qquad (\phi_0, \phi_1) \in \mathbf{C}_{n+1}.
\end{aligned}$$

Let us denote the *boundary* of $\mathbf{\Gamma}_{n+1}$ by

$$(9.9) \quad \partial \mathbf{\Gamma}_{n+1} \triangleq \{(x, \gamma_{n+1}(x)) : x \in [0, \xi_{n+1}]\},$$

and define the *entrance* and *exit* boundaries of $\mathbf{\Gamma}_{n+1}$ by

$$\begin{aligned}
\partial \mathbf{\Gamma}_{n+1}^e &\triangleq \{(x(r_n(\phi_0, \phi_1), \phi_0), \gamma_{n+1}(y(r_n(\phi_0, \phi_1), \phi_1))), \\
&\qquad\qquad\qquad\qquad \text{for some } (\phi_0, \phi_1) \in \mathbf{C}_{n+1}\}, \\
(9.10) \\
\partial \mathbf{\Gamma}_{n+1}^x &\triangleq \{(\phi_0, \phi_1) \in \mathbf{\Gamma}_{n+1} : (x(t, \phi_0), y(t, \phi_1)) \in \mathbf{C}_{n+1}, \\
&\qquad\qquad\qquad t \in (0, \delta] \text{ for some } \delta > 0\},
\end{aligned}$$

respectively. The path $t \mapsto (x(t, \phi_0), y(t, \phi_1))$ starts at some $(\phi_0, \phi_1) \in \mathbf{C}_{n+1}$ and enters the region $\mathbf{\Gamma}_{n+1}$ (for the first time) at the entrance boundary $\partial \mathbf{\Gamma}_{n+1}^e$. Similarly, for every $(\phi_0, \phi_1) \in \partial \mathbf{\Gamma}_{n+1}^x$, the path $t \mapsto (x(t, \phi_0), y(t, \phi_1))$ exits $\mathbf{\Gamma}_{n+1}$ immediately.

REMARK 9.3. By Lemma 9.5 below, the entrance boundary $\partial \mathbf{\Gamma}_{n+1}^e$ is a subset of the boundary $\partial A_n$ of the region $A_n$ in (9.1). Clearly, the curve $t \mapsto (x(t, \phi_0), y(t, \phi_1))$ starting at any $(\phi_0, \phi_1) \in \partial \mathbf{\Gamma}_{n+1}^e \subseteq \partial A_n$ cannot return immediately into the region $A_n$ [otherwise $Jv_n(t, \phi_0, \phi_1) < 0$ for some $t > 0$ and $(\phi_0, \phi_1) \in \mathbf{C}_{n+1}$]. In the theory of Markov processes, every element of $\partial \mathbf{\Gamma}_{n+1}^e$ (resp. $\partial \mathbf{\Gamma}_{n+1}^x$) is a regular boundary point of the domain $A_n$ (resp. the interior of $\mathbf{\Gamma}_{n+1}$) with respect to the process $\widetilde{\mathbf{\Phi}}$.

REMARK 9.4. Observe that for every $(\phi_0, \phi_1) \in \partial\mathbf{\Gamma}_{n+1}^x$, the quantity $r_n(\phi_0, \phi_1)$ in (5.14) is the *return time* of the curve $t \mapsto (x(t, \phi_0), y(t, \phi_1))$ to the stopping region $\mathbf{\Gamma}_{n+1}$ and is also strictly positive. Therefore, the first-order necessary optimality condition in (9.8) also holds on the exit boundary $\partial\mathbf{\Gamma}_{n+1}^x$. Thus,

$$0 = [g + \mu \cdot v_n \circ S](x(t, \phi_0), \gamma_{n+1}(x(t, \phi_0)))|_{t=r_n(\phi_0, \phi_1)},$$
(9.11)
$$(\phi_0, \phi_1) \in \mathbf{C}_{n+1} \cup \partial\mathbf{\Gamma}_{n+1}^x.$$

9.1. *The entrance boundary* $\partial\mathbf{\Gamma}_{n+1}^e$. Since all of the functions in (9.11) are continuous, (9.11) and the definition of the entrance boundary $\partial\mathbf{\Gamma}_{n+1}^e$ in (9.10) imply

(9.12) $\quad [g + \mu \cdot v_n \circ S](x, y) = 0, \qquad (x, y) \in \partial\mathbf{\Gamma}_{n+1}^e.$

The next lemma follows from (9.12), Lemma 9.2 and the continuity of the function $[g + \mu \cdot v_n \circ S](\cdot, \cdot)$.

LEMMA 9.5. *For every* $n \in \mathbb{N}$, *let* $\alpha_n \in \mathbb{R}_+$ *and* $a_n : \mathbb{R}_+ \mapsto \mathbb{R}_+$ *be the same as in Lemma* 9.2. *Then* $\mathrm{cl}(\partial\mathbf{\Gamma}_{n+1}^e) \subseteq \{(x, a_n(x)) : x \in [0, \alpha_n]\}$, $n \in \mathbb{N}$.

COROLLARY 9.6. *For any* $n \in \mathbb{N}$, *if the equality* $\partial\mathbf{\Gamma}_{n+1} = \mathrm{cl}(\partial\mathbf{\Gamma}_{n+1}^e)$ *holds, then*

(9.13) $\quad \partial\mathbf{\Gamma}_{n+1} = \{(x, a_n(x)) : x \in [0, \alpha_n]\}.$

*In other words,* $\xi_{n+1} = \alpha_n$, *and* $\gamma_{n+1}(x) = a_n(x)$ *for every* $x \in [0, \xi_{n+1}] \equiv [0, \alpha_n]$, *and*

(9.14) $\quad \mathbf{C}_{n+1} = \{(x, y) : [g + \mu \cdot v_n \circ S](x, y) < 0\}.$

PROOF. By Lemma 9.5, $\{(x, \gamma_{n+1}(x)) : x \in [0, \xi_{n+1}]\} = \partial\mathbf{\Gamma}_{n+1} \subseteq \{(x, a_n(x)) : x \in [0, \alpha_n]\}$. Since $\gamma_{n+1}(\cdot)$ and $a_n(\cdot)$ are strictly decreasing, continuous functions which equal zero at the right-hand point of their domains, they must be identical. Finally,

$$\mathbf{C}_{n+1} = \mathbb{R}_+^2 \setminus \mathbf{\Gamma}_{n+1} = \{(x, y) \in [0, \xi_{n+1}) \times \mathbb{R}_+ : y < \gamma_{n+1}(x)\}$$
$$= \{(x, y) \in [0, \alpha_n) \times \mathbb{R}_+ : y < a_n(x)\}$$
$$= \{(x, y) \in \mathbb{R}_+^2 : [g + \mu \cdot v_n \circ S](x, y) < 0\},$$

where the last equality follows from (9.6). $\square$

If the disorder arrival rate $\lambda$ is large, then every point on the boundary $\partial\mathbf{\Gamma}_{n+1}$ of the stopping region $\mathbf{\Gamma}_{n+1}$ belongs to the entrance boundary $\partial\mathbf{\Gamma}_{n+1}^e$;

see Section 10. Therefore, the stopping boundary $\partial \mathbf{\Gamma}_{n+1}$ for the value function $v_{n+1}(\cdot, \cdot)$ is determined as in Corollary 9.6, as soon as the value function $v_n(\cdot, \cdot)$ is calculated. Using this observation, the main solution method described at the beginning of Section 7 can be tailored into a more efficient algorithm; see Section 10 and Figure 4.

The exit boundary $\partial \mathbf{\Gamma}^x_{n+1}$ may not always be nonempty. If it is nonempty, it is also determined by the entrance boundary $\partial \mathbf{\Gamma}^e_{n+1}$, and the general solution method can be similarly enhanced in this case; see Section 11.

9.2. *The exit boundary* $\partial \mathbf{\Gamma}^x_{n+1}$. Using the semigroup property in (4.9) of the functions $x(\cdot, \cdot)$ and $y(\cdot, \cdot)$, along with a change of variable, we obtain

$$(9.15) \quad Jv_n(t, \phi_0, \phi_1) = -e^{-(\lambda+\mu)t} Jv_n(-t, x(t, \phi_0), y(t, \phi_1)),$$
$$t \geq 0, (\phi_0, \phi_1) \in \mathbb{R}^2_+.$$

Substituting in (5.12) $w(\cdot, \cdot) = v_n(\cdot, \cdot)$, and using the identity above, we obtain

$$(9.16) \quad J_t v_n(\phi_0, \phi_1) = e^{-(\lambda+\mu)t}[v_{n+1}(x(t, \phi_0), y(t, \phi_1)) - Jv_n(-t, x(t, \phi_0), y(t, \phi_1))]$$

for $t \geq 0, (\phi_0, \phi_1) \in \mathbb{R}^2_+$. Since $v_{n+1}(x(r_n(\phi_0, \phi_1), \phi_0), y(r_n(\phi_0, \phi_1), \phi_1)) = 0$, and $J_{r_n(\phi_0,\phi_1)} v_n(\phi_0, \phi_1) = v_{n+1}(\phi_0, \phi_1)$, the equality in (9.16) at $t = r_n(\phi_0, \phi_1)$ gives

$$(9.17) \quad v_{n+1}(\phi_0, \phi_1) = [-e^{-(\lambda+\mu)t} Jv_n(-t, x(t, \phi_0), y(t, \phi_1))]|_{t=r_n(\phi_0, \phi_1)}$$

for $(\phi_0, \phi_1) \in \mathbb{R}^2_+$. Recall from Section 9.1 that $(x(r_n(\phi_0, \phi_1), \phi_0), y(r_n(\phi_0, \phi_1), \phi_1)) \in \partial \mathbf{\Gamma}^e_{n+1}$ for every $(\phi_0, \phi_1) \in \mathbf{C}_{n+1} \cup \partial \mathbf{\Gamma}^x_{n+1}$. Therefore, (9.17) implies that we can both calculate the value function $v_{n+1}(\cdot, \cdot)$ and find the continuation region $\mathbf{C}_{n+1}$ by backtracking the curves $t \mapsto (x(-t, \phi_0), y(-t, \phi_1))$ from every point $(\phi_0, \phi_1) \in \partial \mathbf{\Gamma}^e_{n+1}$ on the entrance boundary. Let us define for every $(\phi_0, \phi_1) \in \mathbb{R}^2_+, n \geq 0$,

$$(9.18) \quad \hat{r}(\phi_0, \phi_1) \triangleq \inf\{t \geq 0 : (x(-t, \phi_0), y(-t, \phi_1)) \notin \mathbb{R}^2_+\},$$
$$\hat{r}_n(\phi_0, \phi_1) \triangleq \inf\{t \in (0, \hat{r}(\phi_0, \phi_1)] : -Jv_n(-t, \phi_0, \phi_1) \geq 0\},$$

where the infimum of an empty set is infinity. Since the mapping $t \mapsto Jv_n(t, \phi_0, \phi_1)$ is continuous, we have $Jv_n(-\hat{r}_n(\phi_0, \phi_1), \phi_0, \phi_1) = 0$ if $0 < \hat{r}_n(\phi_0, \phi_1) < \infty$.

LEMMA 9.7. *The entrance boundary* $\partial \mathbf{\Gamma}^e_{n+1}$ *determines the exit boundary* $\partial \mathbf{\Gamma}^x_{n+1}$, *the continuation region* $\mathbf{C}_{n+1}$, *and the value function* $v_{n+1}(\cdot, \cdot)$

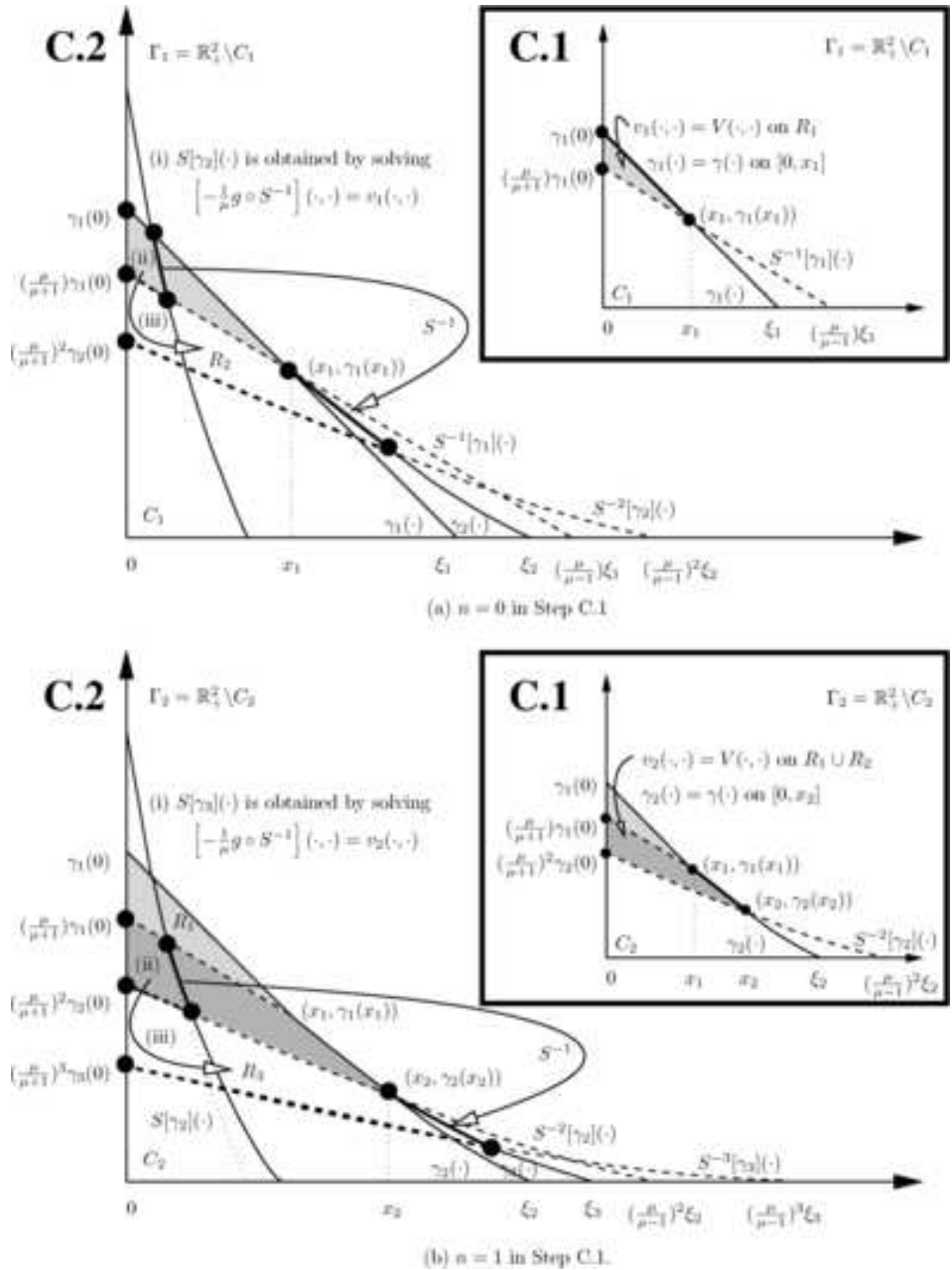

FIG. 4. *Case* I: $\lambda$ *is "large." The illustration of Method* C: *Steps* C.1 *and* C.2 *when* (a) $n=0$, *and* (b) $n=1$ *in Step* C.1.

*on* $\mathbf{C}_{n+1}$:

$$\partial\mathbf{\Gamma}^x_{n+1} = \{(x(-t,\phi_0), y(-t,\phi_1))|_{t=\hat{r}_n(\phi_0,\phi_1)} : (\phi_0,\phi_1) \in \partial\mathbf{\Gamma}^e_{n+1},$$
$$\hat{r}_n(\phi_0,\phi_1) \le \hat{r}(\phi_0,\phi_1)\},$$
$$\mathbf{C}_{n+1} = \{(x(-t,\phi_0), y(-t,\phi_1)) : (\phi_0,\phi_1) \in \partial\mathbf{\Gamma}^e_{n+1}, t \in (0, \check{r}_n(\phi_0,\phi_1)]\} \setminus \partial\mathbf{\Gamma}^x_{n+1},$$

*where* $\check{r}_n(\phi_0,\phi_1) = \hat{r}_n(\phi_0,\phi_1) \wedge \hat{r}(\phi_0,\phi_1)$. *For every* $(\phi_0,\phi_1) \in \partial\mathbf{\Gamma}^e_{n+1}$

$$v_{n+1}(x(-t,\phi_0), y(-t,\phi_1)) = -e^{-(\lambda+\mu)t} J v_n(-t,\phi_0,\phi_1),$$
$$t \in (0, \check{r}_n(\phi_0,\phi_1)].$$

**10. Case I revisited: efficient methods for "large" post-disorder arrival rates.** This is Case I on page 13 where $\lambda \geq [1 - (1+m)(c/2)]^+$ is "large," and the sample-paths of both components of the process $\widetilde{\mathbf{\Phi}} = [\widetilde{\Phi}^{(0)}, \widetilde{\Phi}^{(1)}]^{\mathrm{T}}$ increase between the jumps; see also Figure 1(a). By the relation (4.10), the deterministic functions $t \mapsto x(t,\phi_0)$, $t \in \mathbb{R}_+$, and $t \mapsto y(t,\phi_1)$, $t \in \mathbb{R}_+$, are strictly increasing for every $(\phi_0,\phi_1) \in \mathbb{R}^2_+$.

By Remark 8.6, the identity in (8.11) between the value functions $V(\cdot,\cdot)$ and $V_n(\cdot,\cdot)$ on the set $S^{-n}(\mathbf{\Gamma}_n) \cap \mathbf{C}_n \cap ([0,x_n] \times \mathbb{R}_+)$ holds for every $n \in \mathbb{N}$. Thus, we can find the value function $V(\cdot,\cdot)$ by calculating $V_n(\cdot,\cdot)$, $n \in \mathbb{N}$, using steps A.1–A.3 after Corollary 8.4. This method can be improved further. We shall show that the optimization in each iteration of (7.1) can be avoided, and the value function $V(\cdot,\cdot)$ may be calculated in one pass over the continuation region $\mathbf{C}$; see Figure 4.

Since the boundary $\partial\mathbf{\Gamma}_{n+1}$ of the stopping region $\mathbf{\Gamma}_{n+1}$ is a (strictly) decreasing continuous curve, every point in the set $\partial\mathbf{\Gamma}_{n+1} \cap \mathrm{int}(\mathbb{R}^2_+)$ is accessible from some point in the continuation region $\mathbf{C}_{n+1}$. Therefore, we have $\partial\mathbf{\Gamma}_{n+1} = \mathrm{cl}(\partial\mathbf{\Gamma}^e_{n+1})$ for every $n \in \mathbb{N}_0$. By Corollary 9.6, the set $A_n$ in (9.1) and the continuation region $\mathbf{C}_{n+1}$ (and their boundaries) coincide for every $n \in \mathbb{N}_0$.

If the value function $v_n(\cdot,\cdot) \equiv V_n(\cdot,\cdot)$ for some $n \in \mathbb{N}_0$ is already calculated, then the boundary of the continuation region $\mathbf{C}_{n+1}$ becomes immediately available as in (9.13). In fact, (9.5) and (9.14) imply

$$(10.1) \quad S(\mathbf{C}_{n+1}) = \left\{ (x,y) \in \mathbb{R}^2_+ : v_n(x,y) < \left[-\frac{1}{\mu} \cdot g \circ S^{-1}\right](x,y) \right\},$$

$$(10.2) \quad S(\partial\mathbf{\Gamma}_{n+1}) = \left\{ (x,y) \in \mathbb{R}^2_+ : v_n(x,y) = \left[-\frac{1}{\mu} \cdot g \circ S^{-1}\right](x,y) \right\}.$$

The set on the right-hand side in (10.2) is a strictly decreasing, convex and continuous curve in $\mathbb{R}^2_+$, and it is the same as

$$S(\partial\mathbf{\Gamma}_{n+1}) = S(\text{the boundary of the set } \mathrm{epi}(\gamma_{n+1}) \cap ([0,\xi_{n+1}] \times \mathbb{R}_+))$$

$$(10.3) \quad = \text{the boundary of the set } \operatorname{epi}(S[\gamma_{n+1}]) \cap \left(\left[0, \frac{\mu-1}{\mu}\xi_{n+1}\right] \times \mathbb{R}_+\right)$$

$$= \left\{(x, S[\gamma_{n+1}](x)); x \in \left[0, \frac{\mu-1}{\mu}\xi_{n+1}\right]\right\}.$$

If we know $v_n(\cdot,\cdot)$, then we can determine the set in (10.2) of all points $(x,y) \in \mathbb{R}_+^2$ satisfying

$$(10.4) \quad v_n(x,y) = \left[-\frac{1}{\mu} \cdot g \circ S^{-1}\right](x,y) \equiv -\frac{x}{\mu-1} - \frac{y}{\mu+1} + \frac{\lambda}{c\mu}\sqrt{2},$$

and obtain the boundary function $\gamma_{n+1}(\cdot)$ after the transformation of this set by $S^{-1}$. Then we can calculate the (smallest) minimizer $r_n(\cdot,\cdot)$ of (7.1) by the relation (9.7), and the value function $v_{n+1}(\cdot,\cdot)$ by (9.2). We can continue in this manner to find the value functions $v_{n+2}(\cdot,\cdot)$, $v_{n+3}(\cdot,\cdot)$, .... This method saves us from an explicit search for the solution $r_n(\phi_0,\phi_1)$ of the minimization problem in (7.1) for every $(\phi_0,\phi_1) \in \mathbf{C}_{n+1}$:

  *Step* B.0. Initialize $n=0$, $v_0(\cdot,\cdot) \equiv 0$.
  *Step* B.1. Find the region

$$(10.5) \quad B_n \triangleq \left\{(x,y) \in \mathbb{R}_+^2 : v_n(x,y) < -\frac{x}{\mu-1} - \frac{y}{\mu+1} + \frac{\lambda}{c\mu}\sqrt{2}\right\}, \qquad n \geq 1.$$

  *Step* B.2. Determine the continuation region $\mathbf{C}_{n+1} = S^{-1}(B_n)$ by the transformation of $B_n$ under $S^{-1}$.
  *Step* B.3. Calculate the value function $v_{n+1}(\cdot,\cdot)$ on $\mathbf{C}_{n+1}$ by using (9.2) and (9.7).
  *Step* B.4. Set $n$ to $n+1$ and go to Step B.1.

In fact, we can do much better than this. After $n \in \mathbb{N}$ iterations, we find both $v_n(\cdot,\cdot)$ and $V(\cdot,\cdot)$, $v_{n+1}(\cdot,\cdot)$, $v_{n+2}(\cdot,\cdot)$,... on the subset

$$(10.6) \quad Q_n \triangleq S^{-n}(\mathbf{\Gamma}_n) \cap \mathbf{C}_n \cap ([0, x_n] \times \mathbb{R}_+)$$

$$(10.7) \quad = \{(x,y) \in [0, x_n] \times \mathbb{R}_+ : S^{-n}[\gamma_n](x) \leq y < \gamma_n(x)\}, \qquad n \in \mathbb{N},$$

of $\mathbf{C}_{n+1}$ by Corollary 8.4. Therefore, we need to determine only the set

$$(10.8) \quad R_{n+1} \triangleq Q_{n+1} \setminus Q_n, \qquad n \in \mathbb{N} \qquad (R_1 \equiv Q_1),$$

in Step B.2, and calculate the value function $v_{n+1}(\cdot,\cdot)$ only on this set in Step B.3. By Lemma 8.5, this modified method calculates $V(\cdot,\cdot)$ [and all $V_n(\cdot,\cdot)$, $n \in \mathbb{N}$, simultaneously] on any given set $\mathbb{R}_+ \times (0, B)$, $B > 0$, in a finite number of iterations. We shall describe this modified method as Steps C.0–C.2 after Proposition 10.1 after establishing a few facts below. Several steps of the method are also illustrated in Figure 4.

Since $v_0 \equiv 0$, setting $n = 0$ in (10.4) gives a straight line; substituting $(x, y) = (x, S[\gamma_1](x))$ and comparing this with $S(\partial \Gamma_1)$ in (10.3) gives

$$S[\gamma_1](x) = -\frac{\mu+1}{\mu-1}x + \frac{\mu+1}{\mu} \cdot \frac{\lambda}{c}\sqrt{2},$$

(10.9)

$$x \in \mathrm{supp}(S[\gamma_1]) = \left[0, \frac{\mu-1}{\mu} \cdot \frac{\lambda}{c}\sqrt{2}\right],$$

and $\xi_1 = (\lambda/c)\sqrt{2}$. Using (8.7), we find

(10.10) $\gamma_1(x) = S^{-1}[S[\gamma_1]](x) = -x + \frac{\lambda}{c}\sqrt{2}, \qquad x \in [0, \xi_1] = \left[0, \frac{\lambda}{c}\sqrt{2}\right].$

The function $S^{-1}[\gamma_1](\cdot)$ is affine and intersects with $\gamma_1(\cdot)$ at $x_1 \equiv x_1(\gamma_1) = (\lambda/c)(\sqrt{2}/2)$; see (8.9). By Corollary 8.4 and Remark 8.6, the boundary of the stopping region on $[0, x_1]$ is

(10.11) $\gamma(x) = \gamma_1(x) = -x + \frac{\lambda}{c}\sqrt{2}, \qquad x \in [0, x_1] \equiv \left[0, \frac{\lambda\sqrt{2}}{c} \cdot \frac{1}{2}\right];$

see the inset in Figure 4(a). Hence, the boundaries of the optimal stopping region $\Gamma$ and the stopping regions $\Gamma_n$, $n \in \mathbb{N}$ stick on the upper half of the hypotenuse of the rectangular triangle $\{(x,y) \in \mathbb{R}_+^2 : g(x,y) \leq 0\} \equiv \mathrm{cl}(\mathbf{C}_0)$.

PROPOSITION 10.1. *Fix any $n \in \mathbb{N}$. The functions in $\mathcal{S}_n \triangleq (S^{-k}[\gamma_k])_{k=1}^n$ do not intersect inside the continuation region $\mathbf{C}_n = \{(x,y) \in \mathbb{R}_+^2 : y < \gamma_n(x)\}$. The function $S[\gamma_{n+1}](\cdot)$ intersects with each function in $\mathcal{S}_n \cup \{\gamma_n\}$ pairwise exactly once.*

The same conclusions hold when every $\gamma_k$, $k = 1, \ldots, n$, in the proposition is replaced with $\gamma$; this can be verified using the elementary properties of convex functions and the affine structure of the boundary function $\gamma(\cdot)$ in (10.11); see [2] for the details. Then the proof of Proposition 10.1 follows easily from Corollary 8.4. We are now ready to give a better version of method B in Section 10 to calculate each $v(\cdot, \cdot)$ and the boundary function $\gamma(\cdot)$. Recall that $S^{-n}[\gamma_n](\cdot)$ and $x_n(\gamma_n)$ are defined by (8.7) and (8.9). The steps C.1 and C.2 below are illustrated in Figure 4 for $n = 0$ and $n = 1$.

*Step* C.0. Initialize $n = 0$, $x_0 = 0$, $v_0(\cdot, \cdot) \equiv 0$ and the region $R_1$ as in (10.8).

*Step* C.1. Calculate the value function $V(\phi_0, \phi_1) = v_{n+1}(\phi_0, \phi_1)$ for every $(\phi_0, \phi_1) \in R_{n+1}$ using (9.2) and (9.7).

*Step* C.2. Set $n$ to $n+1$.

(i) Determine the set

(10.12) $\left\{(x,y) \in R_n : v_n(x,y) = -\frac{x}{\mu-1} - \frac{y}{\mu+1} + \frac{\lambda}{c\mu}\sqrt{2}\right\}$

of points in $R_n$ which satisfy (10.4). This is the intersection of the set in (10.3) and $R_n$. Namely, it is the section of the strictly decreasing, convex and continuous curve $x \mapsto S[\gamma_{n+1}](x)$ contained in $R_n$.

(ii) Find the subset of $R_n$ enclosed between the vertical $y$-axis and the curve in (10.12). This is the intersection $R_n \cap B_n$ of the sets $R_n$ and $B_n$ in (10.5).

(iii) Find the set $R_{n+1} = S^{-1}(R_n \cap B_n)$ in (10.8) by applying the transformation $S^{-1}(\cdot, \cdot)$ to the set found in (ii).

The region $R_{n+1}$ is enclosed between the $y$-axis from left, the $S^{-1}$-transformation of the curve in (10.12) from right. This right boundary of $R_{n+1}$ extends the boundary $\gamma(\cdot) \equiv \gamma_{n+1}(\cdot)$ from the previous iteration into the region $S^{-(n+1)}(\mathbf{\Gamma}) \equiv S^{-(n+1)}(\mathbf{\Gamma}_{n+1})$.

(iv) Go to Step C.1.

**11. The smoothness of the value functions and the stopping boundaries.** The general method described at the beginning of this section evaluates the integrals $Jv_n(\cdot, \phi_0, \phi_1)$ in (7.2) of the function $[g + \mu \cdot v_n \circ S](\cdot, \cdot)$ along the curves $(x(\cdot, \phi_0), y(\cdot, \phi_1))$ in $\mathbb{R}_+^2$ in order to calculate the value function $v_{n+1}(\phi_0, \phi_1)$ as in (7.1). For an accurate implementation of this method it may be useful to know how smooth the integrand, or essentially the value function $v_n(\cdot, \cdot)$, is.

The smoothness of the value function $V(\cdot, \cdot)$ may also allow us to formulate the original optimal stopping problem in (4.12) as a free-boundary problem. Then, in principle, we can calculate the value function $V(\cdot, \cdot)$ directly, by solving a partial differential equation, as the next proposition suggests.

PROPOSITION 11.1. *Suppose that there is a bounded and continuous function* $w : \mathbb{R}_+^2 \mapsto (-\infty, 0]$ *which is continuously differentiable on* $\mathbb{R}_+^2 \setminus \partial\mathbf{\Gamma}$, *and whose first-order derivatives are locally bounded near the boundary* $\partial\mathbf{\Gamma} = \{(x, \gamma(x)) : x \in [0, \xi]\}$. *Moreover,*

$$(11.1) \qquad (\widetilde{\mathcal{A}} - \lambda)w(x, y) + g(x, y) = 0, \qquad (x, y) \in \mathbf{C},$$

$$(11.2) \qquad w(x, y) = 0, \qquad (x, y) \in \mathbf{\Gamma},$$

$$(11.3) \qquad (\widetilde{\mathcal{A}} - \lambda)w(x, y) + g(x, y) > 0, \qquad (x, y) \in \mathbf{\Gamma} \setminus \partial\mathbf{\Gamma},$$

$$(11.4) \qquad w(x, y) < 0, \qquad (x, y) \in \mathbf{C},$$

*where* $\widetilde{\mathcal{A}}$ *is the infinitesimal generator in* (A.4) *of the process* $\widetilde{\mathbf{\Phi}}$ *acting on the continuously differentiable functions.*

*Suppose also that the sample-paths of the process* $\widetilde{\mathbf{\Phi}} = (\widetilde{\Phi}^{(0)}, \widetilde{\Phi}^{(1)})$ *spend zero time on the boundary* $\partial\mathbf{\Gamma}$ *almost surely, that is,*

$$(11.5) \qquad \mathbb{E}_0^{\phi_0, \phi_1}\left[\int_0^\infty \mathbf{1}_{\partial\mathbf{\Gamma}}(\widetilde{\mathbf{\Phi}}_t)\, dt\right] = 0, \qquad (\phi_0, \phi_1) \in \mathbb{R}_+^2.$$

*If the convex function $\gamma(\cdot)$ is also Lipschitz continuous on $[0,\xi]$, then $w(\cdot,\cdot) = V(\cdot,\cdot)$ on $\mathbb{R}_+^2$.*

PROOF. Similar to the proof of Theorem 10.4.1 in [14], page 215. □

Under certain conditions, we are able to show that the bounded, concave and continuous value functions $v_n(\cdot,\cdot)$, $n \in \mathbb{N}$, and $V(\cdot,\cdot)$ are continuously differentiable on $\mathbb{R}_+^2 \setminus \partial \mathbf{\Gamma}_{n+1}^x$ and $\mathbb{R}_+^2 \setminus \partial \mathbf{\Gamma}^x$, respectively, and are *not* differentiable on the exit boundaries $\partial \mathbf{\Gamma}_n^x$, and $\partial \mathbf{\Gamma}^x$ in (9.10), respectively. The exit boundaries $\partial \mathbf{\Gamma}_n^x$, $n \in \mathbb{N}$, and $\partial \mathbf{\Gamma}^x$, and the entrance boundaries $\partial \mathbf{\Gamma}_n^e$, and $\partial \mathbf{\Gamma}^e$, are connected subsets of $\mathbb{R}_+^2$, and we have

(11.6) $\quad \partial \mathbf{\Gamma}_n = \partial \mathbf{\Gamma}_n^x \cup \mathrm{cl}(\partial \mathbf{\Gamma}_n^e), \qquad n \in \mathbb{N}, \quad \text{and} \quad \partial \mathbf{\Gamma} = \partial \mathbf{\Gamma}^x \cup \mathrm{cl}(\partial \mathbf{\Gamma}^e).$

Moreover, the boundary functions $\gamma_n(\cdot)$ and $\gamma(\cdot)$ are continuously differentiable on their supports.

The hypotheses of Proposition 11.1 are satisfied with $w(\cdot,\cdot) \triangleq v(\cdot,\cdot)$ in (5.10). Thus, the function $v(\cdot,\cdot) \equiv V(\cdot,\cdot)$ may be obtained by solving the *variational inequalities* (11.1)–(11.4). This may be a challenging problem since, as we already pointed out above, the *smooth-fit principle* is guaranteed *not to hold* on some part of the free boundary. We shall not investigate the variational problem here, but give a concrete example with this interesting boundary behavior and describe our solution method for it.

The main result is Proposition 11.17 below, and is proven by induction. Here, we shall study the basis of the induction by breaking it down in several lemmas. The proof of the induction hypothesis is very similar, and we shall point out only the major differences after the proposition's statement on page 46.

Let us define the continuous mapping $G_n : \mathbb{R}_+^3 \mapsto \mathbb{R}$ by

(11.7) $$G_n(t,\phi_0,\phi_1) \triangleq [g + \mu \cdot v_n \circ S](x(t,\phi_0), y(t,\phi_1)),$$
$$(t,\phi_0,\phi_1) \in \mathbb{R}_+^3, n \geq 0.$$

Note that (9.2) gives

(11.8) $$v_{n+1}(\phi_0,\phi_1) = \int_0^{r_n(\phi_0,\phi_1)} e^{-(\lambda+\mu)t} G_n(t,\phi_0,\phi_1)\,dt,$$
$$(\phi_0,\phi_1) \in \mathbf{C}_{n+1}, n \geq 0.$$

Using (9.11), (9.12) and Lemmas 9.2 and 9.5, we obtain

(11.9) $\quad (x(t,\phi_0), y(t,\phi_1))|_{t=r_n(\phi_0,\phi_1)} \in \partial \mathbf{\Gamma}_{n+1}^e \subseteq \{(x, a_n(x)) : x \in [0, \alpha_n]\}$
$\quad \equiv \{(x,y) \in \mathbb{R}_+^2 : [g + \mu \cdot v_n \circ S](x,y) = 0\},$
$$(\phi_0,\phi_1) \in \mathbf{C}_{n+1} \cup \partial \mathbf{\Gamma}_{n+1}^x.$$

Therefore, for every $n \in \mathbb{N}$,

(11.10) $\quad 0 = G_n(t, \phi_0, \phi_1)|_{t=r_n(\phi_0, \phi_1)}, \qquad (\phi_0, \phi_1) \in \mathbf{C}_{n+1} \cup \partial \mathbf{\Gamma}_{n+1}^x.$

Under certain local smoothness and nondegeneracy conditions on the function $G_n(\cdot, \cdot, \cdot)$, the *implicit function theorems* guarantee that the equation $G_n(t, \phi_0, \phi_1) = 0$ determines $t = t_n(\phi_0, \phi_1)$ in an open neighborhood of every point $(r_n(\phi_0, \phi_1), \phi_0, \phi_1)$ in $\mathbb{R}_+ \times \mathbf{C}_{n+1}$, as a smooth function of the variables $(\phi_0, \phi_1)$. Since the continuation region $\mathbf{C}_{n+1}$ has compact closure, a finite subcovering of these open neighborhoods exists. Patching the solutions $t_n(\cdot, \cdot)$ in the finite subcovering gives the global solution, which is smooth and must coincide with $r_n(\cdot, \cdot)$ on $\mathbf{C}_{n+1}$.

In the remainder, we shall use the following version of the implicit function theorem; see [16], Chapter 14, and Conjecture 11.18.

THEOREM 11.2 (Implicit function theorem). *Let $A \subseteq \mathbb{R}^m$ be an open set, $F: A \mapsto \mathbb{R}$ be a continuously differentiable function, and $(\bar{t}, \bar{x}) \in \mathbb{R} \times \mathbb{R}^{m-1}$ be a point in $A$ such that*

$$F(\bar{t}, \bar{x}) = 0 \quad \text{and} \quad \left.\frac{\partial}{\partial t} F(t, x)\right|_{(t,x)=(\bar{t}, \bar{x})} \neq 0.$$

*Then there exist an open set $B \subseteq \mathbb{R}^{m-1}$ containing the point $\bar{x}$, and a unique continuously differentiable function $f: B \mapsto \mathbb{R}$, such that $\bar{t} = f(\bar{x})$ and $F(f(x), x) = 0$ for all $x \in B$.*

Since $v_0(\cdot, \cdot) \equiv 0$, we have

(11.11)
$$G_0(t, \phi_0, \phi_1) = x(t, \phi_0) + y(t, \phi_1) - \frac{\lambda}{c}\sqrt{2},$$
$$(t, \phi_0, \phi_1) \in \mathbb{R}_+ \times \mathbf{C}_1.$$

The function $G_0(\cdot, \cdot, \cdot)$ is continuously differentiable on $\mathbb{R}_+ \times \mathbf{C}_1$. By (6.8) in Section 6, its partial derivative

(11.12) $\qquad D_t G_0(t, \phi_0, \phi_1) = \frac{d}{dt}[x(t, \phi_0) + y(t, \phi_1)]$

with respect to $t$-variable may vanish at most once; if this happens, this derivative is strictly negative before, and strictly positive after, the derivative vanishes; otherwise, it is strictly positive everywhere (see also Figure 2). Namely, the function $t \mapsto G_0(t, \phi_0, \phi_1)$ has at most one local minimum for every $(\phi_0, \phi_1) \in \mathbb{R}_+^2$.

LEMMA 11.3. *Fix any $(\phi_0, \phi_1) \in \mathbb{R}_+^2$. The function $t \mapsto G_0(t, \phi_0, \phi_1)$ from $\mathbb{R}_+$ into $\mathbb{R}$ has at most one local minimum. It is strictly increasing if there is no local minimum. If there is a local minimum, then the function $G_0(\cdot, \phi_0, \phi_1)$ is strictly decreasing before the minimum and strictly increasing after the minimum.*

LEMMA 11.4. *The smallest minimizer $r_0(\phi_0, \phi_1)$ in (5.13) is continuously differentiable at every $(\phi_0, \phi_1) \in \mathbf{C}_1$.*

PROOF. The result will follow from Theorem 11.2 applied to the function $G_0(\cdot, \cdot, \cdot)$ on $\mathbb{R} \times \mathbf{C}_1$ at the point $(\bar{t}, \bar{x}) = (r_0(\phi_0, \phi_1), \phi_0, \phi_1) \in \mathbb{R} \times \mathbf{C}_1$. We only need to establish that

$$D_t G_0(t, \phi_0, \phi_1)|_{t=r_0(\phi_0,\phi_1)} \neq 0, \qquad (\phi_0, \phi_1) \in \mathbf{C}_1.$$

Let us fix $(\phi_0, \phi_1) \in \mathbf{C}_1$ and assume that $D_t G_0(r_0(\phi_0, \phi_1), \phi_0, \phi_1) = 0$. Then the function $t \mapsto G_0(t, \phi_0, \phi_1)$ is strictly decreasing on $t \in [0, r_0(\phi_0, \phi_1)]$ by Lemma 11.3 and

$$G_0(t, \phi_0, \phi_1) > G_0(r_0(\phi_0, \phi_1), \phi_0, \phi_1) = 0, \qquad t \in [0, r_0(\phi_0, \phi_1)).$$

Therefore, (11.8) implies that $v_1(\phi_0, \phi_1) > 0$. This contradicts our choice of $(\phi_0, \phi_1)$ in the continuation region $\mathbf{C}_1$ as well as the bound $v_1(\cdot, \cdot) \leq 0$. □

COROLLARY 11.5. *The value function $v_1(\phi_0, \phi_1)$ is continuously differentiable at every $(\phi_0, \phi_1) \in \mathbf{C}_1$. For every $(\phi_0, \phi_1) \in \mathbf{C}_1$,*

$$\begin{aligned}
D_{\phi_0} v_1(\phi_0, \phi_1) &= \int_0^{r_0(\phi_0,\phi_1)} e^{-(\lambda+\mu)t} D_{\phi_0} G_0(t, \phi_0, \phi_1)\, dt \\
&= \frac{1 - e^{-(\mu-1)r_0(\phi_0,\phi_1)}}{\mu - 1}, \\
D_{\phi_1} v_1(\phi_0, \phi_1) &= \int_0^{r_0(\phi_0,\phi_1)} e^{-(\lambda+\mu)t} D_{\phi_1} G_0(t, \phi_0, \phi_1)\, dt \\
&= \frac{1 - e^{-(\mu+1)r_0(\phi_0,\phi_1)}}{\mu + 1}.
\end{aligned} \tag{11.13}$$

PROOF. By (11.8) and Lemma 11.4, the value function $v_1(\cdot, \cdot)$ is continuously differentiable. Using (11.9) after applying the chain-rule to (11.8) with $n = 0$ gives the integrals in (11.13). These integrals can be calculated explicitly by using (4.7) or (4.8). □

COROLLARY 11.6. *The entrance boundary $\partial \mathbf{\Gamma}_1^e$ in (9.10) is connected. More precisely,*

$$(11.14) \quad \partial \mathbf{\Gamma}_1^e = \{(x, \gamma_1(x)) : x \in (\xi_1^e, \xi_1)\} \qquad \textit{for some } 0 \leq \xi_1^e < \xi_1,$$

*where $\xi_1$ is the same as in $[0, \xi_1] = \mathrm{supp}(\gamma_1)$, the support of the boundary function $\gamma_1(\cdot)$; see Proposition 8.1.*

COROLLARY 11.7. *The restriction of the boundary function $\gamma_1(\cdot)$ to the interval $(\xi_1^e, \xi_1)$ is continuously differentiable. In fact,*

$$\gamma_1(x) = a_0(x), x \in [\xi_1^e, \xi_1], \qquad [0, \xi_1] \equiv \mathrm{supp}(\gamma_1) = \mathrm{supp}(a_0) \equiv [0, \alpha_0],$$

*where*

$$(11.15) \qquad a_0(x) = \begin{cases} -x + \dfrac{\lambda}{c}\sqrt{2}, & x \in \left[0, \dfrac{\lambda}{c}\sqrt{2}\right), \\ 0, & \text{elsewhere,} \end{cases}$$

*is the continuously differentiable boundary function of the region $A_0 = \{(x,y) \in \mathbb{R}_+^2 : [g + \mu \cdot v_0 \circ S](x,y) < 0\}$ in (9.6).*

PROOF. The function $a_0(\cdot)$ in (11.15) is continuously differentiable on its support $\mathrm{supp}(a_0) = [0, \alpha_0]$, and the result follows from Lemma 9.5 and Corollary 11.6. □

The entrance boundary $\partial\mathbf{\Gamma}_1^e$ always exists. However, the exit boundary $\partial\mathbf{\Gamma}_1^x$ may not exist all the time. Next we shall identify the geometry of the exit boundary $\partial\mathbf{\Gamma}_1^x$ when it exists.

LEMMA 11.8. *For every $(\phi_0, \phi_1) \in \partial\mathbf{\Gamma}_1^x$, we have $[g + \mu \cdot v_0 \circ S](\phi_0, \phi_1) > 0$. Therefore, $\mathrm{cl}(\partial\mathbf{\Gamma}_1^e) \cap \partial\mathbf{\Gamma}_1^x = \varnothing$.*

PROOF. Suppose that $(\phi_0, \phi_1) \in \partial\mathbf{\Gamma}_1^x$. Let us assume that $[g + \mu \cdot v_0 \circ S](\phi_0, \phi_1) \leq 0$. Then $G_0(0, \phi_0, \phi_1) = [g + \mu \cdot v_0 \circ S](\phi_0, \phi_1) \leq 0 = G_0(r_0(\phi_0, \phi_1), \phi_0, \phi_1)$, and Lemma 11.3 implies

$$G_0(t, \phi_0, \phi_1) = [g + \mu \cdot v_0 \circ S](x(t, \phi_0), y(t, \phi_1)) < 0, \qquad t \in (0, r_0(\phi_0, \phi_1)).$$

This inequality and (11.8) for $n = 0$ give

$$v_1(\phi_0, \phi_1) = \int_0^{r_0(\phi_0, \phi_1)} e^{-(\lambda + \mu)t} G_0(t, \phi_0, \phi_1) \, dt < 0,$$

which contradicts $v_1(\phi_0, \phi_1) = 0$. This proves that $[g + \mu \cdot v_0 \circ S](\phi_0, \phi_1) > 0$ for every $(\phi_0, \phi_1) \in \partial\mathbf{\Gamma}_1^x$. Since the mapping $(x,y) \mapsto [g + \mu \cdot v_0 \circ S](x,y)$ is continuous, we have $[g + \mu \cdot v_0 \circ S](\phi_0, \phi_1) = 0$ for every $(\phi_0, \phi_1) \in \mathrm{cl}(\partial\mathbf{\Gamma}_1^e)$ by Lemmas 9.2 and 9.5. Therefore, $\mathrm{cl}(\partial\mathbf{\Gamma}_1^e) \cap \partial\mathbf{\Gamma}_1^x = \varnothing$. □

The next corollary is helpful in determining the point $(\xi_1^e, \gamma_1(\xi_1^e)) \equiv (\xi_1^e, a_0(\xi_1^e))$. The region $A_n$ was introduced in Section 9.

COROLLARY 11.9. *The parametric curve*

$$(11.16) \qquad \mathcal{C}_1 \triangleq \mathbb{R}_+^2 \cap \{(x(t, \xi_1^e), y(t, a_0(\xi_1^e))) : t \in \mathbb{R}\}$$

*is the smallest among all parametric curves* $\mathbb{R}_+^2 \cap \{(x(t,\phi_0), y(t,\phi_1)) : t \in \mathbb{R}\}$, $(\phi_0, \phi_1) \in \mathbb{R}_+^2$ *that majorize the boundary function* $a_0(\cdot)$ *of the region* $A_0 = \{(x,y) : [g + \mu \cdot v_0 \circ S](x,y) < 0\} = \{(x,y) \in \mathbb{R}_+^2 : y < a_0(x)\}$.

*The curve* $\mathcal{C}_1$ *and the boundary* $\partial A_0 = \{(x, a_0(x)) : x \in [0, \alpha_0]\}$ *touch exactly at* $(\xi_1^e, a_0(\xi_1^e)) \equiv (\xi_1^e, \gamma_1(\xi_1^e))$ *and nowhere else.*

PROOF. By Corollary 11.7 and Lemma 11.8, we have $(\xi_1^e, a_0(\xi_1^e)) = (\xi_1^e, \gamma_1(\xi_1^e)) \in \operatorname{cl}(\partial \mathbf{\Gamma}_1) \setminus \partial \mathbf{\Gamma}_1^e$. Therefore $(\xi_1^e, a_0(\xi_1^e)) \notin \partial \mathbf{\Gamma}_1^e \cup \partial \mathbf{\Gamma}_1^x$. Hence there exists some $\delta > 0$ such that

$$(11.17) \quad (x(t, \xi_1^e), y(t, a_0(\xi_1^e))) \in \mathbf{\Gamma}_1 \subseteq \mathbb{R}_+^2 \setminus A_0, \qquad t \in (-\delta, +\delta).$$

Recall from (9.18) that $\hat{r}(\xi_1^e, a_0(\xi_1^e))$ is the exit time of $t \mapsto (x(-t, \xi^e), y(-t, a_0(\xi_1^e)))$ from $\mathbb{R}_+^2$. Then the function

$$t \mapsto G_0(t, \xi_1^e, a_0(\xi_1^e)) \equiv [g + \mu \cdot v_0 \circ S](x(t, \xi_1^e), y(t, a_0(\xi_1^e))),$$
$$t \in [-\hat{r}(\xi_1^e, a_0(\xi_1^e)), \infty)$$

has a zero at $t = 0$ and is nonnegative for every $t \in (-\delta, \delta)$ by (11.17). Hence it has a local minimum at $(\xi_1^e, a_0(\xi_1^e))$. By Lemma 11.3, the function $G_0(t, \xi_1^e, a_0(\xi_1^e))$ is strictly positive for $t \neq 0$. Therefore, $y(t, a_0(\xi_1^e)) > a_0(x(t, \xi_1^e))$ for $t \neq 0$. $\square$

REMARK 11.10. Since $\mathcal{C}_1 \subset \mathbb{R}_+^2 \setminus A_0$, we have $Jv_0(t, \phi_0, \phi_1) > 0$ for $t > 0$ and $(\phi_0, \phi_1) \in \mathcal{C}_1$. This implies $v_1(\phi_0, \phi_1) = 0$ for $(\phi_0, \phi_1) \in \mathcal{C}_1$. Therefore, $\mathcal{C}_1 \subset \mathbf{\Gamma}_1$.

The curve $\mathcal{C}_1$ divides $\mathbb{R}_+^2$ into two components. Since the continuation region $\mathbf{C}_1$ is connected and contains $A_0$, the region $\mathbf{C}_1$ is contained in the (lower) component which lies between $\mathcal{C}_1$ and $x$-axis. Thus the boundary $\partial \mathbf{\Gamma}_1$ is completely below the curve $\mathcal{C}_1$, and they touch at the point $(\xi_1^e, a_0(\xi_1^e)) \equiv (\xi_1^e, \gamma_1(\xi_1^e))$. See Figure 5(a).

The next corollary shows that no point on the boundary $\{(x, a_0(x)) : x \in [0, \xi_1^e)\}$ over the interval $[0, \xi_1^e)$ of the region $A_0 = \{(x,y) \in \mathbb{R}_+^2 : [g + \mu \cdot v_0 \circ S](x,y) < 0\}$ is a boundary point for the stopping region $\mathbf{\Gamma}_1$.

COROLLARY 11.11. *For every* $x \in [0, \xi_1^e)$, *we have* $\gamma_1(x) > a_0(x)$ *and* $[g + \mu \cdot v_0 \circ S](x, \gamma_1(x)) > 0$.

PROOF. If $[0, \xi_1^e) = \varnothing$, there is nothing to prove. Otherwise, fix any $\phi_0 \in [0, \xi_1^e)$. Assume that $(\phi_0, a_0(\phi_0)) \in \partial \mathbf{\Gamma}_1$. By Corollary 11.6 and Lemma 11.8, we have $(\phi_0, a_0(\phi_0)) \notin \partial \mathbf{\Gamma}_1^e \cup \partial \mathbf{\Gamma}_1^x$. The same argument as in the proof of

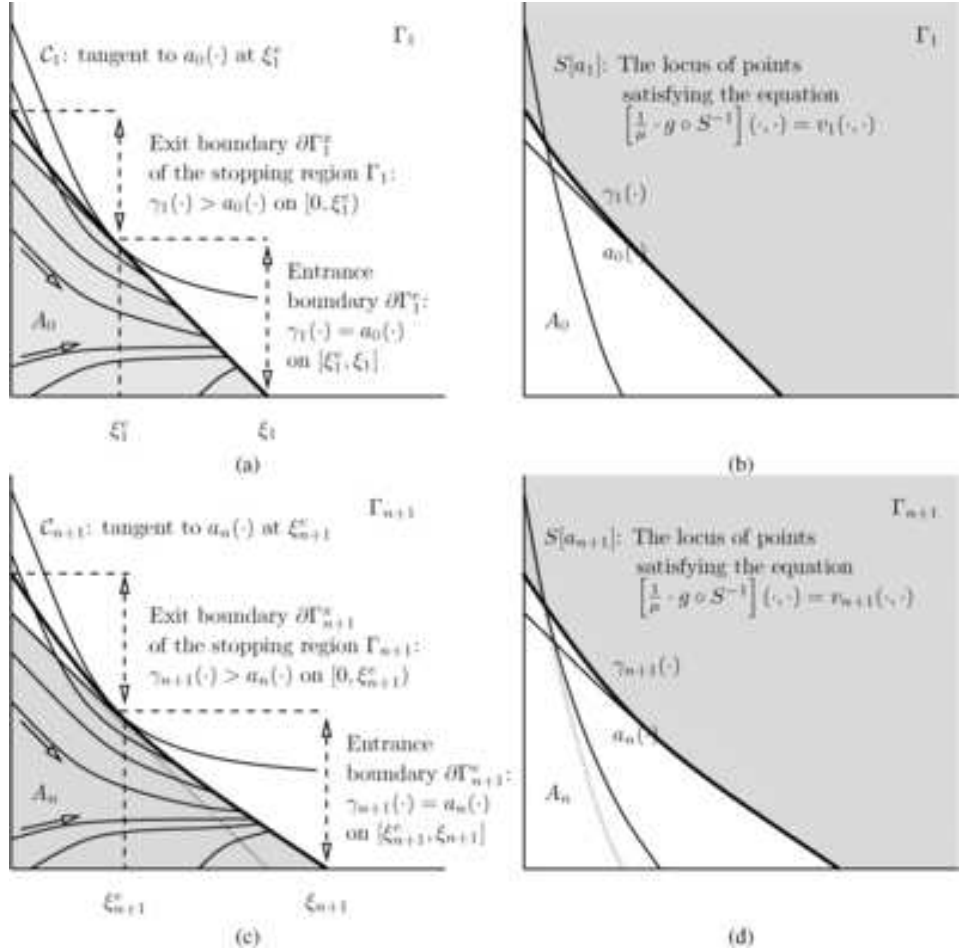

FIG. 5. (a) *and* (b) *illustrate Steps* D.1 *and* D.2 *of Method* D *for* $n = 0$, *and* (c) *and* (d) *for a general* $n$.

Corollary 11.9, with $(\phi_0, a_0(\phi_0))$ instead of $(\xi_1^e, a_0(\xi_1^e))$, gives that the parametric curve $\{(x(-t, \phi_0), y(-t, a_0(\phi_0))) : t \in [-\hat{r}(\phi_0, a_0(\phi_0)), \infty)\}$ is the smallest majorant of the boundary function $a_0(\cdot)$, and both curves touch at the point $(\phi_0, a_0(\phi_0))$. But this implies $\phi_0 = \xi_1^e$, a contradiction with our choice of $\phi_0$. □

COROLLARY 11.12. *If* $\phi_0 \in [0, \xi_1^e)$, *then* $(\phi_0, \gamma_1(\phi_0)) \in \partial \mathbf{\Gamma}_1^x$ *has an open neighborhood, on the intersection with the continuation region* $\mathbf{C}_1$ *of which the function* $r_0(\cdot, \cdot)$ *is bounded and bounded away from zero.*

*On the other hand, the function* $r_0(\cdot, \cdot)$ *is continuous on the entrance boundary* $\partial \mathbf{\Gamma}_1^e$: *for every* $(\phi_0, \phi_1) \in \partial \mathbf{\Gamma}_1^e$ *and every sequence* $\{(\phi_0^{(n)}, \phi_1^{(n)})\}_{n \in \mathbb{N}} \subseteq$ $\mathbf{C}_1$ *converging to the boundary point* $(\phi_0, \phi_1)$, *we have* $\lim_{n\to\infty} r_0(\phi_0^{(n)}, \phi_1^{(n)}) = 0$.

LEMMA 11.13. *If* $\xi_1^e = 0$, *then* $\partial \mathbf{\Gamma}_1 = \mathrm{cl}(\partial \mathbf{\Gamma}_1^e)$. *If* $\xi_1^e > 0$, *then the exit boundary* $\partial \mathbf{\Gamma}_1^x$ *is not empty, and* $\partial \mathbf{\Gamma}_1 = \partial \mathbf{\Gamma}_1^x \cup \mathrm{cl}(\partial \mathbf{\Gamma}_1^e)$.

If $\partial \mathbf{\Gamma}_1 \neq \mathrm{cl}(\partial \mathbf{\Gamma}_1^e)$, then $\xi_1^e > 0$, and the exit boundary $\partial \mathbf{\Gamma}_1^x$ is not empty by Lemma 11.13. The characterization of the exit boundary in Lemma 9.7 can be cast as

$$\mathbf{\Gamma}_1^x = \{(x(-t, \phi_0), y(-t, a_0(\phi_0)))|_{t=\hat{r}_0(\phi_0, a_0(\phi_0))} : \phi_0 \in (\xi_1^e, \xi_1^x]\}$$

for some $\xi_1^x \in (\xi_1^e, \xi_1)$. More precisely,

(11.18) $\qquad \xi_1^x = \inf\{\phi_0 \in [\xi_1^e, \xi_1] : \hat{r}_0(\phi_0, \gamma_1(x)) \leq \hat{r}(\phi_0, \gamma_1(x))\}.$

For every $\phi_0 \in (\xi_1^e, \xi_1^x]$, we have $\hat{r}_0(\phi_0, \gamma_1(\phi_0)) \leq \hat{r}(\phi_0, \gamma_1(\phi_0))$. See (9.18).

Our next result shows that the exit boundary $\partial \mathbf{\Gamma}_1^x = \{(\phi_0, \gamma_1(\phi_0)) : \phi_0 \in [0, \xi_1^e)\}$ is on a continuously differentiable curve, if it is not empty.

LEMMA 11.14. *The restriction of the boundary function* $\gamma_1(\cdot)$ *to the interval* $[0, \xi_1^e)$ *is continuously differentiable.*

If the value function $v_1(\cdot, \cdot)$ were continuously differentiable on the exit boundary $\partial \mathbf{\Gamma}_1^x$, then the result would follow from an application of the implicit function theorem to the identity $v_1(\phi_0, \phi_1) = 0$ near the point $(\phi_0, \phi_1) = (\phi_0, \gamma_1(\phi_0))$. Unfortunately, $v_1(\cdot, \cdot)$ is not differentiable on $\partial \mathbf{\Gamma}_1^x$; see Lemma 11.16. Therefore, we shall first extend the restriction to the set $\mathbf{C}_1 \cup \partial \mathbf{\Gamma}_1^x$ of the value function $v_1(\cdot, \cdot)$ to a new function $\tilde{v}_1(\cdot, \cdot)$ on an open set $B_1 \supset \mathbf{C}_1 \cup \partial \mathbf{\Gamma}_1^x$ such that $\tilde{v}_1(\cdot, \cdot)$ is continuously differentiable on $B_1$. We shall then use the identity $\tilde{v}_1(\phi_0, \gamma_1(\phi_1)) = 0$ as above.

LEMMA 11.15. *The boundary function* $\gamma_1(\cdot)$ *is continuously differentiable on the interior of its support* $[0, \xi_1]$.

The next result shows that the value function is not differentiable on the exit boundary $\partial \mathbf{\Gamma}_1^x$. In fact, as the proof reveals, the left and right partial derivatives are different along the exit boundary. Therefore, the *smooth-fit principle* does not apply to the value function $v_1(\cdot, \cdot)$ along (some part of) the boundary, if the exit boundary $\partial \mathbf{\Gamma}_1^x$ is not empty.

LEMMA 11.16. *The value function* $v_1(\cdot, \cdot)$ *is continuously differentiable on the entrance boundary* $\partial \mathbf{\Gamma}_1^e$, *but is* not *differentiable on the exit boundary* $\partial \mathbf{\Gamma}_1^x$.

The techniques used above in the analysis of the value function $v_1(\cdot,\cdot)$ and the boundary function $\gamma_1(\cdot)$ can be extended by induction to every function $v_n(\cdot,\cdot)$ and the boundary function $\gamma_n(\cdot)$, if the following are true for every $n \in \mathbb{N}$:

$\mathrm{A}_1(n)$: For every $(\phi_0, \phi_1) \in \mathbb{R}_+^2$, the function $t \mapsto G_n(t, \phi_0, \phi_1)$ in (11.7) from $\mathbb{R}_+$ into $\mathbb{R}$ has at most one local minimum. It is strictly increasing if there is no local minimum. If there is a local minimum, then the function $G_n(\cdot, \phi_0, \phi_1)$ is strictly decreasing before the minimum and strictly increasing after the minimum.

$\mathrm{A}_2(n)$: The function $(x,y) \mapsto [g + \mu \cdot v_n \circ S](x,y)$ is (continuously) differentiable on the entrance boundary $\partial \mathbf{\Gamma}^e_{n+1}$ of the stopping region $\mathbf{\Gamma}_{n+1} = \{(x,y) : v_{n+1}(x,y) = 0\}$.

PROPOSITION 11.17. *If $\mathrm{A}_1(k)$ and $\mathrm{A}_2(k)$ above are valid for every $0 \leq k \leq n$, then the following hold:*

1. *The value function $v_{n+1}(\cdot,\cdot)$ is continuously differentiable on $\mathbb{R}_+^2 \setminus \mathbf{\Gamma}^x_{n+1}$ everywhere except the exit boundary $\partial \mathbf{\Gamma}^x_{n+1}$. For every $(\phi_0, \phi_1) \in \mathbf{C}_{n+1}$,*

$$D_{\phi_0} v_{n+1}(\phi_0, \phi_1) = \int_0^{r_n(\phi_0, \phi_1)} e^{-(\lambda+\mu)t} D_{\phi_0} G_n(t, \phi_0, \phi_1)\,du$$

$$= \int_0^{r_n(\phi_0, \phi_1)} e^{-(\lambda+\mu)t} [1 + (\mu - 1) D_{\phi_0} v_n \circ S]$$

$$\times (x(u, \phi_0), y(u, \phi_1))\,du$$

*and*

$$D_{\phi_1} v_{n+1}(\phi_0, \phi_1) = \int_0^{r_n(\phi_0, \phi_1)} e^{-(\lambda+\mu)t} D_{\phi_1} G_n(t, \phi_0, \phi_1)\,du$$

$$= \int_0^{r_n(\phi_0, \phi_1)} e^{-(\lambda+\mu)t} [1 + (\mu + 1) D_{\phi_1} v_n \circ S]$$

$$\times (x(u, \phi_0), y(u, \phi_1))\,du.$$

2. *The entrance boundary $\partial \mathbf{\Gamma}^e_{n+1}$ is connected. More precisely,*

$$\partial \mathbf{\Gamma}^e_{n+1} = \{(x, a_n(x)) : x \in (\xi^e_{n+1}, \xi_{n+1})\} \qquad \textit{for some } \xi^e_{n+1} \in [0, \xi_{n+1}).$$

*The boundary function $a_n(\cdot)$ of the region $A_n = \{(x,y) \in \mathbb{R}_+^2 : [g + \mu \cdot v_n \circ S](x,y) < 0\}$ is continuously differentiable on $(\xi^e_{n+1}, \xi_{n+1})$. Therefore, the boundary function $\gamma_{n+1}(\cdot) \equiv a_n(\cdot)$ on $(\xi^e_{n+1}, \xi_{n+1})$ is continuously differentiable.*

3. *The parametric curve*

$$\mathcal{C}_{n+1} \triangleq \mathbb{R}_+^2 \cap \{(x(t, \xi^e_{n+1}), y(t, a_n(\xi^e_{n+1}))) : t \in \mathbb{R}\}$$

*is the smallest among all the parametric curves* $\mathbb{R}_+^2 \cap \{(x(t,\phi_0), y(t,\phi_1)) : t \in \mathbb{R}\}$, $(\phi_0, \phi_1) \in \mathbb{R}_+^2$ *that majorize the boundary function* $a_n(\cdot)$ *of the region* $A_n = \{(x,y) : [g + \mu \cdot v_n \circ S](x,y) < 0\} = \{(x,y) \in \mathbb{R}_+^2 : y < a_n(x)\}$.

*The curve* $\mathcal{C}_{n+1}$ *and the boundary* $\partial A_n = \{(x, a_n(x)) : x \in [0, \alpha_n]\}$ *touch exactly at* $(\xi_{n+1}^e, a_n(\xi_{n+1}^e)) \equiv (\xi_{n+1}^e, \gamma_{n+1}(\xi_{n+1}^e))$ *and nowhere else.*

4. *If* $\xi_{n+1}^e = 0$, *then* $\partial \mathbf{\Gamma}_{n+1} = \text{cl}(\partial \mathbf{\Gamma}_{n+1}^e)$. *If* $\xi_{n+1}^e > 0$, *then the exit boundary* $\partial \mathbf{\Gamma}_{n+1}^x$ *is not empty, and* $\partial \mathbf{\Gamma}_{n+1} = \partial \mathbf{\Gamma}_{n+1}^x \cup \text{cl}(\partial \mathbf{\Gamma}_{n+1}^e)$.
5. *The boundary function* $\gamma_{n+1}(\cdot)$ *is continuously differentiable on the interior of its support* $[0, \xi_{n+1}]$.

The proof of the proposition is by induction on $n \in \mathbb{N}_0$. The suppositions $A_1(0)$ and $A_2(0)$ are always valid; see Lemma 11.3, and note that $[g + \mu \cdot v_0 \circ S](\cdot, \cdot) \equiv g(\cdot, \cdot)$ is continuously differentiable everywhere. All of the claims are proved for the basis of the induction $n = 0$ before the statement of the proposition. For $n \geq 1$, the proofs are the same with obvious changes, with the exception of the differentiability of $a_n(\cdot)$ in part 2 of the proposition.

For $n = 0$, the differentiability of $a_0(x) = -x + (\lambda/c)\sqrt{2}$, $x \in (0, \xi_1)$ was obvious. For $n \geq 1$, the function $a_n(\cdot)$ is not available explicitly, only through

$$[g + \mu \cdot v_n \circ S](x, a_n(x)) = 0, \qquad x \in [0, \xi_{n+1}].$$

By $A_2(n)$, the function $[g + \mu \cdot v_n \circ S](\cdot, \cdot)$ is continuously differentiable on $\partial \mathbf{\Gamma}_{n+1}^e = \{(x, a_n(x)) : x \in (\xi_{n+1}^e, \xi_{n+1})\}$. Since $y \mapsto [g + \mu \cdot v_n \circ S](x, y)$ is strictly increasing for every $x \in \mathbb{R}_+$, we have

$$\frac{\partial}{\partial y}[g + \mu \cdot v_n \circ S](x, y)|_{(x,y)=(x, a_n(x))} > 0, \qquad x \in (\xi_{n+1}^e, \xi_{n+1}).$$

Thus, the function $a_n(\cdot)$ is continuously differentiable on $(\xi_{n+1}^e, \xi_{n+1})$ by the implicit function theorem.

11.1. *The interplay between the exit and entrance boundaries* We have been unable to identify fully all the cases where the hypotheses $A_1(n)$ and $A_2(n)$ on page 46 are satisfied for every $n \in \mathbb{N}$ (see, though, Section 11.2 for the important case of "large" disorder arrival rate $\lambda$, and Section 11.3 for another interesting example, where these hypotheses *are* satisfied). However, they are sufficient for Proposition 11.17 to hold, and Proposition 11.17 shows the crucial interplay between the exit and entrance boundaries. We would like to illustrate this interplay briefly; it may be very useful in designing efficient detection algorithms for general Poisson disorder problems. We suggest how the gap may be closed as an interesting open research problem.

In Section 9.2, we showed that both the value functions and the exit boundaries are determined by the entrance boundaries; see Lemma 9.7. More explicitly, once the entrance boundary $\partial \mathbf{\Gamma}_{n+1}^e$ has been obtained, one can calculate the value function $v_{n+1}(\cdot,\cdot)$ and the exit boundary $\partial\mathbf{\Gamma}_{n+1}^x$ by running backward in time the parametric curves $t \mapsto (x(t,\phi_0), y(t,\phi_1))$ from every point $(\phi_0,\phi_1)$ on the entrance boundary $\partial\mathbf{\Gamma}_{n+1}^e$ and by evaluating the explicit expressions of Lemma 9.7 along the way. On the other hand, the entrance boundary $\partial\mathbf{\Gamma}_{n+1}^e$ can be found when the value function $v_n(\cdot,\cdot)$ has already been calculated. Since $v_0 \equiv 0$ is readily available, the following iterative algorithm will give us every $v_n(\cdot,\cdot)$, $n \in \mathbb{N}_0$, and the boundary functions $\gamma_n(\cdot)$; see also Figure 5.

*Step* D.0. Initialize $n = 0$, $v_0(\cdot,\cdot) \equiv 0$ on $\mathbb{R}_+^2$. Let $a_0(\cdot)$ be the boundary function of the region $A_0 = \{(\phi_0,\phi_1) \in \mathbb{R}_+^2 : [g + \mu \cdot v_0 \circ S](\phi_0,\phi_1) < 0\}$; see (11.15).

*Step* D.1. There is unique number $\phi_0 = \xi_{n+1}^e$ in the bounded support $\phi_0 \in [0, \xi_{n+1}]$ of the function $a_n(\cdot)$ such that, for every small $\delta > 0$, we have

$$a_n(x(t,\phi_0)) \leq y(t, a_n(\phi_0)), \qquad t \in [0,\delta) \qquad \text{if } \phi_0 = 0, \quad \text{and}$$

$$t \in (-\delta, \delta) \qquad \text{if } \phi_0 > 0.$$

Equivalently, the parametric curve $\mathcal{C}_{n+1} \triangleq \mathbb{R}_+^2 \cap \{(x(t,\xi_{n+1}^e), y(t, a_n(\xi_{n+1}^e))) : t \in \mathbb{R}\}$ in (3) of Proposition 11.17 majorizes the boundary $\{(x, a_n(x)) : x \in \text{supp}(a_n)\}$ of the region $A_n = \{(x,y) \in \mathbb{R}_+^2 : [g + \mu \cdot v_n \circ S](x,y) < 0\}$ everywhere. The entrance boundary of the stopping region $\mathbf{\Gamma}_{n+1} = \{(\phi_0,\phi_1) \in \mathbb{R}_+^2 : v_{n+1}(\phi_0, \phi_1) = 0\}$ is given by $\partial\mathbf{\Gamma}_{n+1}^e = \{(\phi_0, a_n(\phi_0)) : \phi_0 \in (\xi_{n+1}^e, \xi_{n+1})\}$.

(i) Find the entrance boundary $\partial\mathbf{\Gamma}_{n+1}^e$.

(ii) For every $(\phi_0,\phi_1) \in \partial\mathbf{\Gamma}_{n+1}^e$, take the following steps to calculate the value function $v_{n+1}(\cdot,\cdot)$ on the continuation region $\mathbf{C}_{n+1}$ and the exit boundary $\partial\mathbf{\Gamma}_{n+1}^x$:

(a) Calculate $\hat{r}(\phi_0,\phi_1) \triangleq \inf\{t \geq 0 : (x(-t,\phi_0), y(-t,\phi_1)) \notin \mathbb{R}_+^2\}$.

(i) If $-Jv_n(-\hat{r}(\phi_0,\phi_1), \phi_0, \phi_1) < 0$, then set $\hat{r}_n(\phi_0,\phi_1) = \infty$. Otherwise, find

$$\hat{r}_n(\phi_0,\phi_1) \triangleq \inf\{t \in (0, \hat{r}(\phi_0,\phi_1)] : -Jv_n(-t, \phi_0, \phi_1) \geq 0\}$$

by a bisection search on $(0, \hat{r}(\phi_0,\phi_1)]$, and add the point

$$(x(-\hat{r}_n(\phi_0,\phi_1), \phi_0), y(-\hat{r}_n(\phi_0,\phi_1), \phi_1)) \in \partial\mathbf{\Gamma}_{n+1}^x$$

to the exit boundary.

(c) Calculate the value function

$$v_{n+1}(x(-t,\phi_0), y(-t,\phi_1)) = -e^{-(\lambda+\mu)t} Jv_n(-t, \phi_0, \phi_1)$$

along the curve $(x(-t,\phi_0), y(-t,\phi_1))$, $t \in (0, \hat{r}(\phi_0,\phi_1) \wedge \hat{r}_n(\phi_0,\phi_1)]$, until it either leaves $\mathbb{R}_+^2$ or hits the exit boundary $\partial\mathbf{\Gamma}_{n+1}^x$.

The union $\partial\mathbf{\Gamma}_{n+1}^x \cup \text{cl}(\partial\mathbf{\Gamma}_{n+1}^e) = \partial\mathbf{\Gamma}_{n+1}^x \cup \{(x, a_n(x)) : x \in [\xi_{n+1}^e, \xi_{n+1}]\}$ gives the boundary $\partial\mathbf{\Gamma}_{n+1} = \{(x, \gamma_{n+1}(x)) : x \in [0, \xi_{n+1}]\}$ and the boundary curve $\gamma_{n+1}(\cdot)$, which is strictly decreasing and convex on its support $[0, \xi_{n+1}]$.

(iii) Set $v_{n+1}(\cdot,\cdot) = 0$ on the stopping region $\mathbf{\Gamma}_{n+1} = \{(x,y) : y \geq \gamma_{n+1}(x)\}$.

*Step* D.2. Set $n$ to $n+1$. Determine the locus of the points $(\phi_0, \phi_1)$ in $\mathbb{R}_+^2$ that satisfy the equation

$$\left[\frac{1}{\mu} \cdot g \circ S^{-1}\right](\phi_0, \phi_1) = v_n(\phi_0, \phi_1).$$

This locus is the same as $\{(x, S[a_{n+1}](x)) : x \in \operatorname{supp}(S[a_{n+1}])\}$; see Notation 8.2. Shift it by the linear transformation $S^{-1}$ of (5.8) to obtain the boundary $\{(x, a_n(x)) : x \in \operatorname{supp}(a_n)\}$ of the region $A_n = \{(x,y) : [g + \mu \cdot v_n \circ S](x,y) < 0\}$. Go to Step D.1.

CONJECTURE 11.18. *The algorithm relies on only two results from Section* 11: (i) *the entrance boundary* $\partial\mathbf{\Gamma}_{n+1}^e$, $n \in \mathbb{N}_0$, *is connected, and* (ii) *the boundary* $\partial\mathbf{\Gamma}_{n+1}^x$, $n \in \mathbb{N}_0$, *is the disjoint union of the exit boundary* $\partial\mathbf{\Gamma}_{n+1}^x$ *and the closure of the entrance boundary* $\partial\mathbf{\Gamma}_{n+1}^e$. *Part* (ii) *was proved by using* (i) *and the first hypothesis* $A_1(n+1)$ *on page* 46; *see Lemma* 11.13. *We conjecture that the hypothesis* $A_1(n+1)$ *always holds for all* $n \in \mathbb{N}_0$.

*On the other hand, part* (i) *was proved by using the continuity of the mapping* $(\phi_0, \phi_1) \mapsto r_n(\phi_0, \phi_1)$ *on the connected continuation region* $\mathbf{C}_{n+1}$; *see Corollary* 11.6. *The continuity of the mapping* $r_n(\cdot,\cdot)$ *followed from its continuous differentiability on* $\mathbf{C}_{n+1}$, *which we proved by using the implicit function theorem (Theorem* 11.2*) under hypothesis* $A_2(n+1)$; *see Lemma* 11.4. *We conjecture that this mapping is always continuous on the continuation region* $\mathbf{C}_{n+1}$. *This may be proved directly by using a weaker version of the implicit function theorem (see, e.g.,* [12]*) or by using nonsmooth analysis (see, e.g.,* [7]*).*

11.2. *The regularity of the value functions and the optimal stopping boundaries when the disorder arrival rate* $\lambda$ *is "large."* One of the cases where both $A_1(n)$ and $A_2(n)$ on page 46 are satisfied for every $n \in \mathbb{N}$, is when the disorder arrival rate $\lambda$ is "large;" see Section 4.3 and Figure 1(a).

Suppose that $\lambda \geq [1 - (1+\mu)(c/2)]^+$. Then the curve $t \mapsto (x(t, \phi_0), y(t, \phi_1))$, and therefore the mapping $t \mapsto G_n(t, \phi_0, \phi_1)$, $t \in \mathbb{R}_+$, are strictly increasing for every $(\phi_0, \phi_1) \in \mathbb{R}_+^2$ and $n \in \mathbb{N}_0$; see Lemma 9.2. Hence, $A_1(n)$ always holds for every $n \in \mathbb{N}_0$.

For the same reason, all of the exit boundaries $\partial\mathbf{\Gamma}_n^x$, $n \in \mathbb{N}$, are empty; see Section 10. Since $\partial\mathbf{\Gamma}_1^x$ is empty, the value function $v_1(\cdot,\cdot)$ is continuously differentiable everywhere. Therefore, $A_2(1)$ holds. Then Proposition 11.17 implies that $v_2(\cdot,\cdot)$ is continuously differentiable everywhere since $\partial\mathbf{\Gamma}_2^x$ is empty. Therefore, $A_2(2)$ holds, and so on.

COROLLARY 11.19 ("Large" disorder arrival rate: smooth solutions of reference optimal stopping problems). *Suppose that* $\lambda \geq [1 - (1+\mu)(c/2)]^+$.

*Then* $A_1(n)$ *and* $A_2(n)$ *hold for every* $n \in \mathbb{N}_0$, *and Proposition* 11.17 *applies. Particularly, for every* $n \in \mathbb{N}_0$:

1. *The value function* $v_{n+1}(\cdot, \cdot)$ *is continuously differentiable everywhere.*
2. *The exit boundary* $\partial \mathbf{\Gamma}_{n+1}^x$ *is empty, and* $\partial \mathbf{\Gamma}_{n+1} = \mathrm{cl}(\partial \mathbf{\Gamma}_{n+1}^e)$.
3. *The boundary function* $\gamma_{n+1}(\cdot)$ *is continuously differentiable on the interior of its support* $[0, \xi_{n+1}]$. *Thus, the function* $\gamma_{n+1}(\cdot)$ *coincides with the boundary*

$$a_0(x) = -x + \frac{\lambda}{c}\sqrt{2} \qquad \text{of the region } A_0 \text{ on the interval } \left[0, \frac{\lambda\sqrt{2}}{2c}\right],$$

*and fits smoothly to this function at the right end-point of the same interval.*

The last part of (3) in the corollary follows from (10.11) in Section 10 and Proposition 11.17. Recall also from Remark 8.6 that, if the disorder arrival rate $\lambda$ is "large," then there is an increasing sequence of sets $\mathbb{R}_+ \times [B_n, \infty)$ whose limit is $\mathbb{R}_+ \times (0, \infty)$, and $v(\cdot, \cdot) = v_n(\cdot, \cdot)$ on $\mathbb{R}_+ \times [B_n, \infty)$ for every $n \in \mathbb{N}$. Therefore, Corollary 11.19 implies immediately that the value function $v(\cdot, \cdot)$ and the boundary function $\gamma(\cdot)$ are continuously differentiable on $\mathbb{R}_+ \times (0, \infty)$ and on the interior of the support $[0, \xi]$ of the function $\gamma(\cdot)$, respectively.

To prove that $v(\cdot, \cdot)$ is continuously differentiable on $(0, \infty) \times \{0\}$, we shall use again the implicit function theorem. By Proposition 5.6 and Remark 5.10,

$$v(\phi_0, 0) = Jv(r(\phi_0, 0), \phi_0, 0) = \int_0^{r(\phi_0, 0)} e^{-(\lambda+\mu)t} G(t, \phi_0, 0)\, dt, \qquad \phi_0 \in \mathbb{R}_+.$$

The function $(t, \phi_0) \mapsto G(t, \phi_0, 0) \triangleq [g + \mu \cdot v \circ S](x(t, \phi_0), y(t, 0))$ is continuously differentiable on $(0, \infty) \times (0, \infty)$ since $v(\cdot, \cdot)$ is continuously differentiable on $\mathbb{R}_+ \times (0, \infty)$ and $(x(t, \phi_0), y(t, 0)) \in (0, \infty) \times (0, \infty)$ for every $t > 0$. Moreover, the partial derivative $(t, \phi_0) \mapsto D_{\phi_0} G(t, \phi_0, 0)$ is locally bounded on $(0, \infty) \times (0, \infty)$ by Corollary 5.4. Therefore, the function $(t, \phi_0) \mapsto Jv(t, \phi_0, 0)$ is continuously differentiable on $(0, \infty) \times (0, \infty)$ and

$$D_{\phi_0} Jv(t, \phi_0, 0) = \int_0^t e^{-(\lambda+\mu)u} D_{\phi_0} G(t, \phi_0, 0)\, du$$

$$= \int_0^t e^{-(\mu+1)u}[1 + (\mu - 1)D_{\phi_0} v \circ S](x(u, \phi_0), y(u, 0))\, du,$$

$$(t, \phi_0) \in \mathbb{R}_+ \times (0, \infty).$$

Since $v(\phi_0, 0) \equiv 0$ for every $\phi_0 \in [\xi, \infty)$, it is continuously differentiable on $(\xi, \infty)$. To show that it is differentiable on $(0, \xi)$, it is enough to prove that

ADAPTIVE POISSON DISORDER PROBLEM 51the mapping $\phi_0 \mapsto r(\phi_0, 0)$ from $(0, \xi)$ to $\mathbb{R}_+$ is continuously differentiable. Observe that, if we define

$$F(t, \phi_0) \triangleq \gamma(x(t, \phi_0)) - y(t, 0), \qquad (t, \phi_0) \in \mathbb{R}_+^2,$$

then $F(r(\phi_0, 0), \phi_0) = 0$ for every $\phi_0 \in [0, \xi]$. For every $\phi_0 \in (0, \xi)$, the function $F(\cdot, \cdot)$ is continuously differentiable in some neighborhood of $(r_0(\phi_0, 0), \phi_0)$ since $x(r(\phi_0, 0), \phi_0) \in (0, \xi)$ and $\gamma(\cdot)$ is continuously differentiable on $[0, \xi)$. Moreover,

$$D_t F(t, \phi_0) = \gamma'(x(t, \phi_0)) D_t x(t, \phi_0) - D_t y(t, 0) < 0$$

at every $(t, \phi_0) \in \mathbb{R}_+^2$ where $D_t F(t, \phi_0)$ exists [because since $\gamma(\cdot)$ is decreasing, $t \mapsto x(t, 0)$ and $t \mapsto x(t, \phi_0)$ are strictly increasing]. Then the implicit function theorem implies that $\phi_0 \mapsto r(\phi_0, 0)$, and therefore, $\phi_0 \mapsto v(\phi_0, 0) = Jv(r(\phi_0, 0), \phi_0, 0)$, is continuously differentiable at $\phi_0 \in (0, \xi)$. An argument similar to that in Corollary 11.12 shows that $\phi_0 \mapsto r(\phi_0, 0)$ is continuous at $\phi_0 = \xi$ and $\lim_{\phi_0 \uparrow \xi} r(\phi_0, 0) = 0$. By the Leibniz rule (see, e.g., [16], Theorem 11.1, page 286), the limit of the derivative

$$D_{\phi_0} v(\phi_0, 0) = D_t Jv(r(\phi_0, 0), \phi_0, 0) + D_{\phi_0} Jv(r(\phi_0, 0), \phi_0, 0)$$
$$= \int_0^{r(\phi_0, 0)} e^{-(\mu+1)u} [1 + (\mu - 1) D_{\phi_0} v \circ S](x(u, \phi_0), y(u, 0)) \, du,$$
$$\phi_0 \in (0, \xi),$$

of the value function $v(\cdot, \cdot)$ at $(\phi_0, 0)$, as $\phi_0$ increases to $\xi$, equals zero. Recall that, since $r(\phi_0, 0) > 0$ for every $\phi_0 \in [0, \xi)$, the derivative of $t \mapsto Jv(t, \phi_0, 0)$ on the right-hand side vanishes at its minimizer $t = r(\phi_0, 0)$. Thus, the left and right derivatives of the concave function $\phi_0 \mapsto v(\phi_0, 0)$ at $\phi_0 = \xi$, are equal:

$$D_{\phi_0}^- v(\xi, 0) = \lim_{\phi_0 \uparrow \xi} D_{\phi_0}^- v(\phi_0, 0) = \lim_{\phi_0 \uparrow \xi} D_{\phi_0} v(\phi_0, 0) = 0 = D_{\phi_0}^+ v(\xi, 0).$$

This shows that the function $(0, \infty) \ni \phi_0 \mapsto v(\phi_0, \cdot)$ is continuously differentiable. Hence the value function $v(\cdot, \cdot)$ is continuously differentiable on $\mathbb{R}_+ \times \{0\}$.

COROLLARY 11.20 ("Large" disorder arrival rate: smooth solution of the main optimal stopping problem). *Suppose that* $\lambda \geq [1 - (1 + \mu)(c/2)]^+$. *Then:*

1. *the value function $v(\cdot, \cdot)$ is continuously differentiable everywhere;*
2. *the boundary function $\gamma(\cdot)$ is continuously differentiable on the interior of its support $[0, \xi]$, coincides with the boundary function*

$$a_0(x) = -x + \frac{\lambda}{c}\sqrt{2} \qquad \textit{of the region } A_0 \textit{ on the interval } \left[0, \frac{\lambda\sqrt{2}}{2c}\right],$$

*and fits smoothly to this function at the right end-point of the same interval;*

3. *the function $v(\cdot,\cdot)$ is the solution of the variational inequalities* (11.1)–(11.4).

COROLLARY 11.21. *If $\lambda \geq [1 - (1+\mu)(c/2)]^+$, then both the sequence $\{v_n\}_{n \in \mathbb{N}}$ of the value functions and the sequences of their partial derivatives $\{D_{\phi_0} v_n\}_{n \in \mathbb{N}}$ and $\{D_{\phi_1} v_n\}_{n \in \mathbb{N}}$ converge uniformly to the value function $v$ and its partial derivatives $D_{\phi_0} v$ and $D_{\phi_1} v$, respectively.*

PROOF. This follows immediately from Theorem 25.7 in [17], page 248.
□

The results obtained in Section 9 for the functions $v_n$, $n \in \mathbb{N}$ can be extended easily to the value function $v(\cdot,\cdot)$, $n \in \mathbb{N}$. As in Lemma 9.2, we have

$$A \triangleq \{(x,y) \in \mathbb{R}_+^2 : [g + \mu \cdot v \circ S](x,y) < 0\} = \{(x,y) \in \mathbb{R}_+^2 : y < a(x)\},$$

$$\{(x,y) \in \mathbb{R}_+^2 : [g + \mu \cdot v \circ S](x,y) = 0\} = \{(x,a(x)) : x \in [0,\alpha]\}$$

for some decreasing function $a : \mathbb{R}_+ \mapsto \mathbb{R}_+$ which is strictly decreasing on its finite support $[0,\alpha]$. We have $A \subseteq \mathbf{C}$, and equality holds if $\lambda \geq [1 - (1+\mu)(c/2)]^+$ since the parametric curves $t \mapsto (x(t,\phi_0), y(t,\phi_1))$ increase and do not come back to the region $A$ after they leave; see Section 10. Therefore, $\gamma(\cdot) \equiv a(\cdot)$ and

(11.19) $\quad [g + \mu \cdot v \circ S](x,y) > 0, \quad (x,y) \in \mathbf{\Gamma} \setminus \partial \mathbf{\Gamma}.$

PROOF OF COROLLARY 11.20. Only (3) remains to be proven. The function $v : \mathbb{R}_+^2 \mapsto (-\infty, 0]$ is bounded and continuously differentiable. By the definition of the continuation region $\mathbf{C} = \{(x,y) \in \mathbb{R}_+^2 : v(x,y) < 0\}$ and the stopping region $\mathbf{\Gamma} = \mathbb{R}_+^2 \setminus \mathbf{C}$, the (in)equalities (11.2) and (11.4) are satisfied. On the other hand, (11.19) implies

$$[(\widetilde{\mathcal{A}} - \lambda)v + g](\phi_0, \phi_1) = [g + \mu \cdot v \circ S](\phi_0, \phi_1), \quad (\phi_0, \phi_1) \in \mathbf{\Gamma},$$

is strictly positive for every $(\phi_0, \phi_1) \in \mathbf{\Gamma} \setminus \partial \mathbf{\Gamma}$, that is, (11.3) is also satisfied. On the other hand, $[(\widetilde{\mathcal{A}} - \lambda)v + g](\phi_0, \phi_1)$ equals

$$D_{\phi_0} v(\phi_0, \phi_1) \left[(\lambda+1)\phi_0 + \frac{\lambda(1-m)}{\sqrt{2}}\right] + D_{\phi_1} v(\phi_0, \phi_1) \left[(\lambda-1)\phi_1 + \frac{\lambda(1+m)}{\sqrt{2}}\right]$$

$$+ \mu \left[v\left(\left(1 - \frac{1}{\mu}\right)\phi_0, \left(1 + \frac{1}{\mu}\right)\phi_1\right) - v(\phi_0, \phi_1)\right] - \lambda v(\phi_0, \phi_1) + g(\phi_0, \phi_1)$$

$$= D_{\phi_0} v(\phi_0, \phi_1) \cdot D_t x(0, \phi_0)$$

$$+ D_{\phi_1}v(\phi_0,\phi_1) \cdot D_t y(0,\phi_1) - (\lambda+\mu)v(\phi_0,\phi_1)$$
$$+ [g + \mu \cdot v \circ S](\phi_0,\phi_1)$$
$$= \frac{\partial}{\partial t}\left[e^{-(\lambda+\mu)t}v(x(t,\phi_0),y(t,\phi_1))\right.$$
$$\left.+ \int_0^t e^{-(\lambda+\mu)t}[g+\mu v \circ S](x(u,\phi_0),y(u,\phi_1))\,du\right]\Big|_{t=0}$$
$$= \frac{\partial}{\partial t}[e^{-(\lambda+\mu)t}v(x(t,\phi_0),y(t,\phi_1)) + Jv(t,\phi_0,\phi_1)]|_{t=0}, \qquad (\phi_0,\phi_1) \in \mathbf{C}.$$

Observe that the expression in the square brackets of the last equation equals $v(\phi_0,\phi_1)$ for every sufficiently small $t>0$, by (5.19) in Remark 5.10. Therefore, the derivative above equals zero, and (11.1) holds. This completes the proof that the function $v(\cdot,\cdot)$ satisfies the variational inequalities (11.1)–(11.4).

The boundary function $\gamma(\cdot)$ is strictly decreasing on its support. The process $\widetilde{\mathbf{\Phi}}$ can have at most countably many jumps, and its paths are strictly increasing between the jumps. Therefore, the time that the process $\widetilde{\mathbf{\Phi}}$ spends on the boundary $\partial\mathbf{\Gamma} = \{(x,\gamma(x)): x \in [0,\xi]\}$ equals zero almost surely. Since the derivative of the convex boundary curve $0 \geq \gamma'(x) \geq \gamma'(0+) = a_0(0+) = -1$ is bounded on $x \in (0,\xi)$, the curve is Lipschitz continuous on its support. □

Finally, Corollary 11.20 also shows that for every $\lambda \geq [1-(1+m)(c/2)]^+$, the smooth restrictions of value function $v_{n+1}(\cdot,\cdot)$ to the continuation region $\mathbf{C}_{n+1}$ and to the stopping region $\mathbf{\Gamma}_{n+1}$ fit to each other smoothly across the smooth boundary $\partial\mathbf{\Gamma}_{n+1} = \{(x,\gamma_{n+1}(x)): x \in [0,\xi_{n+1}]\}$.

However, if $0 < \lambda < 1-(1+m)(c/2)$ is small, then the corresponding value function does not have to have the same *smooth-fit property*.

11.3. *Failure of the smooth-fit principle*: *a concrete example.* Here we shall give a concrete example of a case where the value function fits smoothly across the entrance boundary, but fails to fit smoothly across the exit boundary of the optimal stopping region; see Figure 6(d).

Suppose that the disorder arrival rate $\lambda$, the pre-disorder arrival rate $\mu$ of the observations, the detection delay cost $c$ per unit time, and the expectation $m = \mathbb{E}_0[\Lambda - \mu]$ of the difference $\Lambda - \mu$ between the arrival rates of the observations after and before the disorder, are chosen so that

$$0 < \lambda < 1 - (1+m)(c/2),$$
$$\frac{\mu+1}{\mu}\phi_d > \bar{\phi}_1,$$

(11.20)

$$y < S[a_0](x) = \frac{\mu+1}{\mu} a_0\left(\frac{\mu}{\mu-1}x\right),$$
$$(x,y) \in \{(\phi_0^*, \phi_1^*), (0, \bar{\phi}_1)\}.$$

Here $\phi_d > 0$ is the mean-reversion level in (4.14) of $y \mapsto y(t, \phi_1)$ for every initial condition $\phi_1 \in \mathbb{R}_+$; see Section 4.4. The point

(11.21) $$(\phi_0^*, \phi_1^*) = \left(\frac{\lambda}{\sqrt{2}}\left(\frac{1-\lambda}{c} - 1\right), \frac{\lambda}{\sqrt{2}}\left(\frac{1+\lambda}{c} + 1\right)\right)$$

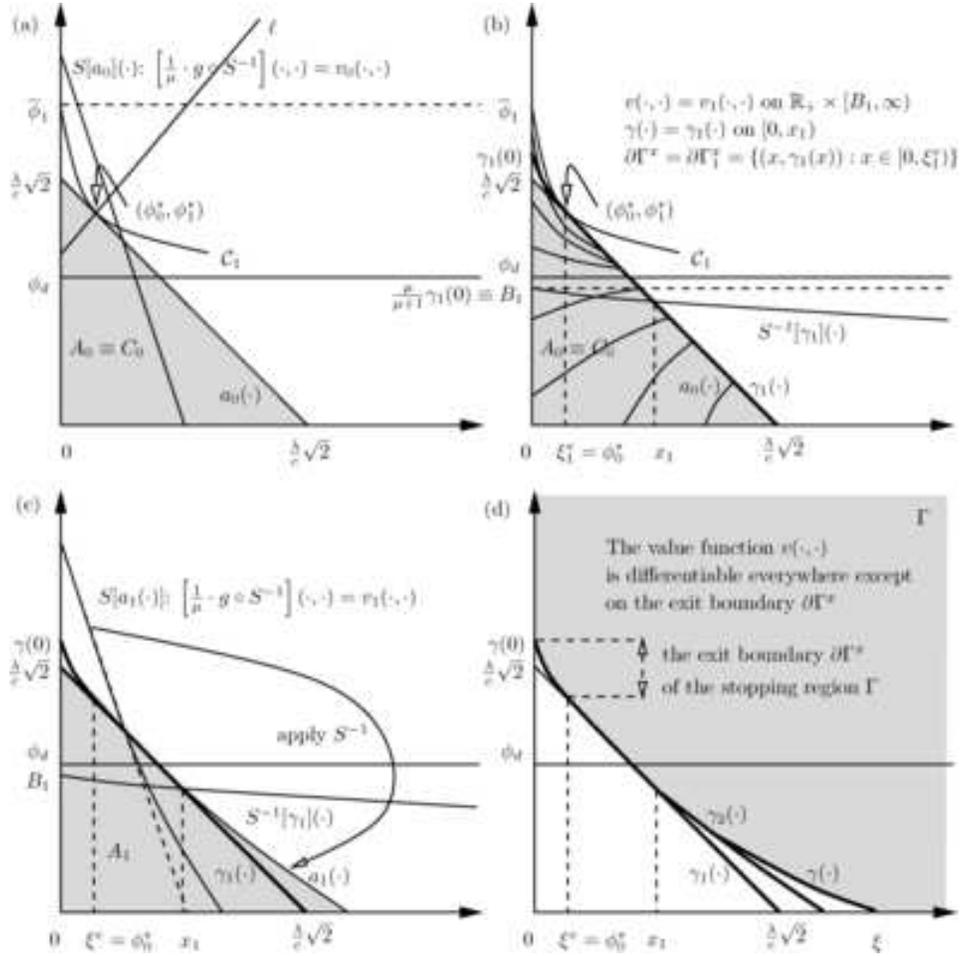

FIG. 6. (a) shows the location of points $(\phi_0^*, \phi_1^*)$ and $(0, \bar{\phi}_1)$ and the line described by the function $S[a_0](\cdot)$. In (b) and (c), we recall how to find the function $v_1(\cdot)$ and region $A_1$, respectively; compare with Figure 4. The boundary function $\gamma(\cdot)$ of the stopping region $\Gamma = \{(x,y) \in \mathbb{R}_+^2 : v(x,y) = 0\}$ is sketched in (d).

is the intersection of the straight lines $\ell$ in (6.7) and $y = a_0(x)$. Recall from (11.15) that $a_0(\cdot)$ is the boundary function of the region $A_0 = \{(x, y) \in \mathbb{R}_+^2 : g(x, y) < 0\} \equiv \{(x, y) \in \mathbb{R}_+^2 : y < a_0(x)\}$. For every initial point $(\phi_0, \phi_1)$ in $\mathbb{R}_+^2$, the sum $t \mapsto x(t, \phi_0) + y(t, \phi_1)$, $t \in \mathbb{R}_+$, of the coordinates of the parametric curve $t \mapsto (x(t, \phi_0), y(t, \phi_1))$, $t \in \mathbb{R}_+$, strictly decreases before the parametric curve meets the line $\ell$, and strictly increases thereafter; see Lemma 11.3 and (6.8). Finally, the point $(0, \bar{\phi}_1)$ with

$$
\begin{aligned}
\bar{\phi}_1 = &-\frac{\lambda(1+m)}{\sqrt{2}(\lambda-1)} \\
&+ \left[\phi_1^* + \frac{\lambda(1+m)}{\sqrt{2}(\lambda-1)}\right]\left[1 + \phi_0^* \frac{\sqrt{2}(\lambda+1)}{\lambda(1-m)}\right]^{-(\lambda-1)/(\lambda+1)}
\end{aligned}
\tag{11.22}
$$

is the initial point on the $y$-axis of the parametric curve $t \mapsto (x(t, 0), y(t, \bar{\phi}_1))$, $t \in \mathbb{R}_+$ which passes through the point $(\phi_0^*, \phi_1^*)$ in (11.21). The coordinate $\bar{\phi}_1$ in (11.22) is found by substituting the solution of $x(t^*, 0) = \phi_0^*$ for $t^*$ into the equation $y(t^*, \bar{\phi}_1) = \phi_1^*$ and solving the latter for $\bar{\phi}_1$; see also Figure 6(a).

Let us show that, under the conditions in (11.20), the "closedness" property in (8.4) holds. By Lemma 11.3 and (6.8), the curve $\mathcal{C}_1$ in Corollary 11.9 becomes

$$\mathcal{C}_1 = \{(x(t, 0), y(t, \bar{\phi}_1)) : t \in \mathbb{R}_+\} \equiv \mathbb{R}_+^2 \cap \{(x(t, \phi_0^*), y(t, \phi_1^*)) : t \in \mathbb{R}\};$$

it is tangent to the broken line $\{(x, a_0(x)) : x \in \mathbb{R}_+\}$ at the point $(\phi_0^*, \phi_1^*)$. Therefore, $\xi_1^e = \phi_0^*$ by the same corollary, and the entrance boundary of the stopping region $\mathbf{\Gamma}_1 = \{(x, y) : v_1(x, y) = 0\}$ is $\partial \mathbf{\Gamma}_1^e = \{(x, a_0(x)) : x \in (\phi_0^*, (\lambda/c)\sqrt{2})\}$ by Corollary 11.6. Moreover, the boundary function $\gamma_1(\cdot)$ of the region $\mathbf{\Gamma}_1 = \{(x, \gamma_1(x)) : \gamma_1(x) \leq y\}$ is supported on $[0, (\lambda/c)\sqrt{2}]$ and satisfies

$$
\begin{aligned}
\gamma_1(x) &= a_0(x), x \in \left[\phi_0^*, \frac{\lambda}{c}\sqrt{2}\right] \quad \text{and} \\
\gamma_1(x) &< y(0, \bar{\phi}_1) = \bar{\phi}_1, \qquad x \in [0, \phi_0^*).
\end{aligned}
\tag{11.23}
$$

The equality follows from Corollary 11.7, and the inequality from Remark 11.10 and the fact that the parametric curve $\mathcal{C}_1$ is decreasing. One can easily see from (11.23) and the second inequality in (11.20) that

$$
\begin{aligned}
(0, \phi_d) \in \mathbb{R}_+ \times [\phi_d, \infty) &\subset S^{-1}(\mathbf{\Gamma}_1) = \{(x, S^{-1}[\gamma_1](x)) : x \in \mathbb{R}_+\}, \\
\phi_d &\geq S^{-1}[\gamma_1](0) = \frac{\mu}{\mu+1}\gamma_1(0).
\end{aligned}
\tag{11.24}
$$

The restrictions of the value functions $v(\cdot, \cdot)$ and $v_1(\cdot, \cdot)$, and therefore those of the boundaries $\partial \mathbf{\Gamma}$ and $\partial \mathbf{\Gamma}_1$, coincide on $\mathbb{R}_+ \times [\phi_d, \infty)$. First, observe

that

$$\left\{\begin{array}{l} x(t,\phi_0) + y(t,\phi_1) \geq x(t,0) + y(t,\bar{\phi}_1) \geq \dfrac{\lambda}{c}\sqrt{2} \\ \text{that is, } (x(t,\phi_0), y(t,\phi_1)) \notin \mathbf{C}_0, \qquad t \in \mathbb{R}_+ \end{array}\right\},$$

$$(\phi_0, \phi_1) \in \mathbb{R}_+ \times [\bar{\phi}_1, \infty),$$

where $\mathbf{C}_0 = \{(x,y) \in \mathbb{R}_+^2 : g(x,y) < 0\}$ is as in (4.13) and coincides with $A_0$. By the second inequality in (11.20) and the properties of the parametric curves $t \mapsto (x(t,\phi_0), y(t,\phi_1))$, $t \in \mathbb{R}_+$ (see Section 4.2), we have

$$(\phi_0, \phi_1) \in \mathbb{R}_+ \times [\phi_d, \infty)$$
$$\implies \left\{\begin{array}{l} S(\phi_0,\phi_1) \in \mathbb{R}_+ \times [\bar{\phi}_1, \infty) \subset \mathbb{R}_+ \times [\phi_d, \infty) \\ (x(t,\phi_0), y(t,\phi_1)) \in \mathbb{R}_+ \times [\phi_d, \infty), \qquad t \in \mathbb{R}_+ \end{array}\right\}.$$

Using the last two displayed equations gives that, if the initial state $\widetilde{\boldsymbol{\Phi}}_0$ on a sample-path of the sufficient statistic $\widetilde{\boldsymbol{\Phi}} = (\widetilde{\Phi}^{(0)}, \widetilde{\Phi}^{(1)})$ is in $\mathbb{R}_+ \times [\phi_d, \infty)$, then the sample-path stays in the region $\mathbb{R}_+^2 \times [\phi_d, \infty)$ and never returns to the *advantageous* region $\mathbf{C}_0$ after the first jump; see Section 4.1. In fact,

$$v \circ S(\phi_0, \phi_1) = \inf_{\tau \in \mathcal{S}} \mathbb{E}_0^{S(\phi_0,\phi_1)}\left[\int_0^\tau e^{-\lambda u} g(\widetilde{\boldsymbol{\Phi}}_u)\, du\right] = 0,$$

$$(\phi_0, \phi_1) \in \mathbb{R}_+ \times [\phi_d, \infty),$$

and therefore,

$$v(\phi_0, \phi_1) = J_0 v(\phi_0, \phi_1)$$

(11.25)
$$= \inf_{t \in [0,\infty]} \int_0^t e^{-(\lambda+\mu)u} \overbrace{[g + \mu \cdot v \circ S]}^{\equiv g(\cdot,\cdot)}(x(u,\phi_0), y(u,\phi_1))\, du$$
$$= J_0 v_0(\phi_0, \phi_1) = v_1(\phi_0, \phi_1), \qquad (\phi_0, \phi_1) \in \mathbb{R}_+ \times [\phi_d, \infty).$$

The stopping region $\boldsymbol{\Gamma} = \{(x,y) \in \mathbb{R}_+^2 : v(x,y) = 0\}$ and its boundary $\partial \boldsymbol{\Gamma}$ are determined by the value function $v(\cdot,\cdot)$. Then (11.25) implies that the restrictions of the boundaries $\partial \boldsymbol{\Gamma}$ and $\partial \boldsymbol{\Gamma}_1$ to the region $\mathbb{R}_+ \times [\phi_d, \infty)$ also coincide. Therefore, the second inequality in (11.20) implies

$$\phi_d < \frac{\lambda}{c}\sqrt{2} \leq \gamma_1(0) = \gamma(0) \quad \text{and} \quad S^{-1}[\gamma](0) \equiv \frac{\mu}{\mu+1}\gamma(0) < \phi_d$$

follows from (11.24). Since the boundary function $S^{-1}[\gamma](\cdot)$ of the region $S^{-1}(\boldsymbol{\Gamma})$ is decreasing [see (8.7)], the second inequality gives

$$\mathbb{R}_+ \times [\phi_d, \infty) \subseteq S^{-1}(\boldsymbol{\Gamma}) = \{(x,y) \in \mathbb{R}_+^2 : S^{-1}[\gamma](x) \leq y\}.$$

But starting at any $(\phi_0, \phi_1) \in \mathbb{R} \times [0, \phi_d]$, the parametric curves $t \mapsto (x(t,\phi_0), y(t,\phi_1))$, $t \in \mathbb{R}_+$, are increasing. Since the boundary functions $S^{-n}[\gamma](\cdot)$ of the regions $S^{-n}(\mathbf{\Gamma}) = \{(x,y) \in \mathbb{R}_+^2 : S^{-n}[\gamma](x) \leq y\}$, $n \in \mathbb{N}$, are also decreasing, every region $S^{-n}(\mathbf{\Gamma})$, $n \in \mathbb{N}$, is "closed" in the sense of (8.4). Therefore, Method A after Corollary 8.4 can be used in order to calculate the value function $v(\cdot,\cdot)$ on $\mathbb{R}_+^2$.

COROLLARY 11.22. *Suppose that* (11.20) *holds. Let* $B_n \triangleq [\mu/(\mu+1)]^n \gamma_1(0)$ *for every* $n \in \mathbb{N}$. *Then the sequence* $\mathbb{R}_+ \times [B_n, \infty)$, $n \in \mathbb{N}$, *increases to* $\mathbb{R}_+ \times (0, \infty)$, *and we have* $v(\cdot,\cdot) = v_n(\cdot,\cdot)$ *on* $\mathbb{R}_+ \times [B_n, \infty)$ *for every* $n \in \mathbb{N}$.

Since $(\phi_0^*, \phi_1^*) \in \mathbb{R}_+ \times [\phi_d, \infty) \subseteq \mathbb{R}_+ \times [B_1, \infty)$, the exit boundaries $\partial \mathbf{\Gamma}_1^x$ and $\partial \mathbf{\Gamma}^x$ of the stopping regions $\mathbf{\Gamma}_1$ and $\mathbf{\Gamma}$ are the same, and

$$\partial \mathbf{\Gamma}^x = \partial \mathbf{\Gamma}_1^x = \{(x, \gamma_1(x)) : x \in [0, \xi_1^e)\} \equiv \{(x, \gamma_1(x)) : x \in [0, \phi_0^*]\}.$$

From the entrance boundary $\partial \mathbf{\Gamma}_1^e = \{(x, a_0(x)) : x \in (\phi_0^*, (\lambda/c)\sqrt{2})\}$ of the stopping region $\mathbf{\Gamma}_1$, we can obtain its exit boundary $\partial \mathbf{\Gamma}_1^x$ and the value function $v_1(\cdot,\cdot)$ on the continuation region $\mathbf{C}_1$ by using Method D in Section 11.1; see Figures 5(a), (b) and 6(b).

Note also that the value function $v(\cdot,\cdot) \equiv v_1(\cdot,\cdot)$ is continuously differentiable on $\mathbb{R}_+ \times [B_1, \infty) \setminus \partial \mathbf{\Gamma}_1^x$ and is not differentiable on $\partial \mathbf{\Gamma}_1^x$ by Corollary 11.5 and Lemma 11.16. Let $x_1 \equiv x_1(\gamma_1) = \min\{x \in \mathbb{R}_+ : S^{-1}[\gamma_1](x) = \gamma_1(x)\}$ be the (smallest) intersection point of the functions $S^{-1}[\gamma_1](\cdot)$ and $\gamma_1(\cdot)$ as in (8.9). Then Corollary 8.4 implies

$$\{(x, \gamma(x)) : x \in \mathbb{R}_+\} \cap S^{-1}(\mathbf{\Gamma}_1) = \{(x, \gamma_1(x)) : x \in [0, x_1]\},$$

and the restriction of the boundary function $\gamma(\cdot) \equiv \gamma_1(\cdot)$ to the interval $[0, x_1)$ is continuously differentiable by Lemma 11.15.

Using Corollary 11.22, we can also show that the restrictions of the value function $v(\cdot,\cdot)$ and the boundary $\partial \mathbf{\Gamma}$ of the stopping region $\mathbf{\Gamma}$ on the complement of the region $\mathbb{R}_+ \times [B_1, \infty)$ are continuously differentiable.

Since the sequence $\{v_n(\cdot,\cdot)\}_{n \in \mathbb{N}}$ of the value functions increases to the function $v(\cdot,\cdot)$, all of them coincide with $v(\cdot,\cdot) \equiv v_1(\cdot,\cdot)$ on the region $\mathbb{R}_+ \times [B_1, \infty)$. On the region $\mathbb{R}_+ \times [0, B_1)$, they differ, but are continuously differentiable.

In fact, since every parametric curve $t \mapsto (x(t, \phi_0), y(t, \phi_1))$, $t \in \mathbb{R}_+$, starting at any point $(\phi_0, \phi_1) \in \mathbb{R}_+ \times [0, \phi_d] \supset \mathbb{R}_+ \times [0, B_1]$ is increasing, the hypothesis $A_1(n)$ on page 46 holds on the region $\mathbb{R}_+ \times [0, \phi_d]$ for every $n \in \mathbb{N}$.

On the other hand, the third inequality in (11.20) guarantees that hypothesis $A_2(n)$ on page 46 also holds on $\mathbb{R}_+ \times [0, \phi_d]$ for every $n \in \mathbb{N}$. Indeed, every entrance boundary $\partial \mathbf{\Gamma}_{n+1}^e$ coincides with some part of the boundary $\partial A_n = \{(x, a_n(x)) : x \in \mathbb{R}_+\}$ of the region $A_n = \{(x,y) \in \mathbb{R}_+^2 : [g + \mu \cdot v_n \circ S](x,y) < 0\}$; see Lemma 9.5. Since the sequence $\{a_n(\cdot)\}_{n \in \mathbb{N}_0}$ of the boundary functions is increasing, the third inequality in (11.20) implies

$$y < S[a_0](x) \leq S[a_n](x), \qquad n \in \mathbb{N}_0, (x,y) \in \{(\phi_0^*, \phi_1^*), (0, \bar{\phi}_1)\}.$$

Thus, by an induction on $n \in \mathbb{N}_0$, we can easily show that the transformation $S(\partial \mathbf{\Gamma}^e_{n+1})$ of the entrance boundary $\partial \mathbf{\Gamma}^e_{n+1}$ of every stopping region $\mathbf{\Gamma}_{n+1}$ is away from the exit boundary $\partial \mathbf{\Gamma}^x_{n+1} \equiv \partial \mathbf{\Gamma}^x_1$. Therefore, the function $(x,y) \mapsto [g + \mu \cdot v_n \circ S](x,y)$ is differentiable on the entrance boundary $\partial \mathbf{\Gamma}^e_{n+1}$. The same induction, as in Section 11.2, will also prove the continuous differentiability of the value functions $v_n(\cdot,\cdot)$, $n \in \mathbb{N}$ and $v(\cdot,\cdot)$ on the region $\mathbb{R}_+ \times [0, \phi_d] \supset \mathbb{R}_+ \times [0, B_1]$, as well as the continuous differentiability of the restrictions of the boundaries $\partial \mathbf{\Gamma}_n$, $n \in \mathbb{N}$, and $\partial \mathbf{\Gamma}$ to the set $\mathbb{R}_+ \times [0, \phi_d]$.

COROLLARY 11.23. *Suppose that* (11.20) *holds. Then the boundary function* $\gamma(\cdot)$ *of the stopping region* $\mathbf{\Gamma} = \{(x,y) \in \mathbb{R}^2_+ : \gamma(x) \leq y\}$ *is continuously differentiable on its support* $[0, \xi]$. *The exit boundary* $\mathbf{\Gamma}^x$ *is not empty. The value function* $v(\cdot,\cdot)$ *is continuously differentiable on* $\mathbb{R}^2_+ \setminus \partial \mathbf{\Gamma}^x$, *but not differentiable on* $\partial \mathbf{\Gamma}^x$.

The interesting feature of the solutions of the problems covered under condition (11.20) is that the smooth-fit principle is satisfied on one connected proper subset of the (connected and continuously differentiable) boundary of the optimal stopping region, and fails on the complement of this subset. Moreover, the value function is continuously differentiable away from the boundary.

The conditions in (11.20) are satisfied, for example, if $\lambda = 0.15$, $\mu = 1.5$ and $c = 0.7$ and $m = 0.9$. In general, the functions $S^{-1}[a_0](\cdot)$ and $a_0(\cdot)$ always intersect on the line $y = x$. Since $\gamma_1(\cdot) \geq a_0(\cdot)$ and $\gamma_1(\cdot)$ is decreasing, we have $x_1 \leq (\lambda/c) \cdot (\sqrt{2}/2)$, with equality if and only if

$$S^{-1}(\phi_0^*, \phi_1^*) \in \{(x,y) \in \mathbb{R}^2_+ : x < y\} \iff 1 < \mu(\lambda + c).$$

This condition is satisfied for the numbers above. As a result, we have $x_1 = (\lambda/c) \cdot (\sqrt{2}/2)$ and $\gamma(x) = a_0(x) = x - (\lambda/c)\sqrt{2}$ for every $x \in [\phi_0^*, x_1]$. The boundary function $\gamma(\cdot)$ is strictly above the function $a_0(\cdot)$ everywhere else.

## APPENDIX: PROOFS OF SELECTED RESULTS

*The* $\mathbb{P}_0$*-infinitesimal generator* $\widetilde{\mathcal{A}}$ *of the process* $\widetilde{\mathbf{\Phi}}$ *in* (4.5)  Let us denote by $\widetilde{\mathcal{A}}$ the infinitesimal generator under $\mathbb{P}_0$ of the process $\widetilde{\mathbf{\Phi}} = [\widetilde{\Phi}^{(0)} \ \widetilde{\Phi}^{(1)}]^{\mathrm{T}}$ in (4.5). For every function $f \in \mathbf{C}^{1,1}(\mathbb{R}_+ \times \mathbb{R}_+)$, we have

$$
\begin{aligned}
f(\widetilde{\mathbf{\Phi}}_t) = {}& f(\widetilde{\mathbf{\Phi}}_0) + \sum_{0 < s \leq t} [f(\widetilde{\mathbf{\Phi}}_s) - f(\widetilde{\mathbf{\Phi}}_{s-})] \\
& + \int_0^t \bigg\{ D_{\phi_0} f(\widetilde{\mathbf{\Phi}}_s) \bigg[ (\lambda+1)\widetilde{\Phi}^{(0)}_s + \frac{\lambda(1-m)}{\sqrt{2}} \bigg] \\
& \qquad + D_{\phi_1} f(\widetilde{\mathbf{\Phi}}_s) \bigg[ (\lambda-1)\widetilde{\Phi}^{(1)}_s + \frac{\lambda(1+m)}{\sqrt{2}} \bigg] \bigg\} ds,
\end{aligned}
\tag{A.1}
$$

ADAPTIVE POISSON DISORDER PROBLEM 59

and $\sum_{0<s\leq t}[f(\widetilde{\boldsymbol{\Phi}}_s) - f(\widetilde{\boldsymbol{\Phi}}_{s-})]$ equals

$$\int_0^t \left[ f\left(\left(1-\frac{1}{\mu}\right)\cdot\widetilde{\Phi}^{(0)}_{s-}, \left(1+\frac{1}{\mu}\right)\cdot\widetilde{\Phi}^{(1)}_{s-}\right) - f(\widetilde{\Phi}^{(0)}_{s-},\widetilde{\Phi}^{(1)}_{s-}) \right] dN_s.$$

Note that $\{N_t - \mu t; t \geq 0\}$ is a $(\mathbb{P}_0, \mathbb{F})$-martingale. Then for every $\mathbb{F}$-stopping time $\tau$ such that

$$\mathbb{E}_0^{\phi_0,\phi_1}|f(\widetilde{\boldsymbol{\Phi}}_\tau)| < \infty \quad \text{and}$$

(A.2) $\quad \mathbb{E}_0^{\phi_0,\phi_1}\left[\int_0^\tau \left| f\left(\left(1-\frac{1}{\mu}\right)\cdot\widetilde{\Phi}^{(0)}_{s-}, \left(1+\frac{1}{\mu}\right)\cdot\widetilde{\Phi}^{(1)}_{s-}\right) \right.\right.$

$$\left.\left. - f(\widetilde{\Phi}^{(0)}_{s-},\widetilde{\Phi}^{(1)}_{s-}) \right| ds \right] < \infty,$$

we have

(A.3) $\quad \mathbb{E}_0 f(\widetilde{\boldsymbol{\Phi}}_\tau) = f(\widetilde{\boldsymbol{\Phi}}_0) + \mathbb{E}_0 \int_0^\tau \widetilde{\mathcal{A}} f(\widetilde{\boldsymbol{\Phi}}_s)\, ds, \qquad t \geq 0,$

and $\widetilde{\mathcal{A}} f(\phi_0, \phi_1)$ equals

$$D_{\phi_0} f(\phi_0, \phi_1) \left[ (\lambda + 1)\phi_0 + \frac{\lambda(1-m)}{\sqrt{2}} \right]$$

$$+ D_{\phi_1} f(\phi_0, \phi_1) \left[ (\lambda - 1)\phi_1 + \frac{\lambda(1+m)}{\sqrt{2}} \right]$$

(A.4)

$$+ \mu \left[ f\left(\left(1 - \frac{1}{\mu}\right)\phi_0, \left(1 + \frac{1}{\mu}\right)\phi_1 \right) - f(\phi_0, \phi_1) \right],$$

$$(\phi_0, \phi_1) \in \mathbb{R}_+ \times \mathbb{R}_+.$$

PROOF OF LEMMA 5.3. Let $w: \mathbb{R}_+^2 \mapsto \mathbb{R}$ be a bounded Borel function. Since $g(\cdot, \cdot) \geq g(0,0) = -\lambda\sqrt{2}/c$ in (4.12) is bounded from below, the function $J_0 w$ is well defined. By (5.7),

$$Jw(t, \phi_0, \phi_1) \geq -\left(\frac{\lambda}{c}\sqrt{2} + \mu\|w\|\right) \int_0^\infty e^{-(\lambda+\mu)u}\, du$$

$$= -\left(\frac{\lambda}{c}\sqrt{2} + \mu\|w\|\right) \frac{1}{\lambda + \mu}$$

for every $t \in [0, \infty]$. Since we also have $J_0 w(\phi_0, \phi_1) \leq Jw(0, \phi_0, \phi_1) = 0$, we obtain (5.9).

Suppose now that $w$ is also concave. For every $u \in \mathbb{R}$, the functions $\phi_0 \mapsto x(u, \phi_0)$ and $\phi_1 \mapsto y(u, \phi_1)$ in (4.8) are linear. The mappings $(\phi_0, \phi_1) \mapsto S(\phi_0, \phi_1)$ in (5.8) and $(\phi_0, \phi_1) \mapsto g(\phi_0, \phi_1)$ in (4.12) are also linear. Therefore, the integrand in (5.7), namely $(\phi_0, \phi_1) \mapsto e^{-(\lambda+\mu)u}(g + \mu \cdot w \circ S)(x(u, \phi_0), y(u, \phi_1))$, is concave for every $u \in [0, \infty)$. Thus, the mappings $(\phi_0, \phi_1) \mapsto Jw(t, (\phi_0, \phi_1))$, $t \in [0, \infty]$ in (5.7) are concave. Then $J_0 w(\phi_0, \phi_1) = \inf_{t \in [0, \infty]} Jw(t, \phi_0, \phi_1)$ is a lower envelope of concave mappings, and therefore, is a concave function of $(\phi_0, \phi_1) \in \mathbb{R}_+^2$. Finally, it is clear from (5.7) that $w_1 \leq w_2$ implies that $J_0 w_1 \leq J_0 w_2$. $\square$

PROOF OF COROLLARY 5.4. The function $v_0 \equiv 0$ has all of the properties. The proof follows from an induction and the properties of concave functions. $\square$

For the proof of Proposition 5.5, we shall need the following result on the characterization of $\mathbb{F}$-stopping times; see [6], Theorem T33, page 308, and [9], Lemma A2.3, page 261.

LEMMA A.1. *For every $\mathbb{F}$-stopping time $\tau$ and every $n \in \mathbb{N}_0$, there is an $\mathcal{F}_{\sigma_n}$-measurable random variable $R_n : \Omega \mapsto [0, \infty]$ such that $\tau \wedge \sigma_{n+1} = (\sigma_n + R_n) \wedge \sigma_{n+1}$ holds $\mathbb{P}_0$-a.s. on $\{\tau \geq \sigma_n\}$.*

PROOF OF PROPOSITION 5.5. First, we shall establish the inequality

$$(A.5) \quad \mathbb{E}_0^{\phi_0, \phi_1} \int_0^{\tau \wedge \sigma_n} e^{-\lambda t} g(\widetilde{\boldsymbol{\Phi}}_t) \, dt \geq v_n(\phi_0, \phi_1), \qquad \tau \in \mathcal{S}, (\phi_0, \phi_1) \in \mathbb{R}_+^2,$$

for every $n \in \mathbb{N}_0$, by proving inductively on $k = 1, \ldots, n+1$ that

$$\mathbb{E}_0^{\phi_0, \phi_1} \int_0^{\tau \wedge \sigma_n} e^{-\lambda t} g(\widetilde{\boldsymbol{\Phi}}_t) \, dt$$

$$(A.6) \quad \geq RHS_{k-1}$$

$$:= \mathbb{E}_0^{\phi_0, \phi_1} \left[ \left( \int_0^{\tau \wedge s} e^{-\lambda t} g(\widetilde{\boldsymbol{\Phi}}_t) \, dt + \mathbf{1}_{\{\tau \geq s\}} e^{-\lambda s} v_{k-1}(\widetilde{\boldsymbol{\Phi}}_s) \right) \bigg|_{s = \sigma_{n-k+1}} \right].$$

Observe that (A.5) follows from (A.6) when we set $k = n + 1$.

If $k = 1$, then the inequality (A.6) is satisfied as an equality, since $v_0 \equiv 0$. Suppose that (A.6) holds for some $1 \leq k < n + 1$. We shall prove that it must also hold when $k$ is replaced with $k + 1$. Let us denote the right-hand side of (A.6) by $RHS_{k-1}$, and rewrite it as

$$RHS_{k-1} = RHS_{k-1}^{(1)} + RHS_{k-1}^{(2)} \triangleq \mathbb{E}_0^{\phi_0, \phi_1} \left[ \int_0^{\tau \wedge \sigma_{n-k}} e^{-\lambda t} g(\widetilde{\boldsymbol{\Phi}}_t) \, dt \right]$$

$$(A.7) \qquad + \mathbb{E}_0^{\phi_0, \phi_1} \left[ \mathbf{1}_{\{\tau \geq \sigma_{n-k}\}} \left( \int_{\sigma_{n-k}}^{\tau \wedge s} e^{-\lambda t} g(\widetilde{\boldsymbol{\Phi}}_t) \, dt \right. \right.$$

$$\left. \left. + \mathbf{1}_{\{\tau \geq s\}} e^{-\lambda s} v_{k-1}(\widetilde{\boldsymbol{\Phi}}_s) \right) \bigg|_{s = \sigma_{n-k+1}} \right].$$

By Lemma A.1, there is an $\mathcal{F}_{\sigma_{n-k}}$-measurable random variable $R_{n-k}$ such that $\tau \wedge \sigma_{n-k+1} = (\sigma_{n-k} + R_{n-k}) \wedge \sigma_{n-k+1}$ holds $\mathbb{P}_0$-almost surely on $\{\tau \geq \sigma_{n-k}\}$. Therefore, the second expectation, denoted by $RHS^{(2)}_{k-1}$, in (A.7) becomes

$$\mathbb{E}_0^{\phi_0,\phi_1} \mathbf{1}_{\{\tau \geq t\}} \left[ \int_t^{(t+R_{n-k}) \wedge s} e^{-\lambda u} g(\widetilde{\boldsymbol{\Phi}}_u)\, du + \mathbf{1}_{\{t+R_{n-k} \geq s\}} e^{-\lambda s} v_{k-1}(\widetilde{\boldsymbol{\Phi}}_s) \right] \Bigg|_{\substack{t=\sigma_{n-k} \\ s=\sigma_{n-k+1}}}$$

$$= \mathbb{E}_0^{\phi_0,\phi_1} \{\mathbf{1}_{\{\tau \geq \sigma_{n-k}\}} e^{-\lambda \sigma_{n-k}} f_{n-k}(R_{n-k}, \widetilde{\boldsymbol{\Phi}}_{\sigma_{n-k}})\}$$

by the strong Markov property of $N$, where

$$f_{k-1}(r,\phi_0,\phi_1) \triangleq \mathbb{E}_0^{\phi_0,\phi_1}\left[ \int_0^{r \wedge \sigma_1} e^{-\lambda t} g(\widetilde{\boldsymbol{\Phi}}_t)\, dt + \mathbf{1}_{\{r \geq \sigma_1\}} e^{-\lambda \sigma_1} v_{k-1}(\widetilde{\boldsymbol{\Phi}}_{\sigma_1}) \right]$$

$$= J v_{k-1}(r,(\phi_0,\phi_1)) \geq J_0 v_{k-1}(\phi_0,\phi_1) = v_k(\phi_0,\phi_1).$$

The (in)equalities follow from (5.3), (5.4) and (5.6), respectively. Thus

$$RHS^{(2)}_{k-1} \geq \mathbb{E}_0^{\phi_0,\phi_1}[\mathbf{1}_{\{\tau \geq \sigma_{n-k}\}} e^{-\lambda \sigma_{n-k}} v_k(\widetilde{\boldsymbol{\Phi}}_{\sigma_{n-k}})].$$

From (A.6) and (A.7), we obtain

$$\mathbb{E}_0^{\phi_0,\phi_1} \int_0^{\tau \wedge \sigma_n} e^{-\lambda t} g(\widetilde{\boldsymbol{\Phi}}_t)\, dt$$

$$\geq RHS_{k-1} = \mathbb{E}_0^{\phi_0,\phi_1}\left[ \int_0^{\tau \wedge \sigma_{n-k}} e^{-\lambda t} g(\widetilde{\boldsymbol{\Phi}}_t)\, dt \right] + RHS^{(2)}_{k-1}$$

$$\geq \mathbb{E}_0^{\phi_0,\phi_1}\left[ \int_0^{\tau \wedge \sigma_{n-k}} e^{-\lambda t} g(\widetilde{\boldsymbol{\Phi}}_t)\, dt + \mathbf{1}_{\{\tau \geq \sigma_{n-k}\}} e^{-\lambda \sigma_{n-k}} v_k(\widetilde{\boldsymbol{\Phi}}_{\sigma_{n-k}}) \right] = RHS_k.$$

This completes the proof of (A.6) by induction on $k$, and (A.5) follows by setting $k = n+1$ in (A.6). When we take the infimum of both sides in (A.5), we obtain $V_n \geq v_n$, $n \in \mathbb{N}$.

The reverse inequality $V_n \leq v_n$, $n \in \mathbb{N}$, follows immediately from (5.11), since every $\mathbb{F}$-stopping time $S_n^\varepsilon$ is less than or equal to $\sigma_n$, $\mathbb{P}_0$-a.s. by construction. Therefore, we only need to establish (5.11). We shall prove it by induction on $n \in \mathbb{N}$. For $n = 1$, the left-hand side of (5.11) becomes

$$\mathbb{E}_0^{\phi_0,\phi_1} \int_0^{S_1^\varepsilon} e^{-\lambda t} g(\widetilde{\boldsymbol{\Phi}}_t)\, dt = \mathbb{E}_0^{\phi_0,\phi_1} \int_0^{r_0^\varepsilon(\phi_0,\phi_1) \wedge \sigma_1} e^{-\lambda t} g(\widetilde{\boldsymbol{\Phi}}_t)\, dt$$

$$= J v_0(r_0^\varepsilon(\phi_0,\phi_1), \phi_0, \phi_1).$$

Since $J v_0(r_0^\varepsilon(\phi_0,\phi_1),\phi_0,\phi_1) \leq J_0 v_0(\phi_0,\phi_1) + \varepsilon$ by Remark 5.2, (5.11) holds for $n = 1$.

Suppose that (5.11) holds for every $\varepsilon > 0$ for some $n \in \mathbb{N}$. We will prove that it also holds when $n$ is replaced with $n+1$. Since $S_{n+1}^\varepsilon \wedge \sigma_1 = r_n^{\varepsilon/2}(\widetilde{\boldsymbol{\Phi}}_0) \wedge \sigma_1$, $\mathbb{P}_0$-a.s., the expectation $\mathbb{E}_0^{\phi_0,\phi_1}[\int_0^{S_{n+1}^\varepsilon} e^{-\lambda t} g(\widetilde{\boldsymbol{\Phi}}_t)\,dt]$ becomes

$$\mathbb{E}_0^{\phi_0,\phi_1}\left[\int_0^{S_{n+1}^\varepsilon \wedge \sigma_1} e^{-\lambda t} g(\widetilde{\boldsymbol{\Phi}}_t)\,dt + \mathbf{1}_{\{S_{n+1}^\varepsilon \geq \sigma_1\}} \int_{\sigma_1}^{S_{n+1}^\varepsilon} e^{-\lambda t} g(\widetilde{\boldsymbol{\Phi}}_t)\,dt\right]$$

$$= \mathbb{E}_0^{\phi_0,\phi_1}\left[\int_0^{r_n^{\varepsilon/2}(\phi_0,\phi_1) \wedge \sigma_1} e^{-\lambda t} g(\widetilde{\boldsymbol{\Phi}}_t)\,dt\right]$$

$$+ \mathbb{E}_0^{\phi_0,\phi_1}[\mathbf{1}_{\{r_n^{\varepsilon/2}(\phi_0,\phi_1) \geq \sigma_1\}} e^{-\lambda \sigma_1} f_n(\widetilde{\boldsymbol{\Phi}}_{\sigma_1})]$$

by the strong Markov property of $N$, where

$$f_n(\phi_0,\phi_1) \triangleq \mathbb{E}_0^{\phi_0,\phi_1}\left[\int_0^{S_n^{\varepsilon/2}} e^{-\lambda t} g(\widetilde{\boldsymbol{\Phi}}_t)\,dt\right] \leq v_n(\phi_0,\phi_1) + \varepsilon/2$$

by the induction hypothesis. Therefore,

$$\mathbb{E}_0^{\phi_0,\phi_1}\left[\int_0^{S_{n+1}^\varepsilon} e^{-\lambda t} g(\widetilde{\boldsymbol{\Phi}}_t)\,dt\right]$$

$$\leq \mathbb{E}_0^{\phi_0,\phi_1}\left[\int_0^{r_n^{\varepsilon/2}(\phi_0,\phi_1) \wedge \sigma_1} e^{-\lambda t} g(\widetilde{\boldsymbol{\Phi}}_t)\,dt\right.$$

(A.8)
$$\left.+ \mathbf{1}_{\{r_n^{\varepsilon/2}(\phi_0,\phi_1) \geq \sigma_1\}} e^{-\lambda \sigma_1} v_n(\widetilde{\boldsymbol{\Phi}}_{\sigma_1})\right] + \varepsilon/2$$

$$= Jv_n(r_n^{\varepsilon/2}(\phi_0,\phi_1),(\phi_0,\phi_1)) + \varepsilon/2.$$

Since $Jv_n(r_n^{\varepsilon/2}(\phi_0,\phi_1),(\phi_0,\phi_1)) \leq v_{n+1}(\phi_0,\phi_1) + \varepsilon/2$ by Remark 5.2, this inequality and (A.8) prove (5.11) when $n$ is replaced with $n+1$. □

PROOF OF PROPOSITION 5.6. Corollary 5.4 and Propositions 5.5 and 5.1 imply that $v(\phi_0,\phi_1) = \lim_{n\to\infty} v_n(\phi_0,\phi_1) = \lim_{n\to\infty} V_n(\phi_0,\phi_1) = V(\phi_0,\phi_1)$ for every $(\phi_0,\phi_1) \in \mathbb{R}_+^2$. Next, let us show that $V = J_0 V$. Since $(v_n)_{n\geq 1}$ is a decreasing sequence,

(A.9) $V(\phi_0,\phi_1) = \lim_{n\to\infty} v_n(\phi_0,\phi_1) = \inf_{n\geq 1} v_n(\phi_0,\phi_1) = \inf_{n\geq 1} J_0 v_{n-1}(\phi_0,\phi_1)$

for every $(\phi_0,\phi_1) \in \mathbb{R}_+^2$. Since $(Jv_n)_{n\geq 1}$ is a decreasing sequence, and $\{v_n\}_{n\in\mathbb{N}}$ are uniformly bounded, we have $V(\phi_0,\phi_1) = \inf_{n\geq 1} J_0 v_{n-1}(\phi_0,\phi_1) = J_0 v(\phi_0,\phi_1) = J_0 V(\phi_0,\phi_1)$ by the dominated convergence theorem and (A.9). Finally, since $U \leq 0$, we have $U \leq v_n$ for every $n$ by induction, and $U \leq \lim_{n\to\infty} v_n = V$.
□

PROOF OF LEMMA 5.7. Let us fix a constant $u \geq t$ and $(\phi_0, \phi_1) \in \mathbb{R}_+^2$. Then

$$
\begin{aligned}
Jw&(u, \phi_0, \phi_1) \\
&= \mathbb{E}_0^{\phi_0, \phi_1}\left[\int_0^{t \wedge \sigma_1} e^{-\lambda s} g(\widetilde{\boldsymbol{\Phi}}_s)\, ds + \mathbf{1}_{\{u \geq \sigma_1\}} e^{-\lambda \sigma_1} w(\widetilde{\boldsymbol{\Phi}}_{\sigma_1})\right] \\
&\quad + \mathbb{E}_0^{\phi_0, \phi_1}\left[\mathbf{1}_{\{\sigma_1 > t\}} \int_t^{u \wedge \sigma_1} e^{-\lambda s} g(\widetilde{\boldsymbol{\Phi}}_s)\, ds\right].
\end{aligned}
\tag{A.10}
$$

On the event $\{\sigma_1 > t\}$, we have $u \wedge \sigma_1 = [t + (u-t)] \wedge [t + (\sigma_1 \circ \theta_t)] = t + [(u-t) \wedge (\sigma_1 \circ \theta_t)]$. Therefore, the strong Markov property of $N$ gives

$$
\begin{aligned}
\mathbb{E}_0^{\phi_0, \phi_1}& \mathbf{1}_{\{\sigma_1 > t\}} \int_t^{u \wedge \sigma_1} e^{-\lambda s} g(\widetilde{\boldsymbol{\Phi}}_s)\, ds \\
&= \mathbb{E}_0^{\phi_0, \phi_1} \mathbf{1}_{\{\sigma_1 > t\}} e^{-\lambda t} \mathbb{E}_0^{\widetilde{\boldsymbol{\Phi}}_t}\left[\int_0^{(u-t) \wedge \sigma_1} e^{-\lambda s} g(\widetilde{\boldsymbol{\Phi}}_s)\, ds\right] \\
&= \mathbb{E}_0^{\phi_0, \phi_1}[\mathbf{1}_{\{\sigma_1 > t\}} e^{-\lambda t}(Jw(u-t, \widetilde{\boldsymbol{\Phi}}_t) - \mathbb{E}_0^{\widetilde{\boldsymbol{\Phi}}_t}[\mathbf{1}_{\{u-t \geq \sigma_1\}} e^{-\lambda \sigma_1} w(\widetilde{\boldsymbol{\Phi}}_{\sigma_1})])] \\
&= e^{-(\lambda+\mu)t} Jw(u-t, (x(t, \phi_0), y(t, \phi_0))) \\
&\quad - \mathbb{E}_0^{\phi_0, \phi_1}[\mathbf{1}_{\{\sigma_1 > t\}} \mathbf{1}_{\{u \geq \sigma_1\}} e^{-\lambda \sigma_1} w(\widetilde{\boldsymbol{\Phi}}_{\sigma_1})].
\end{aligned}
$$

The third equality follows from the definition of $Jw$ in (5.3) and the last from (4.10) and the strong Markov property. Substituting the last equation into (A.10) and simplifying the rest give

$$Jw(u, \phi_0, \phi_1) = Jw(t, (\phi_0, \phi_1)) + e^{-(\lambda+\mu)t} Jw(u-t, (x(t, \phi_0), y(t, \phi_0))).$$

Finally, taking the infimum of both sides over $u \in [t, +\infty]$ gives (5.12). $\square$

PROOF OF PROPOSITION 5.11. First, let us show (5.22) for $n = 1$. Fix $\varepsilon \geq 0$ and $(\phi_0, \phi_1) \in \mathbb{R}_+^2$. By Lemma A.1, there exists a constant $u \in [0, \infty]$ such that $U_\varepsilon \wedge \sigma_1 = u \wedge \sigma_1$. Then

$$
\begin{aligned}
\mathbb{E}_0^{\phi_0, \phi_1} M_{U_\varepsilon \wedge \sigma_1} &= \mathbb{E}_0^{\phi_0, \phi_1}\left[e^{-\lambda(u \wedge \sigma_1)} V(\widetilde{\boldsymbol{\Phi}}_{u \wedge \sigma_1}) + \int_0^{u \wedge \sigma_1} e^{-\lambda s} g(\widetilde{\boldsymbol{\Phi}}_s)\, ds\right] \\
&= \mathbb{E}_0^{\phi_0, \phi_1}\left[\int_0^{u \wedge \sigma_1} e^{-\lambda s} g(\widetilde{\boldsymbol{\Phi}}_s)\, ds + \mathbf{1}_{\{u \geq \sigma_1\}} e^{-\lambda \sigma_1} V(\widetilde{\boldsymbol{\Phi}}_{\sigma_1})\right] \\
&\quad + \mathbb{E}_0^{\phi_0, \phi_1}[\mathbf{1}_{\{u < \sigma_1\}} e^{-\lambda u} V(\widetilde{\boldsymbol{\Phi}}_u)] \\
&= JV(u, (\phi_0, \phi_1)) + e^{-(\lambda+\mu)u} V(x(u, \phi_0), y(u, \phi_1)) \\
&= J_u V(\phi_0, \phi_1),
\end{aligned}
\tag{A.11}
$$

where the third equality follows from (5.3) and (4.10), and the fourth from (5.16).

Fix any $t \in [0, u)$. By (5.16) and (4.10) once again, we have

$$JV(t, \phi_0, \phi_1) = J_t V(\phi_0, \phi_1) - e^{-(\lambda+\mu)t} V(x(t, \phi_0), y(t, \phi_1))$$
$$\geq J_0 V(\phi_0, \phi_1) - \mathbb{E}_0^{\phi_0, \phi_1}[\mathbf{1}_{\{\sigma_1 > t\}} e^{-\lambda t} V(\widetilde{\boldsymbol{\Phi}}_t)].$$

On the event $\{\sigma_1 > t\}$, we have $U_\varepsilon > t$ (otherwise, $U_\varepsilon \leq t < \sigma_1$ would imply $U_\varepsilon = u \leq t$, contradicts with our initial choice of $t < u$). Thus, $V(\widetilde{\boldsymbol{\Phi}}_t) < -\varepsilon$ on $\{\sigma_1 > t\}$. Hence,

$$JV(t, \phi_0, \phi_1) > J_0 V(\phi_0, \phi_1) + \varepsilon e^{-(\lambda+\mu)u} \geq J_0 V(\phi_0, \phi_1)$$

for every $t \in [0, u)$. Therefore, $J_0 V(\phi_0, \phi_1) = J_u V(\phi_0, \phi_1)$, and (A.11) implies

$$\mathbb{E}_0^{\phi_0, \phi_1}[M_{U_\varepsilon \wedge \sigma_1}] = J_u V(\phi_0, \phi_1) = J_0 V(\phi_0, \phi_1) = V(\phi_0, \phi_1) = \mathbb{E}_0^{\phi_0, \phi_1}[M_0].$$

This completes the proof of (5.22) for $n = 1$. Now suppose that (5.22) holds for some $n \in \mathbb{N}$, and let us show the same equality for $n+1$. Note that

$$\mathbb{E}_0^{\phi_0, \phi_1}[M_{U_\varepsilon \wedge \sigma_{n+1}}]$$
$$= \mathbb{E}_0^{\phi_0, \phi_1}[\mathbf{1}_{\{U_\varepsilon < \sigma_1\}} M_{U_\varepsilon}] + \mathbb{E}_0^{\phi_0, \phi_1}\left[\mathbf{1}_{\{U_\varepsilon \geq \sigma_1\}} \int_0^{\sigma_1} e^{-\lambda s} g(\widetilde{\boldsymbol{\Phi}}_s) \, ds\right]$$
$$+ \mathbb{E}_0^{\phi_0, \phi_1}\left[\mathbf{1}_{\{U_\varepsilon \geq \sigma_1\}} \left\{\int_{\sigma_1}^{U_\varepsilon \wedge \sigma_{n+1}} e^{-\lambda s} g(\widetilde{\boldsymbol{\Phi}}_s) \, ds \right.\right.$$
$$\left.\left. + e^{-\lambda(U_\varepsilon \wedge \sigma_{n+1})} V(\widetilde{\boldsymbol{\Phi}}_{U_\varepsilon \wedge \sigma_{n+1}})\right\}\right].$$

Since $U_\varepsilon \wedge \sigma_{n+1} = \sigma_1 + [(U_\varepsilon \wedge \sigma_n) \circ \theta_{\sigma_1}]$ on the event $\{U_\varepsilon \geq \sigma_1\}$, the strong Markov property of $\widetilde{\boldsymbol{\Phi}}$ at the stopping time $\sigma_1$ will complete the proof. $\square$

PROOF OF PROPOSITION 5.12. Note that the sequence of random variables

$$\int_0^{U_\varepsilon \wedge \sigma_n} e^{-\lambda s} g(\widetilde{\boldsymbol{\Phi}}_s) \, ds + e^{-\lambda(U_\varepsilon \wedge \sigma_n)} V(\widetilde{\boldsymbol{\Phi}}_{U_\varepsilon \wedge \sigma_n}) \geq -2 \int_0^\infty e^{-\lambda s} \frac{\lambda}{c} \sqrt{2} \, ds = -\frac{2\sqrt{2}}{c}$$

is bounded from below; see (4.12). By (5.22) and Fatou's lemma, we have

$$V(\phi_0, \phi_1) \geq \mathbb{E}_0^{\phi_0, \phi_1}\left[\liminf_{n \to \infty} \left(\int_0^{U_\varepsilon \wedge \sigma_n} e^{-\lambda s} g(\widetilde{\boldsymbol{\Phi}}_s) \, ds + e^{-\lambda(U_\varepsilon \wedge \sigma_n)} V(\widetilde{\boldsymbol{\Phi}}_{U_\varepsilon \wedge \sigma_n})\right)\right]$$
$$\geq \mathbb{E}_0^{\phi_0, \phi_1}\left[\int_0^{U_\varepsilon} e^{-\lambda s} g(\widetilde{\boldsymbol{\Phi}}_s) \, ds\right] - \varepsilon$$

for every $(\phi_0, \phi_1) \in \mathbb{R}_+^2$. The second inequality follows from (5.21). $\square$

PROOF OF PROPOSITION 8.3. Let us prove (8.5) for $n = 1$. Take $(\phi_0, \phi_1) \in S^{-1}(\mathbf{\Gamma})$. By (8.4), the curve $u \mapsto (x(u, \phi_0), y(u, \phi_1))$, $u \geq 0$, does not leave $S^{-1}(\mathbf{\Gamma})$. Hence,

$$S(x(u, \phi_0), y(u, \phi_1)) \in \mathbf{\Gamma} \quad \text{and} \quad (V \circ S)(x(u, \phi_0), y(u, \phi_1)) = 0, \qquad u \in \mathbb{R}_+.$$

Then Lemma 5.6, (5.4), (5.6) and Proposition 5.5 imply that

$$\begin{aligned} V(\phi_0, \phi_1) &= J_0 V(\phi_0, \phi_1) \\ &= \inf_{t \in [0, \infty]} \int_0^t e^{-(\lambda+\mu)u} [g + \mu \cdot V \circ S](x(u, \phi_0), y(u, \phi_1))\, du \\ &= \inf_{t \in [0, \infty]} \int_0^t e^{-(\lambda+\mu)u} g(x(u, \phi_0), y(u, \phi_1))\, du \\ &= J_0 V_0(\phi_0, \phi_1) \\ &= V_1(\phi_0, \phi_1). \end{aligned}$$

Since $V$ is the limit of the decreasing sequence $\{V_n\}_{n \in \mathbb{N}}$, the equalities $V = V_1 = V_2 = \cdots$ on $S^{-1}(\mathbf{\Gamma})$ follow.

On $S^{-1}(\mathbf{\Gamma}) \cap \mathbf{C}$, we have $0 > V = V_1 = V_2 = \cdots$. Therefore, $S^{-1}(\mathbf{\Gamma}) \cap \mathbf{C} \subseteq \mathbf{C}_k$ for every $k \geq 1$. Taking intersection of both sides with $S^{-1}(\mathbf{\Gamma})$ gives $S^{-1}(\mathbf{\Gamma}) \cap \mathbf{C} \subseteq S^{-1}(\mathbf{\Gamma}) \cap \mathbf{C}_k$ for every $k \geq 1$. To prove the opposite inclusion, note that $V = V_k < 0$ on $S^{-1}(\mathbf{\Gamma}) \cap \mathbf{C}_k$ for every $k \geq 1$. Therefore, $S^{-1}(\mathbf{\Gamma}) \cap \mathbf{C}_k \subseteq \mathbf{C}$, $k \geq 1$. Intersecting both sides with the set $S^{-1}(\mathbf{\Gamma})$ gives $S^{-1}(\mathbf{\Gamma}) \cap \mathbf{C}_k \subseteq S^{-1}(\mathbf{\Gamma}) \cap \mathbf{C}$, $k \geq 1$.

The proof of $S^{-n}(\mathbf{\Gamma}) \cap \mathbf{\Gamma} = S^{-n}(\mathbf{\Gamma}) \cap \mathbf{\Gamma}_n = S^{-n}(\mathbf{\Gamma}) \cap \mathbf{\Gamma}_{n+1} = \cdots$ reads as in the previous paragraph, after every "$\mathbf{C}$" is replaced by "$\mathbf{\Gamma}$," and every strict inequality by an equality. This completes the proof of (8.5) for $n = 1$.

Suppose that (8.5) holds for some $n \in \mathbb{N}$, and let us prove it for $n + 1$. Take $(\phi_0, \phi_1) \in S^{-(n+1)}(\mathbf{\Gamma})$. Since the curve $u \mapsto (x(u, \phi_0), y(u, \phi_1))$, $u \in \mathbb{R}_+$, does not leave the region $S^{-(n+1)}(\mathbf{\Gamma})$ by (8.4), we have $S(x(u, \phi_0), y(u, \phi_1)) \in S^{-n}(\mathbf{\Gamma})$, $u \in \mathbb{R}_+$, and

$$(V \circ S)(x(u, \phi_0), y(u, \phi_1)) = (V_n \circ S)(x(u, \phi_0), y(u, \phi_1)), \qquad u \in \mathbb{R}_+,$$

by the induction hypothesis. Lemma 5.6, (5.4), (5.6) and Proposition 5.5 imply

$$\begin{aligned} V(\phi_0, \phi_1) &= J_0 V(\phi_0, \phi_1) \\ &= \inf_{t \in [0, \infty]} \int_0^t e^{-(\lambda+\mu)u} [g + \mu \cdot V \circ S](x(u, \phi_0), y(u, \phi_1))\, du \\ &= \inf_{t \in [0, \infty]} \int_0^t e^{-(\lambda+\mu)u} [g + \mu \cdot V_n \circ S](x(u, \phi_0), y(u, \phi_1))\, du \\ &= V_{n+1}(\phi_0, \phi_1). \end{aligned}$$

Since $V$ is the limit of the decreasing sequence $\{V_n\}_{n\in\mathbb{N}}$, we have $V = V_{n+1} = V_{n+2} = \cdots$ on $S^{-(n+1)}(\mathbf{\Gamma})$. From these equalities follows the proof of the equalities of the regions in (8.5) for $n+1$, by arguments similar to those presented for $n=1$ above. $\square$

PROOF OF LEMMA 9.2. The obvious choices are the function $a_n : \mathbb{R}_+ \mapsto \mathbb{R}_+$ and the number $\alpha_n$ in (9.3) and (9.4), respectively. By the discussion above,

$$\begin{aligned}\{(x,y) \in \mathbb{R}_+^2 ; [g + \mu \cdot v_n \circ S](x,y) < 0\} \\ = A_n = \mathbb{R}_+^2 \setminus \mathrm{epi}(a_n) \\ = \{(x,y) \in \mathbb{R}_+^2 ; y < a_n(x)\} \\ = \{(x,y) \in [0,\alpha_n) \times \mathbb{R}_+ ; y < a_n(x)\},\end{aligned}$$

and (9.6) follows. The proof will be complete if we show the equality in (9.5).

Since $[g + \mu \cdot v_n \circ S](x,y)$, $x \in \mathbb{R}_+$, is continuous, we have $[g + \mu \cdot v_n \circ S](x, a_n(x)) \geq 0$ for every $x \in \mathbb{R}_+$, and the equality holds for every $x \in [0, \alpha_n)$ because $a_n(x) > 0$, $x \in [0, \alpha_n)$. Because $a_n(\cdot)$ is also continuous, the equality also holds for $(x,y) = (\alpha_n, a(\alpha_n))$, and

$$[g + \mu \cdot v_n \circ S](x, a_n(x)) = 0, \qquad x \in [0, \alpha_n]. \tag{A.12}$$

The identity (9.5) will follow immediately if we show for the same $A_n$ in (9.1) that

$$\begin{aligned}[g + \mu \cdot v_n \circ S](x,y) > 0, \\ (x,y) \in (\mathbb{R}_+^2 \setminus A_n) \setminus \{(x, a_n(x)) : x \in [0, \alpha_n]\}.\end{aligned} \tag{A.13}$$

The nonpositive function $v_n(\cdot,\cdot)$ is concave and equal to zero outside the bounded region $\mathbf{C}_n$. Therefore, the functions $y \mapsto v_n(x,y)$, $x \in \mathbb{R}_+$, and $x \mapsto v_n(x,y)$, $y \in \mathbb{R}_+$, are nonpositive, concave and equal zero for every large real $y$ and $x$, respectively. This implies that the functions $y \mapsto v_n(x,y)$, $x \in \mathbb{R}_+$, and $x \mapsto v_n(x,y)$, $y \in \mathbb{R}_+$, are nondecreasing. Therefore, the functions $y \mapsto [g + \mu \cdot v_n \circ S](x,y)$, $x \in \mathbb{R}_+$, and $x \mapsto [g + \mu \cdot v_n \circ S](x,y)$, $y \in \mathbb{R}_+$, are strictly increasing since both $S(x,y)$ and $g(x,y)$ are strictly increasing in both $x$ and $y$. Now (A.13) follows from (A.12). $\square$

PROOF OF LEMMA 9.7. Fix any $(\phi_0, \phi_1) \in \partial \mathbf{\Gamma}_{n+1}^e$. Then $v_{n+1}(\phi_0, \phi_1) = 0$, and substituting $(\phi_0^{-t}, \phi_1^{-t}) \triangleq (x(-t, \phi_0), y(-t, \phi_1))$ into (9.15) for $t \in [0, \hat{r}(\phi_0, \phi_1)]$ gives

$$J_t v_n(x(-t, \phi_0), y(-t, \phi_1)) = -e^{-(\lambda+\mu)t} J v_n(-t, \phi_0, \phi_1), \tag{A.14}$$
$$t \in [0, \hat{r}(\phi_0, \phi_1)],$$

thanks to the semigroup property of $x(\cdot,\cdot)$ and $y(\cdot,\cdot)$.

By the definition of the entrance boundary $\partial\mathbf{\Gamma}^e_{n+1}$ in (9.10), the point $(\phi_0,\phi_1)$ is reachable from the inside of the continuation region $\mathbf{C}_{n+1}$. Namely, there exists some $\delta>0$ such that $(x(-t,\phi_0),y(-t,\phi_1))\in\mathbf{C}_{n+1}$ and $r_n(x(-t,\phi_0),y(-t,\phi_1))=t$ for every $t\in(0,\delta]$. Then (9.16) implies

$$0 > v_{n+1}(x(-t,\phi_0),y(-t,\phi_1)) = -e^{-(\lambda+\mu)t}Jv_n(-t,\phi_0,\phi_1)$$

for every $t\in(0,\delta]$. Since $\hat{r}_n(\phi_0,\phi_1)$ is the first time when the last function on the right may change its sign, we obtain

$$-Jv_n(-t,\phi_0,\phi_1)<0, \qquad t\in(0,\hat{r}_n(\phi_0,\phi_1)\wedge\hat{r}(\phi_0,\phi_1)).$$

Using (9.16) once again, we conclude

$$v_{n+1}(x(-t,\phi_0),y(-t,\phi_1))$$
$$\leq J_t v_n(x(-t,\phi_0),y(-t,\phi_1))$$
$$= -e^{-(\lambda+\mu)t}Jv_n(-t,\phi_0,\phi_1)<0, \qquad t\in(0,\hat{r}_n(\phi_0,\phi_1)\wedge\hat{r}(\phi_0,\phi_1)).$$

Thus $\{(x(-t,\phi_0),y(-t,\phi_1)); t\in(0,\hat{r}_n(\phi_0,\phi_1)\wedge\hat{r}(\phi_0,\phi_1))\}\subseteq\mathbf{C}_{n+1}$,

$$r_n(x(-t,\phi_0),y(-t,\phi_1)) = t, v_{n+1}(x(-t,\phi_0),y(-t,\phi_1))$$
$$= -e^{-(\lambda+\mu)t}Jv_n(-t,\phi_0,\phi_1),$$

for $t\in(0,\hat{r}_n(\phi_0,\phi_1)\wedge\hat{r}(\phi_0,\phi_1))$. The third equation follows from the second and (A.14), and the second from the first and the fact $(x(t,x(-t,\phi_0)),y(t,y(-t,\phi_1)))=(\phi_0,\phi_1)\in\mathbf{\Gamma}$. Taking the limit in the third equation as $t$ increases to $\hat{r}_n(\phi_0,\phi_1)$ gives

$$\begin{cases} v_{n+1}(x(-t,\phi_0),y(-t,\phi_1))|_{t=\hat{r}_n(\phi_0,\phi_1)}=0 \text{ and} \\ (x(-\hat{r}_n(\phi_0,\phi_1),\phi_0),y(-\hat{r}_n(\phi_0,\phi_1),\phi_1))\in\partial\mathbf{\Gamma}^x_{n+1} \end{cases}$$
$$\text{if } \hat{r}_n(\phi_0,\phi_1)\leq\hat{r}(\phi_0,\phi_1).$$

Finally, every $(\tilde{\phi}_0,\tilde{\phi}_1)\in\mathbf{C}_{n+1}\cup\partial\mathbf{\Gamma}^x_{n+1}$ is reachable from $(\phi_0,\phi_1)\equiv r_n(\tilde{\phi}_0,\tilde{\phi}_1)\in\partial\mathbf{\Gamma}^e_{n+1}$ on the entrance boundary by $\{(x(t,\tilde{\phi}_0),y(t,\tilde{\phi}_1)); t\in[0,r_n(\tilde{\phi}_0,\tilde{\phi}_1)]\}$ which is contained (possibly, excluding end-points) in the continuation region $\mathbf{C}_{n+1}$. $\square$

PROOF OF COROLLARY 11.6. By Lemma 11.4, the function $(\phi_0,\phi_1)\mapsto r_0(\phi_0,\phi_1)$ is continuous on the continuation region $(\phi_0,\phi_1)\in\mathbf{C}_1$. Therefore, the entrance boundary $\partial\mathbf{\Gamma}^e_1$ is the image of the *continuous* mapping [see the definition in (9.10)]

$$(\phi_0,\phi_1)\mapsto(x(r_0(\phi_0,\phi_1),\phi_1),\gamma_1(y(r_0(\phi_0,\phi_1),\phi_1))), \qquad (\phi_0,\phi_1)\in\mathbf{C}_1,$$

from the *connected* $\mathbf{C}_1$ into $\mathbb{R}^2_+$. Thus the set $\partial\mathbf{\Gamma}^e_1$ is a connected subset of $\mathbb{R}^2_+$.

Since the parametric curves $t \mapsto (x(t,\phi_0), y(t,0))$, $\phi_0 \in \mathbb{R}_+$ starting on the $x$-axis are increasing, the points on the boundary $\partial \mathbf{\Gamma}_1$ where these curves meet the boundary belong to the entrance boundary $\partial \mathbf{\Gamma}_1^e$; see also Figure 1. Hence $\{(x, \gamma_1(x)) : x \in [\delta, \xi_1)\} \subseteq \partial \mathbf{\Gamma}_1^e$ for some $0 \leq \delta < \xi_1$. Then the connectedness of $\partial \mathbf{\Gamma}_1^e$ gives (11.14) with $\xi_1^e \triangleq \inf\{x \in \mathbb{R}_+ : (x, \gamma_1(x)) \in \partial \mathbf{\Gamma}_1^e\}$.

Indeed, the point $(\xi_1^e, \gamma_1(\xi_1^e))$ does not belong to the entrance boundary $\partial \mathbf{\Gamma}_1^e$. Suppose it does. Then $\{(x(-t, \xi_1^e), y(-t, \gamma_1(\xi_1^e))); t \in (0, \delta]\} \subset \mathbf{C}_1$ for some $\delta > 0$. Let $(\phi_0, \phi_1) \in \mathbf{C}_1$ be the point in the middle of the vertical line-segment connecting the points $(x(-\delta, \xi_1^e), y(-\delta, \gamma_1(\xi_1^e)))$ and $(x(-\delta, \xi_1^e), \gamma_1(x(-\delta, \xi_1^e)))$. Then we have $x(\delta, \phi_0) = x(\delta, x(-\delta, \xi_1^e)) = \xi_1^e$ and $y(\delta, \phi_1) > y(\delta, y(-\delta, \gamma_1(\xi_1^e))) = \gamma_1(\xi_1^e) = \gamma_1(x(\delta, \phi_0))$ since the mapping $\phi \mapsto y(t, \phi)$ is increasing for every $t \in \mathbb{R}$. Therefore, $(x(\delta, \phi_0), y(\delta, \phi_1)) \in \mathbf{\Gamma}_1$ and $0 < r_0(\phi_0, \phi_1) < \delta$. Thus we have $(x(r_0(\phi_0, \phi_1), \phi_0), y(r_0(\phi_0, \phi_1), \phi_1)) \in \partial \mathbf{\Gamma}_1^e$, but $x(r_0(\phi_0, \phi_1), \phi_0) < x(\delta, \phi_0) = \xi_1^e$ [the mapping $\phi \mapsto x(t, \phi)$ is increasing for every $t \in \mathbb{R}$]. This contradicts the minimality of $\xi_1^e$. □

PROOF OF COROLLARY 11.12. Suppose $\xi_1^e > 0$ and fix any $\phi_0 \in [0, \xi_1^e)$. Let $\bar{\phi}_0 \triangleq (1/2)(\phi_0 + \xi_1^e)$. Then $(\bar{\phi}_0, \gamma_1(\bar{\phi}_0)) \in \partial \mathbf{\Gamma}_1^x$, and $\gamma_1(\phi_0) > \gamma_1(\bar{\phi}_0) > \gamma_1(\xi_1^e)$ since $\gamma_1(\cdot)$ is strictly decreasing on its support. Then the set $B \triangleq [0, \bar{\phi}_0) \times (\gamma_1(\bar{\phi}_0), \infty)$ is an open neighborhood of $(\bar{\phi}_0, \gamma_1(\bar{\phi}_0))$ such that for every $(\tilde{\phi}_0, \tilde{\phi}_1) \in B \cap \mathbf{C}_1$, we have

$$0 < \underline{r} \leq r_0(\tilde{\phi}_0, \tilde{\phi}_1) \leq \bar{r} < \infty,$$

where $\underline{r} \triangleq \inf\{t \geq 0 : y(t, \gamma_1(\bar{\phi}_0)) \leq \gamma_1(\xi_1^e)\}$ and $\bar{r} \triangleq \inf\{t \geq 0 : x(t, 0) \geq \xi_1\}$. This completes the proof of the first part.

Now let $(\phi_0, \phi_1) \in \partial \mathbf{\Gamma}_1^e$ be a point on the entrance boundary. Take any convergent sequence $\{(\phi_0^{(n)}, \phi_1^{(n)})\}_{n \in \mathbb{N}}$ in the continuation region $\mathbf{C}_1$ whose limit is the boundary point $(\phi_0, \phi_1)$. Since $r_0(\cdot, \cdot) \leq \bar{r}$ (see above) on $\mathbf{C}_1$, the sequence $\{r_0(\phi_0^{(n)}, \phi_1^{(n)})\}_{n \in \mathbb{N}}$ is bounded and has a convergent subsequence. We shall conclude the proof of the second part by showing that every convergent subsequence of the sequence $\{r_0(\phi_0^{(n)}, \phi_1^{(n)})\}_{n \in \mathbb{N}}$ has the same limit 0.

Without changing the notation, suppose that $\{r_0(\phi_0^{(n)}, \phi_1^{(n)})\}_{n \in \mathbb{N}}$ converges to some finite number $r_0 \geq 0$. Since the functions $(\phi_0, \phi_1) \mapsto v_1(\phi_0, \phi_1)$ and $(t, \phi_0, \phi_1) \mapsto Jv_0(t, \phi_0, \phi_1)$ are continuous, we have

$$0 = v_1(\phi_0, \phi_1) = \lim_{n \to \infty} v_1(\phi_0^{(n)}, \phi_1^{(n)}) = \lim_{n \to \infty} Jv_0(r_0(\phi_0^{(n)}, \phi_1^{(n)}), \phi_0^{(n)}, \phi_1^{(n)})$$

$$= Jv_0(r_0, \phi_0, \phi_1) = \int_0^{r_0} e^{-(\lambda+\mu)t} G_0(t, \phi_0, \phi_1) \, dt.$$

If we show that $G_0(t, \phi_0, \phi_1) > 0$ for every $t > 0$, then $r_0 = 0$ follows.

However, $t = 0$ is a point of increase for the function $t \mapsto G_0(t, \phi_0, \phi_1)$. Since $(\phi_0, \phi_1) \in \partial \mathbf{\Gamma}_1^e = \{(x, a_0(x)) : x \in (\xi_1^e, \xi_1)\}$ by Corollary 11.7, and the boundary function $a_0(\cdot)$ of the region $A_0 = \{(x,y) \in \mathbb{R}_+^2 : [g + \mu \cdot v_0 \circ S](x,y) < 0\}$ is strictly decreasing, there exists some $\delta > 0$ such that $(x(t,\phi_0), y(t,\phi_1)) \in A_0 \subseteq \mathbf{C}_1$ for every $t \in [-\delta, 0)$. Therefore,

$$G(t, \phi_0, \phi_1) = [g + \mu \cdot v_0 \circ S](x(t,\phi_0), y(t,\phi_1)) < 0 = G_0(0, \phi_0, \phi_1),$$
$$t \in [-\delta, 0).$$

Then Lemma 11.3 implies that $G_0(t, \phi_0, \phi_1) > 0$ for every $t > 0$ and completes the proof of $r_0 = 0$. □

PROOF OF LEMMA 11.13. If $\xi_1^e = 0$, then $\mathrm{cl}(\partial \mathbf{\Gamma}_1^e) = \{(x, \gamma_1(x)) : x \in [0, \xi_1]\} = \partial \mathbf{\Gamma}_1$ by Corollary 11.6. In the remainder, suppose that $\xi_1^e > 0$ and fix any $\phi_0 \in [0, \xi_1^e)$. The boundary point $(\phi_0, \gamma_1(\phi_0))$ is not included in the entrance boundary $\partial \mathbf{\Gamma}_1^e$. We shall prove that it is an exit boundary point; namely, there exists some $\delta > 0$ such that [see (9.10)]

(A.15) $\qquad (x(t, \phi_0), y(t, \gamma_1(\phi_0))) \in \mathbf{C}_1 \qquad \forall t \in (0, \delta].$

Since the boundary $\gamma_1(\cdot)$ is strictly decreasing on its support $[0, \xi_1]$, we have

$$0 \le \phi_0 < \xi_1^e \quad \Longrightarrow \quad \gamma_1(\phi_0) > \gamma_1(\xi_1^e).$$

Then there is always a sequence of points $\{(\phi_0^{(n)}, \phi_1^{(n)})\}_{n \in \mathbb{N}} \subseteq \mathbf{C}_1$ such that

$$\phi_0^{(n)} = \phi_0 \quad \text{and} \quad \phi_1^{(n)} > \gamma_1(\xi_1^e) \qquad \text{for every } n \in \mathbb{N}, \text{ and}$$
$$\lim_{n \to \infty} \phi_1^{(n)} = \uparrow \gamma_1(\phi_0).$$

Namely, the sequence $\{(\phi_0^{(n)}, \phi_1^{(n)})\}_{n \in \mathbb{N}}$ "increases" to the point $(\phi_0, \gamma_1(\phi_0))$ along the vertical line passing through the point $(\phi_0, \gamma_1(\phi_0))$. For every $n \in \mathbb{N}$, we have

$$v_1(\phi_0^{(n)}, \phi_1^{(n)}) = J v_0(r_0(\phi_0^{(n)}, \phi_1^{(n)}), \phi_0^{(n)}, \phi_1^{(n)})$$

and

$$(x(r_0(\phi_0^{(n)}, \phi_1^{(n)}), \phi_0^{(n)}), y(r_0(\phi_0^{(n)}, \phi_1^{(n)}), \phi_1^{(n)})) \in \partial \mathbf{\Gamma}_1^e.$$

By Corollary 11.12, the sequence $\{r_0(\phi_0^{(n)}, \phi_1^{(n)})\}_{n \in \mathbb{N}}$ is bounded. Therefore, it has a convergent subsequence; we shall denote it by the same notation and its limit by $r_0$. The functions $J v_0(\cdot, \cdot, \cdot)$, $x(\cdot, \cdot)$, $y(\cdot, \cdot)$ and $v_1(\cdot, \cdot)$ are continuous, and $v_1(\phi_0, \gamma_1(\phi_0)) = 0$. Therefore, taking limits of the displayed equations above gives

(A.16) $0 = J v_0(r_0, \phi_0, \gamma_1(\phi_0)) \quad \text{and} \quad (x(r_0, \phi_0), y(r_0, \gamma_1(\phi_0))) \in \mathrm{cl}(\partial \mathbf{\Gamma}_1^e).$

The second expression implies that $x(r_0, \phi_0) \geq \xi_1^e$. We shall prove that the inequality is strict, and therefore,

(A.17) $$(x(r_0, \phi_0), y(r_0, \gamma_1(\phi_0))) \in \partial \mathbf{\Gamma}_1^e.$$

Let us assume that $x(r_0, \phi_0) = \xi_1^e$. Then the second expression in (A.16) implies that $(x(r_0, \phi_0), y(r_0, \gamma_1(\phi_0))) = (\xi_1^e, \gamma_1(\xi_1^e))$. Thus $(\phi_0, \gamma_1(\phi_0))$ is on the curve $\mathcal{C}_1$ given by (11.16). Then Corollary 11.9 implies that $G(t, \phi_0, \gamma_1(\phi_0)) > 0$ for every $t \neq r_0$. Since $r_0 > 0$, this implies that $Jv_0(r_0, \phi_0, \gamma_1(\phi_0)) = \int_0^{r_0} e^{-(\lambda+\mu)s} G_0(s, \phi_0, \gamma_1(\phi_0)) \, ds$ is strictly positive. But this contradicts the first equality in (A.16). Therefore, we must have $x(r_0, \phi_0) > \xi_1^e$, and (A.17) is correct.

Now we are ready to prove (A.15). Since $\phi_0 < \xi_1^e$, we have $[g + \mu \cdot v_0 \circ S](\phi_0, \gamma_1(\phi_0)) > 0$ by Corollary 11.11. Because the mapping $[g + \mu \cdot v_0 \circ S](\cdot, \cdot)$ is continuous, there exists some $r_0 > \delta > 0$ such that

$$G_0(t, \phi_0, \gamma_1(\phi_0)) = [g + \mu \cdot v_0 \circ S](x(t, \phi_0), y(t, \gamma_1(\phi_0))) > 0, \qquad t \in [0, \delta].$$

Then for every $t \in (0, \delta]$ we have

$$\begin{aligned} &v_0(x(t, \phi_0), y(t, \gamma_1(\phi_0))) \\ &\leq Jv_0(r_0 - t, x(t, \phi_0), y(t, \gamma_1(\phi_0))) \\ &= \int_0^{r_0 - t} e^{-(\lambda+\mu)u}[g + \mu \cdot v_0 \circ S](x(u, x(t, \phi_0)), y(u, y(t, \gamma_1(\phi_0)))) \, du \\ &= e^{(\lambda+\mu)t} \int_t^{r_0} e^{-(\lambda+\mu)u}[g + \mu \cdot v_0 \circ S](x(t, \phi_0), y(t, \gamma_1(\phi_0))) \, du \\ &= e^{(\lambda+\mu)t} \left[ \underbrace{Jv_0(r_0, \phi_0, \gamma_1(\phi_0))}_{=0} - \int_0^t e^{-(\lambda+\mu)u} G_0(u, \phi_0, \gamma_1(\phi_0)) \, du \right] < 0. \end{aligned}$$

Therefore, (A.15) holds and $(\phi_0, \gamma_1(\phi_0)) \in \partial \mathbf{\Gamma}_1^x$. □

PROOF OF LEMMA 11.14. There is nothing to prove if $\xi_1^e = 0$. Therefore, suppose $\xi_1^e > 0$. Let $B_1$ be the union of the continuation region $\mathbf{C}_1$ and the open subset of $[0, \xi_1^e) \times \mathbb{R}_+$ strictly below the curve $\mathcal{C}_1$ in Corollary 11.9. Then $B_1$ is open and $\mathbf{C}_1 \cup \partial \mathbf{\Gamma}_1^x \subset B_1$; see Remark 11.10. Define

$$\left\{ \begin{aligned} &\tilde{r}_0(\phi_0, \phi_1) \triangleq \inf\{t > 0 : (x(t, \phi_0), y(t, \phi_1)) \in \partial \mathbf{\Gamma}_1^e\} \\ &\tilde{v}_1(\phi_0, \phi_1) \triangleq Jv_0(\tilde{r}_0(\phi_0, \phi_1), \phi_0, \phi_1) \end{aligned} \right\}$$

for every $(\phi_0, \phi_1) \in B_1$.

Then

(A.18) $$r_0(\phi_0, \phi_1) = \tilde{r}_0(\phi_0, \phi_1) \quad \text{and} \quad v_1(\phi_0, \phi_1) = \tilde{v}_1(\phi_0, \phi_1),$$
$$(\phi_0, \phi_1) \in \mathbf{C}_1 \cup \partial \mathbf{\Gamma}_1^x.$$

Let us show that $\tilde{r}_0(\cdot,\cdot)$, and therefore, $\tilde{v}_1(\cdot,\cdot)$ are continuously differentiable on $B_1$. The infimum $\tilde{r}_0(\phi_0,\phi_1)$ is finite and strictly positive for every $(\phi_0,\phi_1) \in B_1$. By (9.12), $G_0(\tilde{r}_0(\phi_0,\phi_1),\phi_0,\phi_1)$ equals

$$[g + \mu \cdot v_0 \circ S](x(\tilde{r}_0(\phi_0,\phi_1),\phi_0), y(\tilde{r}_0(\phi_0,\phi_1),\phi_1)) = 0,$$
(A.19)
$$(\phi_0,\phi_1) \in B_1.$$

The mapping $(t,\phi_0,\phi_1) \mapsto G_0(t,\phi_0,\phi_1)$ is continuously differentiable. If

(A.20) $\quad D_t G_0(t,\phi_0,\phi_1)|_{t=\tilde{r}_0(\phi_0,\phi_1)} \neq 0, \qquad (\phi_0,\phi_1) \in B_1,$

then Theorem 11.2 implies that, in an open neighborhood in $B_1$ of every $(\phi_0,\phi_1)$, the equation $G_0(t,\phi_0,\phi_1) = 0$ determines $t = t(\phi_0,\phi_1)$ implicitly as a function of $(\phi_0,\phi_1)$, and this function is continuously differentiable. In every neighborhood, these solutions must then coincide with $\tilde{r}_0(\phi_0,\phi_1)$. Therefore, $\tilde{r}_0(\phi_0,\phi_1)$ is continuously differentiable on $B_1$. Then the function $\tilde{v}_1(\phi_0,\phi_1)$ is continuously differentiable on $B_1$ since $Jv_0(\cdot,\cdot,\cdot)$ is continuously differentiable on $\mathbb{R}^3_+$.

Now fix any $(\phi_0,\phi_1) \in B_1$ and assume $D_t G_0(\tilde{r}_0(\phi_0,\phi_1),\phi_0,\phi_1) = 0$. Then the function $t \mapsto G_0(t,\phi_0,\phi_1)$ has a local minimum at $t = \tilde{r}_0(\phi_0,\phi_1)$. Lemma 11.3 and (A.19) imply that $G_0(t,\phi_0,\phi_1) > 0$ for every $t \neq \tilde{r}_0(\phi_0,\phi_1)$. Therefore, the parametric curve

$$\{(x(t,\phi_0), y(t,\phi_1)) : t \in \mathbb{R}\} \cap \mathbb{R}^2_+ \subseteq \mathbb{R}^2_+ \setminus A_0$$

does not intersect $A_0$, but touches the boundary $\partial A_0$. Then this curve has to be the same as $\mathcal{C}_1$ in Corollary 11.16, and $(\phi_0,\phi_1) \in \mathcal{C}_1$. But this contradicts $(\phi_0,\phi_1) \in B_1$, since Remark 11.10 and the description of $B_1$ show that $\mathcal{C}_1 \cap B_1 = \varnothing$. Therefore, (A.20) holds.

Now let us show that $\gamma_1(\cdot)$ is continuously differentiable on $[0,\xi_1^e)$. Fix any $\phi_0 \in [0,\xi_1^e)$. Then $(\phi_0,\gamma_1(\phi_0)) \in \partial \mathbf{\Gamma}_1^x \subset B_1$ and $\tilde{v}_1(\phi_0,\gamma_1(\phi_0)) = 0$ by (A.18). The function $\tilde{v}_1(\cdot,\cdot)$ is continuously differentiable on $B_1$. Therefore, the result will again follow from the implicit function theorem (Theorem 11.2) if we show that $D_{\phi_1} \tilde{v}_1(\phi_0,\gamma_1(\phi_0)) \neq 0$. However, $D_{\phi_1} \tilde{v}_1(\phi_0,\gamma_1(\phi_0))$ equals

$$D^-_{\phi_1} \tilde{v}_1(\phi_0,\gamma_1(\phi_0)) = D^-_{\phi_1} v_1(\phi_0,\gamma_1(\phi_0)) = \lim_{\phi_1 \uparrow \gamma_1(\phi_0)} D^-_{\phi_1} v_1(\phi_0,\phi_1)$$

$$= \lim_{\phi_1 \uparrow \gamma_1(\phi_0)} D_{\phi_1} v_1(\phi_0,\phi_1)$$

$$= \lim_{\phi_1 \uparrow \gamma_1(\phi_0)} \frac{1 - e^{-(\mu+1)r_0(\phi_0,\phi_1)}}{\mu + 1} > 0.$$

The second equality follows from (A.18), and the third from the concavity of $v_1(\cdot,\cdot)$. The fourth and the fifth follow from Corollary 11.5. Finally, the limit at the end is strictly positive since $r_0(\cdot,\cdot)$ is bounded away from zero in the intersection of $\mathbf{C}_1$ with some neighborhood of $(\phi_0, \gamma_1(\phi_0))$ by Corollary 11.12. □

PROOF OF LEMMA 11.15. The result follows from Corollary 11.7 if $\xi_1^e = 0$. Therefore, suppose $\xi_1^e > 0$. Then the boundary function $\gamma_1(\cdot)$ is continuously differentiable on $[0, \xi_1^e) \cup (\xi_1^e, \xi_1)$ by Corollary 11.7 and Lemma 11.14. We need to show that $x \mapsto \gamma_1(x)$ is continuously differentiable at $x = \xi_1^e$.

Recall that the function $\gamma_1(\cdot)$ is convex. Therefore, the left derivative $D^-\gamma_1(\cdot)$ and the right derivative $D^+\gamma_1(\cdot)$ of the function $\gamma_1(\cdot)$ exist and are left- and right-continuous, respectively, at $x = \xi_1^e$. Thus

(A.21)
$$\lim_{x \uparrow \xi_1^e} D\gamma_1(x) = \lim_{x \uparrow \xi_1^e} D^-\gamma_1(x) = D^-\gamma_1(\xi_1^e)$$
$$\leq D^+\gamma_1(\xi_1^e) = \lim_{x \downarrow \xi_1^e} D^+\gamma_1(x) = \lim_{x \downarrow \xi_1^e} D\gamma_1(x).$$

The continuity of the derivative $D\gamma_1(\cdot)$ of the function $\gamma_1(\cdot)$ at $x = \xi_1^e$ will follow immediately from the existence of the derivative of $\gamma_1(\cdot)$ at $x = \xi_1^e$.

Now recall from Corollary 11.9 and Remark 11.10 that the point $(\xi_1^e, \gamma_1(\xi_1^e))$ is on the parametric curve $\mathcal{C}_1$, which lies above $\{(x, \gamma_1(x)) : x \in \mathbb{R}_+\}$ and touches it at the point $(\xi_1^e, \gamma_1(\xi_1^e))$. Therefore, for every $t > 0$ and $s > 0$

$$\frac{y(0, \gamma_1(\xi_1^e)) - y(-t, \gamma_1(\xi_1^e))}{x(0, \xi_1^e) - x(-t, \xi_1^e)} \leq \frac{\gamma_1(x(0, \xi_1^e)) - \gamma_1(x(-t, \xi_1^e))}{x(0, \xi_1^e) - x(-t, \xi_1^e)}$$
$$\leq \frac{\gamma_1(x(s, \xi_1^e)) - \gamma_1(x(0, \xi_1^e))}{x(s, \xi_1^e) - x(0, \xi_1^e)}$$
$$\leq \frac{y(s, \gamma_1(\xi_1^e)) - y(0, \gamma_1(\xi_1^e))}{x(s, \xi_1^e) - x(0, \xi_1^e)}.$$

When we take the limit as $t \downarrow 0$ and $s \downarrow 0$, we obtain

$$\frac{D_t y(0, \gamma_1(\xi_1^e))}{D_t x(0, \xi_1^e)} \leq D^-\gamma_1(\xi_1^e) \leq D^+\gamma_1(\xi_1^e) \leq \frac{D_t y(0, \gamma_1(\xi_1^e))}{D_t x(0, \xi_1^e)}.$$

Note that the terms on the far left and far right are the same. Therefore, $D^-\gamma_1(\xi_1^e) = D^+\gamma_1(\xi_1^e)$ and the derivative of the boundary function $\gamma_1(\cdot)$ at $x = \xi_1^e$ exists. □

PROOF OF LEMMA 11.16. Since $v_1(\cdot, \cdot)$ is concave, the left derivatives $D^-_{\phi_0} v_1(\cdot, \cdot)$, $D^-_{\phi_1} v_1(\cdot, \cdot)$ and the right derivatives $D^+_{\phi_0} v_1(\cdot, \cdot)$, $D^+_{\phi_1} v_1(\cdot, \cdot)$ exist and are left- and right-continuous on the boundary $\partial \mathbf{\Gamma}$, respectively. Because $v_1(\cdot, \cdot)$ vanishes on $\mathbf{\Gamma}_1$, and the function $\gamma_1(\cdot)$ is strictly decreasing, we have

(A.22)
$$D^-_{\phi_0} v_1(\phi_0, \phi_1) \geq D^+_{\phi_0} v_1(\phi_0, \phi_1) = 0,$$

$$(\phi_0, \phi_1) \in \partial \mathbf{\Gamma}_1 \setminus \{(0, \gamma_1(0))\}.$$

(A.23)
$$D^-_{\phi_1} v_1(\phi_0, \phi_1) \geq D^+_{\phi_1} v_1(\phi_0, \phi_1) = 0,$$

$$(\phi_0, \phi_1) \in \partial \mathbf{\Gamma}_1 \setminus \{(\xi_1, 0)\}.$$

For every boundary point $(\phi_0, \phi_1) \in \partial \mathbf{\Gamma}_1 \setminus \{(0, \gamma_1(0))\}$ and any $\{(\phi_0^{(n)}, \phi_1^{(n)})\}_{n \in \mathbb{N}} \subset \mathbf{C}_1$ such that $\lim_{n \to \infty} \phi_0^{(n)} = \uparrow \phi_0$ and $\phi_1^{(n)} = \phi_1$ for every $n \in \mathbb{N}$, we have

(A.24)
$$\begin{aligned} D^-_{\phi_0} v_1(\phi_0, \phi_1) &= \lim_{n \to \infty} D^-_{\phi_0} v_1(\phi_0^{(n)}, \phi_1^{(n)}) \\ &= \lim_{n \to \infty} D_{\phi_0} v_1(\phi_0^{(n)}, \phi_1^{(n)}) \\ &= \lim_{n \to \infty} \frac{1 - \exp\{-(\mu - 1) r_0(\phi_0^{(n)}, \phi_1^{(n)})\}}{\mu - 1}. \end{aligned}$$

The second and the third equalities follow from Corollary 11.5. The function $r_0(\cdot, \cdot)$ is continuous on the entrance boundary $\partial \mathbf{\Gamma}_1^e$ and is bounded away from zero in some neighborhood of every point on the exit boundary $\partial \mathbf{\Gamma}_1^x$; see Corollary 11.12. Therefore, the limit on the right in (A.24) equals zero for every point $(\phi_0, \phi_1)$ on the entrance boundary $\partial \mathbf{\Gamma}_1^e$ and is strictly positive for every point $(\phi_0, \phi_1)$ on the exit boundary $\partial \mathbf{\Gamma}_1^x$.

Thus, for every $(\phi_0, \phi_1) \in \partial \mathbf{\Gamma}_1^e$, the equality in (A.22), and as a result of a similar argument, the equality in (A.23) are attained. Therefore, the partial derivatives $D_{\phi_0} v_1(\cdot, \cdot)$ and $D_{\phi_1} v_1(\cdot, \cdot)$ exist at every $(\phi_0, \phi_1) \in \partial \mathbf{\Gamma}_1^e$ and are continuous since $D_{\phi_0} v_1(\cdot, \cdot) = D^\pm_{\phi_0} v_1(\cdot, \cdot)$ is both left- and right-continuous near the entrance boundary $\partial \mathbf{\Gamma}_1^e$.

However, if $(\phi_0, \phi_1)$ is a point on the exit boundary $\partial \mathbf{\Gamma}_1^x$, then the inequalities in (A.23) and (A.24) are strict. Namely, the $v_1(\cdot, \cdot)$ is not differentiable on the exit boundary $\partial \mathbf{\Gamma}_1^x$. $\square$

**Acknowledgment** We are grateful to the referee for detailed comments that helped us improve the manuscript.

E. BAYRAKTAR
DEPARTMENT OF MATHEMATICS
UNIVERSITY OF MICHIGAN
ANN ARBOR, MICHIGAN 48109
USA
E-MAIL: erhan@umich.edu

S. DAYANIK
DEPARTMENT OF OPERATIONS RESEARCH
  AND FINANCIAL ENGINEERING
AND THE BENDHEIM CENTER FOR FINANCE
PRINCETON UNIVERSITY
PRINCETON, NEW JERSEY 08544
USA
E-MAIL: sdayanik@princeton.edu

I. KARATZAS
DEPARTMENTS OF MATHEMATICS
  AND STATISTICS
COLUMBIA UNIVERSITY
NEW YORK, NEW YORK 10027
USA
E-MAIL: ik@math.columbia.edu